\documentclass[12pt,reqno]{amsart}

\usepackage[LGR,T1]{fontenc}
\usepackage[utf8]{inputenc}
\usepackage[foot]{amsaddr}

\title[]{Monge-Ampère gravitating fluids:\\ Least action principles and particle systems}
\date{February, 2025. Second revised version:  May, 2026}

\author[]{Christian Léonard $^{1,2}$}
\address{$^1$MODAL'X, Univ.\ Paris Nanterre, CNRS, F92000 Nanterre France}
\address{$^2$INRIA, Paris, France}

\author[]{Roya Mohayaee $^{3,4}$}
\address{$^3$Sorbonne Université, CNRS, Institut d’Astrophysique de Paris, France}
\address{$^4$Rudolf Peierls Centre for Theoretical Physics, Univ.\ of Oxford, United Kingdom}
\email{christian.leonard@math.cnrs.fr, mohayaee@iap.fr}

\usepackage{amssymb, ams math, amsfonts, latexsym, enumerate}
\usepackage[usenames,dvipsnames]{pstricks}
\usepackage{epsfig}
\RequirePackage[colorlinks,linkcolor=blue,citecolor=blue,urlcolor=blue]{hyperref}
\usepackage[english]{babel} 
\usepackage{graphicx,caption}
\usepackage{subcaption}

\usepackage{tikz}
\usetikzlibrary{arrows.meta, positioning, shadows}

\newtheorem{theorem}[equation]{Theorem}
\newtheorem{lemma}[equation]{Lemma}
\newtheorem{proposition}[equation]{Proposition}
\newtheorem{corollary}[equation]{Corollary}

\newtheorem{statement}[equation]{Statement}
\newtheorem{conjecture}[equation]{Conjecture}

\newtheorem{definition}[equation]{Definition}
\newtheorem{definitions}[equation]{Definitions}

\newtheorem{properties}[equation]{Properties}
\newtheorem{paraset}[equation]{Parameter setting}

\newtheorem{phyp}[equation]{Physics hypothesis}

\theoremstyle{remark}
\newtheorem{remark}[equation]{Remark}
\newtheorem{remarks}[equation]{Remarks}

\numberwithin{equation}{section}

\newcommand{\RR}{\mathbb{R}}
\newcommand{\Rn}{\mathbb{R}^n}  \newcommand{\Rd}{\mathbb{R}^d}
\newcommand{\PP}{\mathbb{P}}
\newcommand{\EE}{\mathbb{E}}
\newcommand{\1}{\mathbf{1}}

\newcommand\pf{_{\#}}

\newcommand{\Leb}{\mathrm{Leb}}

\newcommand{\as}{\textrm{-}\mathrm{a.s.}}
\renewcommand{\ae}{\textrm{-}\mathrm{a.e.}}
\newcommand{\Id}{\mathrm{Id}}
\newcommand{\scal}{\!\cdot\!}

\newcommand{\ddt}{ \frac{d}{dt}}

\newcommand{\mm}{ \mathsf{m}}

\newcommand{\pp}{ \mathsf{p}}
\newcommand{\qq}{ \mathsf{q}}
\newcommand{\rr}{ \mathsf{r}}
\newcommand{\vv}{ \mathsf{v}}

\DeclareMathOperator{\cl}{cl}
\DeclareMathOperator{\cv}{cv}

\DeclareMathOperator{\supp}{supp}

\DeclareMathOperator{\argmin}{argmin}

\DeclareMathOperator{\proj}{proj}

\DeclareMathOperator{\grad}{grad}
\DeclareMathOperator{\Hess}{Hess}

\newcommand{\sbt}{\,\begin{picture}(-1,1)(-1,-3)\circle*{3}\end{picture}\ } 

\newcommand\Glim[1]{\Gamma\textrm{-}\lim_{#1\rightarrow\infty}}

\newcommand\lime{\lim_{\epsilon\rightarrow0}}

\newcommand{\ud}{\frac{1}{2}}
\newcommand{\uud}{\scalebox{1.2}{$\ud$}}

\newcommand{\MAG}{{\small MAG}}

\newcommand{\ZZ}{\Rd}
\newcommand\ZZZ{\ZZ\times\ZZ}

\newcommand\OO{\Omega}

\newcommand\PZ{\mathrm{P}(\ZZ)}

\newcommand\PO{\mathrm{P}(\OO)}

\newcommand\IZ{\int_{\ZZ}}

\newcommand{\xx}{ \mathsf{x}}
\newcommand{\yy}{ \mathsf{y}}
\newcommand{\zz}{ \mathsf{z}}

\newcommand{\Tf}{{\overrightarrow{T}}}
\newcommand{\Tb}{\overleftarrow{T}}
\newcommand{\Proba}{\mathrm{Proba}}
\DeclareMathOperator{\Proj}{Proj}
\newcommand{\IS}{\int _{ s_0} ^{ s_1}}
\newcommand{\ITT}{\int _{ - \infty} ^{ t_1}}
\newcommand{\IT}{\int _{ t_0} ^{ t_1}}

\newcommand{\IRT}{\int _{\Rd\times [t_0,t_1]}}
\newcommand{\uu}{ \mathsf{u}}
\newcommand{\Xh}{\widehat X^N\!}
\newcommand{\lk}{\lambda ^{(k)}}
\newcommand{\vce}{{\dot\rr^\epsilon}}

\newcommand{\vvf}[1]{\overrightarrow{\vv _{ #1}}}
\newcommand{\vvc}[1]{\vv ^{ \mathrm{cu},#1}}

\newcommand{\Rdk}{ (\ZZ)^k}
\newcommand{\Xe}{ \mathsf{X}^{\epsilon}}

\newcommand{\XtN}{ \widetilde{X}^{N}}
\newcommand{\grOW}[1]{\mathrm{grad} ^{ \mathrm{OW}} _{ #1}}

\newcommand{\Ok}{\OO ^{ (k)}}
\newcommand{\POk}{ \mathrm{P}(\Ok)}
\newcommand{\OkP}{\Ok _{ \mathrm{P}}}
\newcommand{\OP}{\OO _{ \mathrm{P}}}
\newcommand{\PRdk}{ \mathrm{P}(\Rdk)}

\newcommand{\Yb}{\overline{Y}}

\newcommand{\Zb}{\overline{Z}}
\newcommand{\Zt}{\widetilde{Z}}
\newcommand{\Xb}{\overline{X}}

\newcommand{\gOW}{\grad^{\mathrm{OW}}}

\newcommand{\EUR}{{\small EUR}}
\newcommand{\OW}{{\small OW}}

\newcommand{\txa}{\tilde{ \mathsf{x}} _{ \gamma}}
\newcommand{\Rem}{R ^{ m, \epsilon}}
\renewcommand{\Re}{R ^{m,V, \epsilon}}

\newcommand{\rb}{\overline{\rho}}
\newcommand{\gcb}{\overline{\nabla \mathcal{C}}}

\keywords{Optimal transport, Otto-Wasserstein manifold, large deviations, Schrödinger problem, early Universe reconstruction}
\subjclass[2010]{49Q22, 60F10, 70F45, 85A40}
\thanks{This research has been conducted within the FP2M federation (CNRS FR 2036)}

\begin{document}

\begin{abstract} 
The Monge-Ampère gravitation theory  ({\tiny MAG}) was introduced by Brenier \cite{Bre11}  to obtain an approximate solution of the early Universe reconstruction problem. It is a modification of  Newtonian gravitation which is based on quadratic optimal transport. Later, Brenier \cite{Bre16}, then  Ambrosio, Baradat and Brenier \cite{ABB20} discovered a double large deviation principle for Brownian particles whose rate function is precisely {\tiny MAG}'s action functional.

In the present article, following Brenier we first recap {\tiny MAG}'s theory. Then,  we slightly  extend it from particles to fluid. This allows us  to revisit the  Ambrosio-Baradat-Brenier particle system.

We propose another particle system which is easier to  interpret in physics  and whose large deviation rate function is half the way to  {\tiny MAG}'s action functional for fluids. 
While the setting of the Schrödinger problem is a system of noninteracting particles, our  particle system   is subject to  some splitting mechanism which regulates the thermal  fluctuations. This gives rise to some conditional Gibbs principle that leaves us with an action functional on the fluid space which is  {\tiny MAG}'s action functional plus an extra term associated with thermal fluctuations.   

In order to recover  {\tiny MAG}'s action functional, we have to remove this extra term. To do so, we propose to add some quantum force field on the Otto-Wasserstein manifold of fluids to  balance the thermal fluctuations.  A microscopic description of a system of particles leading to a conditional Gibbs principle whose action functional generates  such a quantum force remains a challenging open problem.
\end{abstract}

\maketitle 
\noindent\emph{This article is dedicated to  my long-time friend Patrick Cattiaux, on the occasion of his retirement. CL.}

\tableofcontents

\section{Introduction}

The Monge-Ampère gravitation (\MAG) theory was introduced by Brenier in \cite{Bre11} to obtain an approximate solution of the early Universe reconstruction (\EUR) problem. It is a modification of  Newtonian gravitation which is based on quadratic optimal transport. Its action functional is, for any path $ \omega=( \omega_s) _{ s_0\le s\le s_1}$ in $(\Rd)^k$,  
\begin{align*}
 A( \omega):=
\IS \uud \| \dot\omega_s-\vv_s( \omega_s)\|^2 k ^{ -1} \kappa_s\, ds,
\end{align*}
where  $\vv_s(x)$ is a  specific velocity vector field which encodes some optimal transport feature,   $\kappa_s$ is a specific  positive number   whose inverse $\kappa ^{ -1}_s$ is a  diffusion coefficient  and  as usual $\dot \omega_s:= d \omega_s/ds$ stands for the velocity of $ \omega$.  
This means that the system evolves in $(\ZZ)^k$ and its  motion $ \omega$ solves the least action problem
\begin{align}\label{eq-zz}
\inf _{ \omega: \omega(s_0)= \mathsf{a}, \omega(s_1)= \mathsf{b}} A (\omega)
\end{align}
where $ \mathsf{a}$ and $ \mathsf{b}$ are prescribed initial and final positions. See \eqref{eq-24} for the exact formulation of the action functional.

A short note about action functionals and their link with the equations of motion is proposed at Appendix \ref{app-b}.

\subsection*{Why mappings rather than clouds?}

  The coexistence of $k$ particles  is essential to convey (semi-discrete) optimal transport into the model. 
Any element of $(\ZZ)^k$ can be interpreted as a cloud of $k$ particles  living in the state  space $\Rd.$ But, one key idea of  \cite{ABB20} is to interpret the $j$-th coordinate of $(z_j) _{ 1\le j\le k}\in (\ZZ)^k$ as the position of one particle starting from   $x_j\in\ZZ$ where $x_1,\dots,x_k$ are prescribed initial positions. With this in mind, we see that $(z_j) _{ 1\le j\le k}$ encodes the mapping $T: \left\{ x_1,\dots, x_k\right\} \mapsto \ZZ$ defined by $T( x_j):= z_j,$ $1\le j\le k.$ Let us call  $(z_j) _{ 1\le j\le k}$ a $k$-mapping.  The reason for interpreting this cloud as some mapping will appear later in relation with optimal transport, see Definition \ref{def-km}.  Indeed, we shall be interested in the Monge-Kantorovich optimal transport problem from the initial probability measure $  k ^{ -1} \sum _{ 1\le j\le k} \delta _{  x_j}$ to some observed target measure. Here and for the rest of this article, $ \delta_a$ is the Dirac measure at $a$. 
As a consequence, we shall use the full information carried by the vector $(z_1,\dots,z_k)$ associated with a mapping, rather than the poorer information given by the set $ \left\{z_1,\dots,z_k\right\}$ attached to a cloud.

\subsection*{A double large deviation principle}

In two  articles by Brenier \cite{Bre16} and Ambrosio, Baradat and Brenier  \cite{ABB20}, \MAG's action functional was interpreted as  the rate function of some double large deviation principle involving Brownian trajectories. More precisely, it happens that 
\[A= \Gamma\textrm{-}\lim_{ \epsilon\rightarrow 0} A^ \epsilon,\] where
\begin{align*}
A^ \epsilon( \omega):=\IS  \uud \| \dot\omega_s-\vv^ \epsilon_s( \omega_s)\|^2 k ^{ -1}  \kappa_s\, ds,
\end{align*}
 and  $\vv^ \epsilon$ is the \emph{current velocity} of a rescaled Brownian motion in $(\Rd)^k$
\begin{align}\label{eq-23b}
 \Xe_s= \mathsf{X} _{ s_0}+ \sqrt{ \epsilon}\ \mathsf{B}_s,
 \quad s_0\le s\le s_1,
\end{align}
see \eqref{eq-30}. 
The forward velocity of $\Xe$ is zero but since the time marginal flow $(r^ \epsilon_s) _{ s_0\le s\le s_1}$ of $\Xe$ solves the heat equation, we have
 \begin{align}\label{eq-114}
  \partial_s r^ \epsilon-\epsilon \Delta r^ \epsilon/2
= \partial_s r^ \epsilon+\nabla\scal ( r^ \epsilon \vv^ \epsilon)=0,
 \end{align}
 with 
 \[\vv^ \epsilon_s=- \epsilon\nabla\log \sqrt{r^ \epsilon_s}.
 \]
This expression is that of a current velocity because it enters the continuity equation $\partial_s r^ \epsilon+\nabla\scal ( r^ \epsilon \vv^ \epsilon)=0$. 
 In \cite{ABB20},  the  stochastic differential equation  in $(\Rd)^k$,
\begin{align}\label{eq-121}
d \mathsf{Z} ^{ \epsilon,\eta}_s
	=\vv^ \epsilon_s( \mathsf{Z} ^{ \epsilon,\eta}_s)\,ds+ \sqrt{\eta k \kappa_s ^{ -1}}\, d \mathsf{W}_s ,\quad 0<s_0\le s\le s_1,
\end{align}
is introduced, where $\mathsf{W}$ is another $(\Rd)^k$-valued Brownian motion and the drift field is the current velocity $\vv^ \epsilon$  of $\Xe$, see \eqref{eq-21}.   It is known that this collection of random evolutions obeys  the Freidlin-Wentzell large deviation principle  when the parameter $\eta$ tends to zero, $ \epsilon$ being fixed,
\begin{align*}
\Proba( \mathsf{Z} ^{ \epsilon,\eta}\in \sbt)
 \underset{\eta\to 0}\asymp \exp \left( - \eta ^{ -1}\ \inf _{ \zz\in \bullet} A^ \epsilon( \zz)\right)
\end{align*}
 with rate function $A^ \epsilon$, \cite{FW,DZ}. 
Therefore, taking into account the already mentioned Gamma limit:  $ \Gamma\textrm{-}\lim_{ \epsilon\rightarrow 0} A^ \epsilon=A,$  we see that \emph{letting $\eta\to 0$, then $ \epsilon\to 0,$ the family of Brownian diffusion processes $ \mathsf{Z} ^{ \epsilon,\eta}$ obeys the  large deviation principle
\begin{align*}
\Proba( \mathsf{Z} ^{ \epsilon,\eta}\in \sbt)
 \underset{\eta\to 0, \epsilon\to 0}\asymp \exp \left( -  \eta ^{ -1}\ \inf _{ \zz\in \bullet} A( \zz)\right)
\end{align*}
 with \MAG's action functional $A$ as its  rate function.} This also implies the Gibbs conditioning principle
\begin{align*}
\Proba( \mathsf{Z} ^{ \epsilon,\eta}\in \sbt\mid  \mathsf{Z} ^{ \epsilon,\eta} _{ s_0}= \mathsf{a},  \mathsf{Z} ^{ \epsilon,\eta} _{ s_1}= \mathsf{b}) \underset{\eta\to 0, \epsilon\to 0}\asymp \exp \left( -  \eta ^{ -1}\ \inf _{ \omega\in \bullet, \omega(s_0)= \mathsf{a}, \omega(s_1)= \mathsf{b}} A( \omega)\right),
\end{align*}
which in turns means that, conditionally on $  \omega(s_0)= \mathsf{a}$ and $  \omega(s_1)= \mathsf{b},$ the most likely path of $\mathsf{Z} ^{ \epsilon,\eta}$ as $\eta\to0$, then $ \epsilon\to0,$ solves \MAG's least action problem \eqref{eq-zz}.

The stochastic evolution of $ \mathsf{Z} ^{ \epsilon,\eta}$ is interpreted in \cite{ABB20} as a cloud of $k$ Brownian particles \emph{surfing the heat wave}.

\subsection*{Aim of the  article}

The main goal of the present article is to revisit this interpretation. Indeed,  the above  model does not provide a clear physical picture.  In particular,  the \emph{forward} velocity vector field $\vv^ \epsilon$ of $ \mathsf{Z} ^{ \epsilon,\eta}$ happens to be the \emph{current} velocity vector field of the \emph{other} diffusion process: $\Xe.$ This enigmatic substitution  is uneasy to interpret.
\\
We propose another particle system inspired by \cite{ABB20} whose limit when the number of  particles tends to infinity describes a fluid keeping \MAG's main properties.

\subsection*{A heuristic and explorative point of view}

This article is not a typical mathematical paper, nor a physics paper. It is not even an hybrid of these two types, since it is not rigorous enough as regards the usual standard in mathematics, but it contains  too much  mathematics and not enough physics to be considered as a physics paper.

It presents a model   whose certain building components require mathematical tools (optimal transport and  Otto-Wasserstein geometry) that are, for the moment,  too specialized for most physicists\footnote{But we think that optimal transport and \OW-geometry are  efficient  tools in physics, and hope that this article  will help promoting them among some physicists.}. On the other hand, it is not clear at present time that this model is physically valid. Maybe some parts of it should be kept and some others should be rejected. We hope however that, in case the model presented in this article fails to be physically plausible, at least some   ideas  could be preserved and could  help building a more established theory. Note that the pioneering articles \cite{Bre11,Bre16,ABB20} as well as \cite{vR11} have, to some extent,  this same flavor. 

To our opinion, it is not necessary at this stage of research, where the physical pertinence of the model still needs to be challenged,  to prove the mathematical results with full rigor. However, we decided either to present  the entire chains of arguments during the ``proofs'', or to rely on Otto's heuristic calculus \cite{Vill09}. 
Our main lack of mathematical rigor consists in \emph{disregarding the regularity issues when they are not critical}. Any incomplete proof will start by \emph{"Proof"} rather than the usual \emph{Proof}.

\subsection*{Acknowledgements}

The authors  warmly thank an anonymous referee for many suggestions and questions that improved significantly the initial version of this article. They also thank the editorial board of  the \emph{Annales de la faculté des sciences de Toulouse} for having considered positively the publication of this long maths-physics heuristic paper in this special volume in honor of Patrick Cattiaux.

\section{Main results}

This section is dedicated to an informal presentation of the main results of the article. The outline of the paper is presented at page \pageref{outline}.

\subsection*{Another particle system}

In the present paper, we  drop the  particle system $ \mathsf{Z} ^{ \epsilon,\eta}$. But we  keep \cite{ABB20}'s key idea of working with the $k$-mapping valued Brownian process $\Xe$, because it establishes a  crucial link with quadratic optimal transport: a central feature of \MAG. We are going to investigate the large deviations of the empirical process   of a family of  splitting $(\ZZ)^k$-valued Brownian  particles, as their number tends to infinity. In absence of an extra force field, each sibling is a copy of $\Xe.$  But to arrive at \MAG, it is necessary that the whole system is   { immersed in some additional  quantum force field. }

Let us give some more details during this introductory section. 

\subsection*{From particles to fluids}

Passing to the limit as the number of particles tends to infinity, one does not work anymore with individual particles, but with fluids, i.e.\ probability measures on the state space. In fact, we shall look at probability measures on $(\Rd)^k,$ i.e.\ fluids of $k$-mappings.   We call these probability measures: $k$-fluids.  Keeping $k$-mappings  is essential to trace the effect of optimal transport. However, the relevant system to be observed is not the $k$-fluid, which contains the full details of the history of the $k$ particles that  are necessary to play with optimal transport, but its one-particle marginal measure on $\Rd.$

We do not perform the limit $ \epsilon\to 0.$ Neither do we look at the limit $k\to \infty,$ nor at the one-particle marginal projection of the $k$-fluid,  leaving these steps to forthcoming investigations.  We only look at the fluid analogue of the above  action functional $A^ \epsilon,$  which is
\begin{align}\label{eq-115}
(p_s) _{ s_0\le s\le s_1}\mapsto 
\IS \uud \|\dot\pp_s-\dot \rr^ \epsilon_s\|^2 _{ p_s}\, \kappa_s\, ds,
\end{align}
where $(p_s) _{ s_0\le s\le s_1}$ is a $k$-fluid-valued path, $\dot\pp$ is its velocity in the Otto-Wasserstein (\OW) manifold, $\|\sbt\|_p$ is the norm on the tangent vector space at $p$ (beware, it is not an $L^p$-norm) and 
\begin{align*}
\dot \rr^ \epsilon_s=\vv^ \epsilon_s=- \epsilon\nabla\log \sqrt{r^ \epsilon_s}
\end{align*}
is the current velocity of the heat flow $(r^ \epsilon_s) _{ s_0\le s\le s_1}$ -- recall \eqref{eq-114} -- in the \OW-manifold. 
We call the corresponding model: $ \epsilon$-\MAG\ for $k$-fluid.
\\
Basic definitions and results about the \OW-geometry are gathered at Appendix \ref{app-OW}.

\subsection*{A preliminary system of independent particles}

On the way to the action functional \eqref{eq-115}, we  first consider the empirical process 
\begin{align}\label{eq-x4}
X^N: s\in [s_0,s_1]\mapsto \frac 1N   \sum _{ 1\le i\le N} \delta _{ \mathsf{X} ^{ \epsilon}_i(s)}
\in \mathrm{P}\big(\Rdk \big)
\end{align}
of  a  sequence  $( \mathsf{X} ^{ \epsilon}_i) _{ i\ge 1}$ of independent copies  of $ \mathsf{X} ^{ \epsilon}.$ The evolution of the random   \emph{path of $k$-mappings} $\Xe$ is defined at \eqref{eq-23b} and its initial law  is such that for any $1\le j\le k,$ its $j$-th coordinate  $ \mathsf{X} ^{ \epsilon,j}$ is the path of  a particle starting from  the $j$-th random draw without replacement  from the set $\{ x_1,\dots,  x_k\}.$ Although   the coordinates $ \mathsf{X} ^{ \epsilon,j},$  $1\le j\le k,$ are correlated, they share the same law with initial distribution $ k ^{ -1}\sum _{ 1\le j\le k} \delta _{  x_j}.$

Relying on the Schrödinger problem approach, we shall establish at Section \ref{sec-nbm} the following already known  result.

\begin{statement}[\cite{Leo12e,CCGL20}]\emph{(Gibbs conditioning principle).} \label{res-x2} For any probability measures $ \alpha$ and $ \beta$ on $\Rdk,$ conditionally on  $ X^N(s_0)\simeq \alpha$ and $ X^N(s_1)\simeq \beta,$  the most likely trajectory $(p_s) _{ s_0\le s\le s_1}\in  C([s_0,s_1], \PRdk) $ of $X^N$ as $N$ tends to infinity solves the least action problem 
\begin{align}\label{eq-x6}
\inf _{ (p)} \IS \Big(\uud   \| \dot\pp_s-\vce_s\|^2_{p_s}  
	+ \epsilon^2 I(p_s| r^ \epsilon_s)\Big)\, ds
\end{align}
where the infimum runs through all the $(p):[s_0,s_1]\to \PRdk $ satisfying $p _{ s_0}=\alpha $ and $ p _{ s_1}=\beta ,$ and\begin{align*}
I(p|r):= \uud  \int \Big|\nabla\log \sqrt{ \frac{dp}{dr}}\Big|^2\, dp
\end{align*}
 is the Fisher information of $p$ with respect to $r$.
 \end{statement}

 \begin{definition}[Entropic interpolation] \label{def-x1}
The entropic interpolation  between $ \alpha$ and $ \beta$, built upon $\Xe$, is the unique (if it exists) solution $(p)$ of the least action problem \eqref{eq-x6}.
\end{definition}

For a generic definition, see Definition \ref{def-x2}. 

 We see that the empirical process $X^N$ of the system of independent copies of $\Xe$ gives us the action 
 \begin{align}\label{eq-117}
 (p_s) _{ s_0\le s\le s_1}\mapsto 
 \IS \Big(\uud   \| \dot\pp_s-\vce_s\|^2_{p_s}  
	+ \epsilon^2 I(p_s| r^ \epsilon_s)\Big)\, ds.
\end{align}

Comparing \eqref{eq-117} with the desired action \eqref{eq-115}, we observe that it is necessary 
\begin{enumerate}[\qquad(1)]
\item
to introduce the coefficient $\kappa_s$ into \eqref{eq-117} and 
\item
to remove from \eqref{eq-117} the Fisher information term.
\end{enumerate}

In order to do so, we are going to

\begin{enumerate}[\qquad(1)]
\item
introduce a splitting mechanism and 
\item
introduce some quantum force.
\end{enumerate}

\subsection*{Splitting $k$-mappings}

Let us first briefly describe the splitting mechanism; the quantum force will come later. We look at an empirical process  $(\Xb^N_s) _{ s_0\le s\le s_1}$ such that
for any $s,$  $\Xb^N_s$ is the empirical measure of $\lfloor \kappa_s N\rfloor$  $k$-mappings in $(\Rd)^k$ ($\lfloor a\rfloor$ denotes the integer part of the real number $a$). The factor $\kappa_s$ is analogous to $\kappa_s$ in formula \eqref{eq-121}. The trajectory of each  $k$-mapping is   a copy of $\Xe,$ but these copies are not independent. As   $s\mapsto\kappa_s$ is an increasing function, during any small time interval $[s,s+h],$ a fraction $\kappa'_s\,h+o _{ h\to 0}(h)$ of the $k$-mappings branches: each of them gives birth to a new $k$-mapping starting at the same place as its genitor and evolving in the future  according to the kinematics of $\Xe$ and independently of the other $k$-mappings. Although the number  of $k$-mappings increases with time, $\Xb^N$ is normalized so that its \emph{total mass remains constant}: $\Xb^N_s(\Rdk)=1$ for all $s.$ As a consequence, the random fluctuations of  $s\mapsto\Xb^N_s$ decrease with time\footnote{This can be quantified for instance by $ \mathrm{Var} (\int _{ \Rdk} f\,d\Xb^N_s) \underset{N\to \infty}\longrightarrow 0$ for any time $s$ and any bounded  function $f.$}. We prefer calling this branching  process a splitting process because the total mass remains constant.

\emph{The splitting mechanism acts as a cooling}.

We shall see that the large deviation rate function of $(\Xb^N) _{ N\ge 1}$ as $N$ tends to infinity  leads us to the action functional 
\begin{align} \label{eq-116}
(p_s) _{ s_0\le s\le s_1}
	\mapsto \IS\uud    \| \dot\pp_s-\dot\rr ^ \epsilon_s\|^2_{p_s}\kappa_s  \,ds
	+ \epsilon^2 \IS I(p_s| r^ \epsilon_s)\kappa_s\, ds.
\end{align}
This splitting mechanism  permits to plug $\kappa_s$ into \eqref{eq-117}, that is to pass from \eqref{eq-117} to \eqref{eq-116}.

\subsection*{Change of time}

As is done in \cite{ABB20}, things are  clearer after choosing $ \kappa_s=2 (s-s_0)$,   
changing  time as: $ s=s_0+ e ^{ 2t}$ (note that this maps $[s_0,s_1]$ onto $[- \infty,t_1]$) and setting \[q_t:=p_s=p _{s_0+ e ^{ 2t}}.\]  
 Indeed, \eqref{eq-115} becomes
$	
 (q_t) _{  t\le t_1}\mapsto  
 \ITT \uud\|\dot\qq_t-\dot\mm^\epsilon _t\|^2 _{ q_t}\,dt,
$	
where $\dot\mm^ \epsilon_t$ is the velocity of the time-changed heat flow \[m^ \epsilon_t:= r^ \epsilon _{ s}=r^ \epsilon _{s_0+ e ^{ 2t}}\] in the \OW-manifold, see \eqref{eq-97} and \eqref{eq-54}.  This leads us to the following 

\begin{definition}[$ \epsilon$-\MAG's action functional for a $k$-fluid]
It is
 \begin{align} \label{eq-118}
 (q_t) _{ t_0\le  t\le t_1}\mapsto  
 \IT \uud\|\dot\qq_t-\dot\mm^\epsilon _t\|^2 _{ q_t}\,dt,
\end{align}
\end{definition}
This is the $k$-fluid analogue of the main action functional of the Ambrosio-Baradat-Brenier particle system in \cite{ABB20}. 

Similarly, the action \eqref{eq-116} becomes
  \begin{align}\label{eq-119}
 (q_t) _{t_0\le  t\le t_1}\mapsto  \IT\uud    \| \dot\qq_t-\dot\mm^\epsilon _t\|^2_{q_t}\,dt
	+ \epsilon^2 \IT \sigma_t I(q_t| m^ \epsilon_t)\,  dt, 
\end{align}
with  $ \sigma_t=4 e ^{ 4t},$  see \eqref{eq-83}.

\begin{statement}[Equation of motion of $\epsilon$-\MAG\ for a $k$-fluid]
We show at Theorem \ref{res-25} that the Newton equation corresponding to $ \epsilon$-\MAG's action functional \eqref{eq-118} is
\begin{align} \label{eq-120}
\ddot \qq_t=\ddot \mm^ \epsilon_t,
\end{align}
where the double dot stands for the acceleration in the sense of the \OW-geometry.
\end{statement}
 
One may wonder why we started with time $s$ rather  than with time $t$. Exactly as in \cite{ABB20}, the reason for this is that, unlike $\dot\rr^ \epsilon,$ the velocity field $\dot \mm^ \epsilon_t$ which is the correct velocity field attached to $ \epsilon$-\MAG\ is not the current velocity of a simple stochastic differential equation. Hence it is easier to start with $s$ until \eqref{eq-116}, and then switch to time $t.$ 

\subsection*{Initial condition for $\Xe$}
One must be aware that the action functional \eqref{eq-118}, and therefore the corresponding dynamics, is not entirely specified by the stochastic differential equation \eqref{eq-23b} or similarly by the heat equation  \eqref{eq-114}. One also needs to fix an initial law 
\begin{align}\label{eq-x7}
\mathrm{Law}(\Xe _{ s_0})=r^ \epsilon _{ s_0}= m^ \epsilon _{ - \infty}=: \gamma \in\PRdk
\end{align}
 to know who is the solution $r^ \epsilon_s$ of the heat equation, and therefore who are $m^ \epsilon_t,$ $\dot\mm^\epsilon _t$ and $\ddot\mm^\epsilon _t$. \emph{The dynamics of \MAG\ crucially depends on the choice of this initial distribution $ \gamma\in\PRdk.$}  This is commented at Remark \ref{rem-02}-(i) and taken into account for instance in the  Definition \ref{def-01} with the notion of \MAG's action \emph{pushed by $ \lambda$}.

Remarkably, unlike \MAG's dynamics,  the  dynamics of any entropic interpolation built upon $\Xe$, recall Definition \ref{def-x1}, does not depend on the choice of the initial law $ \gamma$ of $\Xe$.  This is the content of  

\begin{statement}\label{res-x3}
For any $ \gamma\in\PRdk$, Problem \eqref{eq-x6} is equivalent to
\begin{align}\label{eq-x6b}
\inf _{ (p)} \IS \Big(\uud   \| \dot\pp_s-\vce_s( \gamma)\|^2_{p_s}  
	+ \epsilon^2 I(p_s| r^ \epsilon_s( \gamma))\Big)\, ds
\end{align}
where as in \eqref{eq-x6} the infimum runs through all the $(p):[s_0,s_1]\to \PRdk $ satisfying $p _{ s_0}=\alpha $ and $ p _{ s_1}=\beta $, but where we denote  $r^ \epsilon_s ( \gamma)$  the probability law of $\Xe_s$ when $\Xe$  solves \eqref{eq-23b} with the intial condition \eqref{eq-x7}. 
\end{statement}
This invariance\footnote{This can be interpreted as a gauge invariance.} holds because the dynamics of the entropic interpolations only depends on $\Xe$ via the collection of all its bridges:  $\Proba (\Xe\in\sbt\mid \Xe _{ s_0}=a, \Xe _{ s_1}=b),$ see  \cite{Leo12e}.

\subsection*{Gibbs conditioning principle}

Doing the above change of time: $ s=s_0+ e ^{ 2t}$, we look at 
\[\Yb^N(t):=\Xb^N(s)=\Xb^N(s_0+e ^{ 2t})\] and  \eqref{eq-116} becomes \eqref{eq-119}.

\begin{statement}\emph{(Gibbs conditioning principle).} \label{res-x1} For any probability measures $ \alpha$ and $ \beta$ on $\Rdk,$ conditionally on  $ \Yb^N(- \infty)\simeq \alpha$ and $ \Yb^N(t_1)\simeq \beta,$  the most likely trajectory $(q_t) _{ t\le t_1}\in  C([- \infty,t_1], \PRdk) $ of $\Yb^N$ as $N$ tends to infinity solves the least action problem associated with \eqref{eq-119}:
\begin{align}\label{eq-x3}
\inf _{ (q)} \ITT\uud    \| \dot\qq_t-\dot\mm^\epsilon _t\|^2_{q_t}\,dt
	+ \epsilon^2 \ITT \sigma_t I(q_t| m^ \epsilon_t)\,  dt
\end{align}
where the infimum runs through all the $(q):[- \infty,t_1]\to \PRdk $ satisfying $q _{ - \infty}=\alpha $ and $ q _{ t_1}=\beta .$ 
 \end{statement}

To arrive at \MAG's action functional \eqref{eq-118}, it remains to subtract 
$
\epsilon^2 \ITT \sigma_t I(q_t| m^ \epsilon_t)\,  dt
$
from  \eqref{eq-119}. Be aware  that this term does not vanish as $ \epsilon$ tends to zero, see \eqref{eq-18} and Appendix \ref{app-d}.  It is convenient to express this operation by means of the corresponding Newton equation of motion in the \OW-manifold because adding a term in an action functional amounts to apply some force field to the system, see Appendix \ref{app-b}. It happens that the force corresponding to the energy potential
\begin{align}\label{eq-x5}
+\epsilon^2  \sigma_t I(\sbt| m^ \epsilon_t),
\end{align}
with a plus sign\footnote{The energy potential in the action \eqref{eq-x3} is $-\epsilon^2  \sigma_t I(\sbt| m^ \epsilon_t)$, see Appendix \ref{app-b}.}, can be interpreted as a quantum force.

\begin{remark}
In view of Statement \ref{res-x3}, one can generate the Lagrangian $\uud   \| \dot\pp-\vce_s( \gamma_o)\|^2_{p}$ for some $ \gamma_o$, by subtracting from the  Lagrangian $\uud   \| \dot\pp\|^2_{p}  
	+ \epsilon^2 I(p| \Leb)$ (which corresponds to $ \gamma=\Leb$) of the entropic interpolation problem, the  "quantum" energy $\epsilon^2 I(p| r^ \epsilon_s( \gamma_o)).$ 
\end{remark}

\subsection*{Quantum versus entropy}

The change of sign in front of the entropic\footnote{The Fisher information is a rate of loss of relative \emph{entropy} along a heat flow.} potential $- \epsilon^2  I(\sbt|m)$ that is necessary to switch from a thermal evolution to a quantum one is another way of describing the standard quantization rule which allows to switch from the heat equation to the Schrödinger equation. 
 Since the very beginning of quantum mechanics, it was noticed that the Schrödinger equation 
\begin{align*}
\big(- i\hbar \partial_t  -   \hbar^2   \Delta  /2 +V\big)\Psi =0,
\end{align*}
where $\hbar$ is Planck's constant and $ \Psi$ is a complex wave function, looks like the parabolic heat equation
\begin{align*}
(- \epsilon \partial_t +   \epsilon^2   \Delta /2 +V ) f=0
\end{align*}
where $f$ is a real function and $ \epsilon>0$ is some parameter which, as $\hbar$, has the physical dimension of an action, i.e.\ [time]$\times$[energy]. One switches from one equation to the other by applying the rule
\begin{align}\label{eq-96b}
\epsilon \leftrightarrow \hbar,
\qquad
V \leftrightarrow -V,
\qquad
t \leftrightarrow -it, 
\qquad
f \leftrightarrow \Psi.
\end{align}
The mysterious part of this quantization rule is $t\! \leftrightarrow\! -it$, for nobody knows what a complex time is\footnote{Putting apart the fact that the physical status of time is still a controversial problem.}. Still, it works as a recipe. 

Nevertheless,  Schrödinger \cite{Sch31, Sch32} revealed    in 1931  a striking analogy between the \emph{dynamics} -- in contrast with the kinematics -- of the thermal and quantum evolutions. This  gave rise to the Schrödinger problem, see Section \ref{sec-nbm} for some details and \cite{Leo12e} for a survey. Inspecting forces and accelerations is what dynamics does. This implies that second order derivatives in time play a major role. Applying the above quantization rule $t \leftrightarrow -it,$ we see that $$ \partial_t^2 \leftrightarrow -\partial_t^2.$$ 
This is the main reason why \emph{a quantum force may be interpreted as minus an entropic force}. This is illustrated at Section \ref{sec-B-k-fluid} by means of the Madelung equations \eqref{eq-86} and \eqref{eq-86c}, relying on the articles \cite{vR11} by von Renesse and \cite{Co18} by Conforti. 

\emph{Annihilating some entropic force requires to apply some quantum force, and vice-versa.}

\subsection*{Quantum force?}

 In terms of an equation of  motion,  subtracting the Fisher information term from the action functional \eqref{eq-119}    to arrive at the desired action  \eqref{eq-118} amounts to apply the force field 
 \begin{align}\label{eq-122}
  - \epsilon^2 \sigma_t \gOW _{ q_t}  I(\sbt|m^ \epsilon_t),
 \end{align}
 where $\gOW$ is the gradient of the \OW-manifold. 
At Section  \ref{sec-B-k-fluid}, we propose the Physical hypothesis \ref{physhyp} that this force is a quantum one. This is supported by an analogy with   two known results about the dynamics of entropic interpolations \cite{Co18} and quantum interpolations \cite{vR11}, see Theorems \ref{res-21} and \ref{res-21b}. 

\subsection*{Open problems}

We do not investigate further the following challenging problems: 
\begin{enumerate}[(a)]
\item
What is the (quantum?) nature of the force field $-\gOW I(\sbt|m)?$
\item
What is the physical meaning of the relation between the opposite forces which are the entropic force: $ +\gOW I(\sbt|m),$ and the quantum  one:  $-\gOW I(\sbt|m)?$
\end{enumerate}

For simplicity, from now on we call any field of the form $-\gOW I(\sbt|m)$ a quantum force and any field of the form $+\gOW I(\sbt|m)$ an entropic force.

We believe that the analogies between the  dynamics of the entropic and quantum interpolations discovered  by Schrödinger \cite{Sch31,Sch32}, together with the \OW-geometry which permits to express them as simple Newton equations, could reveal some interesting interplay between entropic disorder and quantum order.

\subsection*{Concentration of matter}

At Appendix \ref{app-d}, some hints are presented to support the statement that the quantum force $- \epsilon^2 \sigma_t \gOW _{ q_t}  I(\sbt|m^ \epsilon_t)$ remains of order 1 as $ \epsilon$ tends to zero. Except for the change of time from $s$ to $t$, without the quantum pressure, the $k$-fluid would follow the thermal evolution of an entropic interpolation. In particular, far from the instants $t_0$ and $t_1$ it would look like the flow of a thermal diffusion. But we know  since Brenier's article \cite{Bre11} that \MAG's evolution exhibits  concentration of matter phenomena at any period of time. This is  illustrated at Appendix \ref{app-c}, see Figure  \ref{fig-matcon}, and numerically verified \cite{LBM24}, see Figure \ref{fig-numsim}.

At Theorem \ref{res-25} it is not only proved that \eqref{eq-120} holds, that is: $\ddot \qq_t=\ddot \mm^ \epsilon_t$, but also that there  exists a negligible subset $ \mathcal{N}$ of $\Rdk$ which is 
 a finite union of  vector subspaces with codimension at least 2, such that  for all $\yy\in \mathcal{N},$ the force field $\ddot \mm^ \epsilon_t(\yy)$ diverges as $ \epsilon$ decreases to zero, see \eqref{eq-64b}. This strong force is responsible for the concentration of matter in certain areas of space. According to our model, it is a direct \emph{consequence of quantum pressure} because without this force field we would observe a forward or/and a backward in time diffusive evolution without any strong concentration areas.

\subsection*{No microscopic picture}

When considering quantum forces, we abandon the idea of giving a microscopic description of a particle system and we stay at a macroscopic level. To give an idea of the problem we are facing when trying to recover the fluid dynamics associated to \eqref{eq-118} from the particle system leading to \eqref{eq-119}, let us spend some time with the entropic interpolations to present an  analogy with our problem, adopting the viewpoint advocated by Schrödinger when addressing his eponym problem \cite{Sch31, Sch32}.  

Any entropic interpolation $(p_s) _{ s_0\le s\le s_1}$ built upon the Brownian motion $\Xe$  between prescribed initial and final distributions, solves the Newton equation \cite{Co18}
\begin{align}\label{eq-x1}
\ddot\pp_s= \epsilon^2 \gOW _{ p_s} I(\sbt|\Leb),
\end{align}
while a McCann displacement interpolation (minimizer of the dynamic quadratic transport problem \cite{McC94,McC97})   between the same endpoint marginals solves the geodesic equation 
\begin{align}\label{eq-x2}
\ddot \pp_s=0
\end{align}
in the \OW-manifold \cite{Otto01,Vill09}. As the entropic interpolation is a thermal perturbation of the displacement interpolation \cite{Leo12a}, one should interpret \eqref{eq-x1} as a thermal force. 

To get a hint for \eqref{eq-x1}, let us go back to Statement \ref{res-x3} and choose $ \gamma=\Leb.$ With this uniform initial distribution, $\Xe$ becomes reversible, see \eqref{eq-89}. In particular, for all $s,$  $r^ \epsilon_s( \Leb)=\Leb$ and the current velocity vanishes, that is: $\vce_s(\Leb)=0$. It follows that \eqref{eq-x6b} writes as
\begin{align*}
\inf _{ (p)} \IS \Big(\uud   \| \dot\pp_s\|^2_{p_s}  
	+ \epsilon^2 I(p_s| \Leb)\Big)\, ds.
\end{align*}
An analogy with the standard least action principle in a finite dimensional Riemannian manifold leads us to \eqref{eq-x1}, see Appendix \ref{app-b} for instance.

Since a displacement interpolation is a mixture of geodesics of the state manifold ($\Rdk$ in the present case) and an entropic interpolation is a mixture of bridges of the reference path measure (Brownian bridges in the present case), we see that we switch from a deterministic evolution  to a stochastic one  when passing from \eqref{eq-x2} to \eqref{eq-x1}. The force field created by $\epsilon^2 \gOW _{ p} I(\sbt|\Leb)$ has a nontrivial effect. In particular, at a microscopic scale it transforms differentiable paths into nowhere differentiable Brownian paths, increasing  the temperature from zero to a positive one. 

Now, passing in the other direction from \eqref{eq-x1} to \eqref{eq-x2} amounts to add a quantum force field as  \eqref{eq-122} with the effect of annihilating the temperature.

This  article is not able to propose a full microscopic description of an interacting particle system leading to \MAG's dynamics. It eventually presents a partial\footnote{Partial, because the relevant quantum force field still needs to be identified.} macroscopic description of a self-interacting fluid. However, some parts of the dynamics are justified by the thermodynamic limit of a large particle system of Brownian particles described at Statement \ref{res-x1}.

\subsection*{Picking up the pieces}

\begin{enumerate}
\item[(i)]
We start with the empirical process $X^N$ given at \eqref{eq-x4} which is built upon a sequence of independent copies of the Brownian $k$-mapping $\Xe$, see \eqref{eq-23b}. It obeys the Gibbs conditioning principle described at Statement \ref{res-x2} and leads us to the action \eqref{eq-117}.
\item[(ii)]
The coefficient $ \kappa_s$ in the actions $A^ \epsilon$ and \eqref{eq-116} is necessary when doing the time change from $s$ to $t$ to obtain the right expression of the $ \epsilon$-\MAG's action \eqref{eq-118}. In order to incorporate $\kappa_s$ into the action \eqref{eq-117} attached to $X^N$, we replace  $X^N$ by  the empirical process $\Xb^N$ of a system of  splitting $k$-mappings. Performing the change of time from $s$ to $t$, $\Xb^N(s)$ becomes $\Yb^N(t).$ It obeys the Gibbs conditioning principle described at Statement \ref{res-x1} and leads us to the action \eqref{eq-119}.
\item[(iii)]
To pass from the action \eqref{eq-119} to the desired action \eqref{eq-118} it is necessary to subtract an action built upon a Fisher information, see \eqref{eq-x5}. Since we did not find any microscopic description of an interacting particle system that emulates the action \eqref{eq-x5}, we step back to a macroscopic model at the level of fluids. We propose the hypothesis that the force field \eqref{eq-122} corresponding to \eqref{eq-x5} could be interpreted as resulting from  a quantum pressure emanating from the fluid itself: it  is created by the matter profile $q=|\Psi|^2$. This quantum pressure partly acts against the entropic force resulting from the Brownian motions of the particles (which are $k$-mappings). It stiffens the fluid.\\ This still needs to be investigated.
\end{enumerate}

Revisiting  model \eqref{eq-121}, the parameter $N$ replaces $\eta,$ the splitting mechanism replaces the parameter $\kappa$ in the factor  $ \sqrt{\eta k\kappa_s ^{ -1}}\,d \mathsf{W}_s,$ and the quantum force field supersedes the enigmatic substitution of the forward velocity of $ \mathsf{Z} ^{ \epsilon,\eta}$ by the the current velocity of $\Xe$.

\subsection*{Outline of the article} \label{outline}
\begin{enumerate}
\item[3\,-\ ]
Section \ref{sec-motdef}  begins with a short description of the physical basis of the Monge-Ampère gravitation theory. We also sketch the early Universe reconstruction problem as the main motivation for this modified  theory of gravitation.
\item[4\,-\ ]
  Section \ref{sec-MAG-mapping} is dedicated to the exposition of the mathematics of \MAG, as introduced by Brenier in \cite{Bre11}.
\item[5\,-\ ]
   Section  \ref{sec-B-k-mapping} gives the details about the model \eqref{eq-121} which is the main object of the article \cite{ABB20} by Ambrosio, Baradat and Brenier. 
\item[6\,-\ ]
   At Section \ref{sec-MAG-fluid}, we introduce the analogues for fluids of  the action functionals $A$ and $A^ \epsilon$, such as \eqref{eq-115}. 
\item[7\,-\ ]
The  action functional \eqref{eq-115}/\eqref{eq-118} we mainly work with is introduced at Section  \ref{sec-fluid}.
\item[8\,-\ ]
 Newton's equation of motion  \eqref{eq-120} is ``proved'' at Section \ref{sec-Neweq}. A detailed description of the acceleration $\ddot \qq$ is also obtained at Theorem \ref{res-25}. It reveals  divergences of the force field at some specific places which are responsible for the concentration of matter, a welcome property when modeling gravitation.
 \end{enumerate}

 The last three sections are dedicated to the construction of the fluid that was presented during this introductory section.
 \begin{enumerate}
\item[9\,-\ ]
  At Section \ref{sec-nbm}, we recall already known facts about the dynamics of the solutions of the Schrödinger problem. This leads us to the action functional \eqref{eq-117}.
\item[10\,-\ ]
   Then, we introduce the factor $\kappa_s$, at Section \ref{sec-kappa} to arrive at the action functional \eqref{eq-116}. This is the place where we prove that the splitting mechanism that was previously described  transforms \eqref{eq-117} into \eqref{eq-116}. 
\item[11\,-\ ]
   Finally,  at Section \ref{sec-B-k-fluid}, we give some clues in favor of the hypothesis that  the  force \eqref{eq-122} has a quantum nature and address a  physics open problem about the relation between entropic disorder and quantum order. 
   \end{enumerate}

Sections  \ref{sec-motdef}, \ref{sec-MAG-mapping}, \ref{sec-B-k-mapping}, \ref{sec-nbm} and a large part of Section  \ref{sec-B-k-fluid} are  expository: most of their content is already known. The rest of the article consists of new material.

 { This program of investigation is far from being complete}: a short list of remaining things to be done is proposed at Appendix \ref{app-a}. We also provide at Appendix \ref{app-b} a basic reminder about the minimization of action functionals, and at Appendix \ref{app-OW} basic definitions and results about the \OW-geometry. A detailed  argument is exposed at Appendix \ref{app-d} to show that the quantum force does not vanish when $ \epsilon$ tends to zero. 
 Some simple analogies are presented at Appendix \ref{app-c}  to illustrate the concentration  of matter that results from \MAG's dynamics.
Finally,  at  Appendix \ref{app-e}, basic properties of the quantum potential are gathered.

\section{\MAG. Motivation and definition}
\label{sec-motdef}

\subsection*{Newtonian gravitation}

The equation of  motion of a test particle in gravitational interaction with a fluid is 
\begin{align} \label{eq-01}
\ddot x_t=-\nabla \varphi_t( x_t),\quad t\ge 0,
\end{align}
where the scalar potential $ \varphi$ solves the Poisson equation
\begin{align}
\Delta \varphi_t= \mu_t. \label{eq-02}
\end{align}
Here, $ x_t\in\Rd$ is the position of the test particle at time $t$, $\ddot x_t$ is its acceleration and  $ \mu_t(x)$ is the  density of the fluid at time $t$ and location $x.$ While Newton's equation  \eqref{eq-01} is physically correct (the mass of the test particle is irrelevant), for simplicity  we do not write the gravitational constant in   \eqref{eq-02}. 

In computational cosmology, the state space $\Rd$ is replaced by the flat torus 
\[ \mathbb{T}^d _L=\Rd/[0,L]^d\] with size $L>0$. Assuming that $ \mu_t( \mathbb{T}^d_L):=\int _{  \mathbb{T}^d_L} \mu_t(x)\, dx< \infty,$ without loss of generality one can  normalize $ \mu_t$ as a probability measure. 
\\
Periodic boundary conditions imply that for any regular $[0,L]^d$-periodic function $ \varphi$  
\begin{multline*}
\int _{[0,L]^d} \Delta \varphi\,d\Leb
	= \int _{ \Rd}  \1 _{[0,L]^d}\nabla\scal (\nabla \varphi)\,d\Leb\\
	=- \int _{ \Rd}\nabla \1 _{[0,L]^d}\scal \nabla \varphi\,d\Leb
	=\int _{ \partial [0,L]^d} \nabla \varphi\cdot \vec n\, d \sigma
	=0
\end{multline*}
where $\Leb$ stands for the Lebesgue measure, $\vec n$ is the outer unit normal vector on the boundary $ \partial [0,L]^d$ of $ [0,L]^d$ equipped with the surface measure $ \sigma=\Leb _{ | \partial [0,L]^d}.$ This is Stokes formula and its vanishing means that the mass of any distribution of matter evolving along the integral curves of the  vector field $\nabla  \varphi$  remains constant  since    $ \mathbb{T}^d _L$ has no physical boundary. 
\\
Hence, for a Poisson equation in $ \mathbb{T}^d _L$ to admit a solution it is necessary that its  right hand  term has a zero mass.  Let us replace  equation \eqref{eq-02} by
 $	
 \Delta \varphi_t= \mu_t- \lambda
 $	
 where $ \lambda$ is some probability measure on $ \mathbb{T}^d _L.$ The best natural choice for $ \lambda$ is the uniform unit volume measure   on the torus
 \begin{align*}
  \lambda_L = L ^{ -d} \ \Leb _{ | \mathbb{T}^d _L}.
 \end{align*}
 Indeed, doing this we see that, with the balanced Poisson equation 
 \begin{align}\label{eq-03}
 \Delta \varphi_t= \mu_t- \lambda_L\quad \textrm{on } \mathbb{T}^d _L,
 \end{align} 
 a uniform distribution of matter: $ \mu_t= \lambda_L,$ does not generate any gravitational force, as expected. Moreover, it successfully passes  the ``large box''-test
\begin{align}\label{eq-105}
\lim _{ L\to \infty} \eqref{eq-03}= \eqref{eq-02}.
\end{align}
Indeed,  letting $L$ tend to infinity, one recovers the standard Poisson equation \eqref{eq-02} in $\Rd$ because the background source term $ \lambda_L$ vanishes as $L$ tends to infinity.  For further physical justification of this model in cosmology, see  \cite{HBS91,BP97}.

\subsection*{Monge-Ampère gravitation (\MAG)}

 Brenier  \cite{Bre11} introduced a modified theory of gravitation  which, although not fundamental,  is effective for solving the \emph{early Universe reconstruction} (\EUR) problem. The  main feature of this modified theory  consists in replacing the Poisson equation \eqref{eq-03} by the Monge-Ampère equation
 \begin{align}\label{eq-04b}
L ^{ -d} \det( \mathbb{I}+ \Hess(L^d  \varphi_t))=\mu_t,
\quad \textrm{on}\ \mathbb{T}^d _L,
 \end{align}
 where $ \mathbb{I}$ is the identity matrix. 
 The couple of equations \eqref{eq-01} and \eqref{eq-04b}   is called \emph{Monge-Ampère gravitation}, \MAG\ for short.
Note that  \eqref{eq-04b} admits \eqref{eq-03} as its linearization in the limit of weak gravitation, that is   
 \begin{align*}
 \lim _{  \|\Hess \varphi\|\to 0} \eqref{eq-04b}= \eqref{eq-03}, 
 \end{align*}
 or more precisely
 $\Delta \varphi_t= \mu_t- \lambda_L+o_{ \|\Hess \varphi_t\|\to 0}(\|\Hess \varphi_t\|),$ $L$ being fixed.
Moreover, it is exactly \eqref{eq-03} in dimension one.  
 But in dimension $d\ge 2,$ it fails the ``large box''-test:
 \begin{align*}
\lim _{ L\to \infty} \eqref{eq-04b}\neq  \eqref{eq-02}.
 \end{align*}
 For instance, in dimension $d=3,$ 
 \begin{align}
 L ^{ -3} \det( \mathbb{I}+  \Hess(L^3 \varphi))
&=L ^{ -3}\det( \mathbb{I}+ L^3 \mathrm{diag}(a,b,c))
=L ^{ -3}(1+L^3a)(1+L^3b)(1+L^3b)\nonumber\\
&=L ^{ -3}+  \Delta \varphi+ L^3(ab+ac+bc)+L^6abc \label{eq-113},
 \end{align}
where $a, b$ and $c$ are the eigenvalues of $\Hess \varphi.$
We see 
that the dominating term as $L$ tends to infinity is not $ \Delta \varphi$ as desired, but $L^6\det(\Hess \varphi).$ This remains true in any dimension where the dominating term is $L ^{ d^2-d}\det(\Hess \varphi)$.
\\
Moreover,  when $d\ge 2$, there is no mixed asymptotic regime $(L\to \infty, \|\Hess \varphi\|\to 0),$ where \eqref{eq-03}\&\eqref{eq-04b} is an approximation of  Newtonian gravity which would be valid for \emph{any} $ \varphi.$
 Let us show it.
Denoting $H:=\|\Hess \varphi\|$, we see that  
\begin{align*}
L ^{ -d} \det( \mathbb{I}+ \Hess(L^d  \varphi))
	= \Delta \varphi+ L ^{ -d} + L ^{ -d}  \sum _{ 2\le n\le d}O ((L^d  H)^n).
\end{align*}
A regime where $ \Delta \varphi$ would be  the dominating term,   must satisfy 
\begin{align*}
H ^{ -1}\Big( L ^{ -d} + L ^{ -d}  \sum _{ 2\le n\le d}O ((L^d  H)^n)\Big)\longrightarrow 0,
\end{align*}
because $\Delta \varphi=O(H)\rightarrow 0.$ 
Since this term is of order $(L^dH) ^{ -1} + \sum _{ 1\le n\le d-1}(L^d  H)^n,$ this would imply that $L^dH$  tends simultaneously to $ \infty$ and $0$, a contradiction. \hfill$\square$

Nevertheless, despite this negative result, \MAG\ reveals to be   effective for approximately solving  the early Universe reconstruction  problem. This is illustrated by Figure \ref{fig-numsim} and will be partly explained  during a short discussion  at page \pageref{discussion}.

\subsection*{Early Universe reconstruction (\EUR)}

\begin{figure}[h]
		\includegraphics[width=8cm]{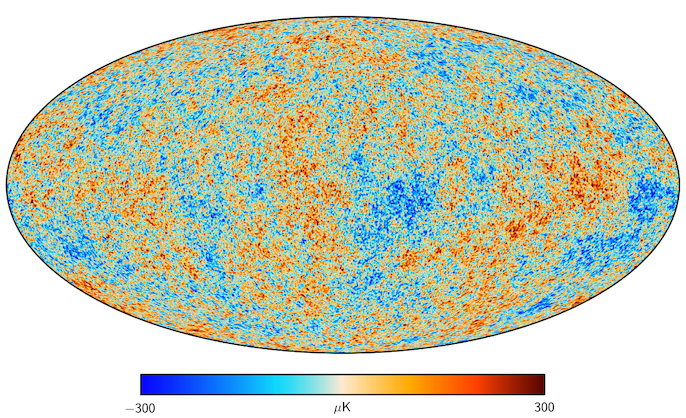}
		\captionof{figure}{Cosmic microwave background} \label{fig-cmb}
\end{figure}

This problem was   addressed by Peebles in  the seminal paper  \cite{Pee89}.
One specific feature with cosmology is that  the distribution of matter/energy of the very early Universe was highly uniform  as 
testified by the observation of the cosmic microwave background which exhibits a relative fluctuation from uniformity of order $10 ^{ -4}$: typically $3K\pm 300 \mu K$ (the unit $K$ is a Kelvin degree), see Figure \ref{fig-cmb}.

This very peculiar property of the initial condition is one   reason for the Monge-Ampère strategy to be successful when solving \EUR. However, the other reasons for  its effectiveness  are not fully understood at present time. We hope that this article  will be a step towards  a  clearer picture. 
\\
One aim of \EUR\ is to give  estimates of the field of very early fluctuations   from this uniformity \[\mu'_0(x):=\lim _{ t\to 0^+} \frac{ \mu_t(x)-L ^{ -d}}{t},\] a  crucial information to test the cosmic inflation theory and provide details about the initial quantum fluctuations.

\subsection*{An effective theory in computational cosmology}

\begin{figure}[h]
  \centering
  \subcaptionbox{Poisson\label{fig-numsim-a}}{%
    \includegraphics[width=5cm]{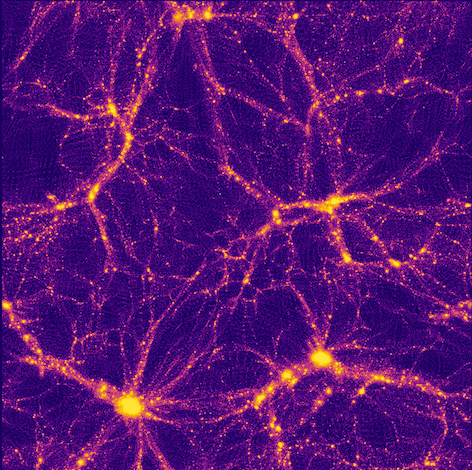}}
    \qquad\qquad
  \subcaptionbox{Monge-Ampère\label{fig-numsim-b}}
  {%
\includegraphics[width=5cm]{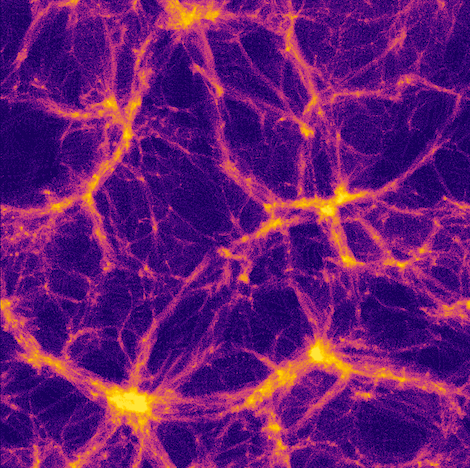}    
	}
  \caption{Typical structures of the actual  Universe\\ {\tiny Courtesy of Bruno Lévy}}\label{fig-numsim}
\end{figure}

As an illustration of the good performance of \MAG\ in cosmology, we provide at Figure \ref{fig-numsim} a couple of images representing typical structures of the actual Universe. Both are obtained by running numerical simulations starting from the same initial condition: a tiny perturbation of the uniform measure on $ \mathbb{T}^3_L.$ The left-hand side image (A)  is obtained using the balanced Poisson equation \eqref{eq-03}, while the right-hand side (B)  results from using the Monge Ampère equation \eqref{eq-04b} instead of the physically more realistic law \eqref{eq-03}.

\subsubsection*{\MAG\ works well with \EUR}\label{discussion}

Let us discuss a little bit about  \MAG\ being a good approximation of the Newtonian gravity in the special setting of the \EUR\ problem. The symmetric operator $\Hess \varphi$ admits an orthogonal basis of eigenvectors. Let $K:= \{\xi\in\RR^3; [\Hess \varphi] \cdot  \xi=0\}$ denotes its kernel, and $K^\perp$ be $K$'s orthogonal subspace in $\RR^3.$ Suppose that $\Hess \varphi$ admits at least one zero eigenvalue and at least one non-zero one. This means that both $K$ and $K^\perp$ are eigenspaces   and that  $( \mathrm{dim}(K), \mathrm{dim}(K^\perp))$ is either $(1,2)$ or $(2,1).$ 
Since $\Hess \varphi$ is the Jacobian matrix $\mathrm{Jac}( \nabla \varphi)$ of the opposite of the acceleration field $\ddot x=- \nabla \varphi,$ this implies that 
the acceleration is locally constant in direction $K$, so that any concentration of matter in direction $K^\perp$ is pushed by a force field in direction $K$, possibly fostering the formation of a structure of dimension $ \mathrm{dim}(K).$ 
\begin{enumerate}[(1)\quad]
\item
At places where $\Hess \varphi$ admits a single zero-eigenvalue, that is $\mathrm{dim}(K)=1,$ then 1D  singular structures (filaments) might appear.

\item
At places where $\Hess \varphi$ admits a double zero-eigenvalue, that is $\mathrm{dim}(K)=2,$ then 2D  singular structures (sheets) might appear.
\end{enumerate}
  Going back to \eqref{eq-113}, it appears that \eqref{eq-04b} is exactly  \eqref{eq-03} if two eigenvalues of $\Hess \varphi$ vanish. 
   In practice, for \MAG\ being close to Newtonian gravitation, one eigenvalue should be very close to zero to kill the leading term $L^6 abc$ and another one should be small enough to control the intermediate term $L^3(ab+ac+bc).$ 
 \\
 As can be seen at Figure \ref{fig-numsim-a}, it happens that $N$-body simulations based on the Newtonian gravitation reveal that, after some time, most of the matter is concentrated  in singular structures: sheets and filaments (Figure \ref{fig-numsim} depicts  2D-slices of  3D-cubes). This is a clue in favor of the good fit between these two theories in this special case. However, estimating the accuracy of \MAG\ for solving \EUR\ still remains a mathematical open problem.
 
 Remark that the above classification of singular structures in terms of $\Hess \varphi$ differs from the standard approach by Zeldovich \cite{Zel70} which is based on the Jacobian of a velocity field rather than an acceleration field. The way our classification complements Zeldovich's one will be explored elsewhere.
 
 More about the approximation of the Poisson equation by the Monge-Ampère equation can be found at Appendix \ref{app-MA-P}.

\subsubsection*{Pros and cons}

From a practical point of view, one  advantage of \MAG\ is that it provides us with  faster computations than the standard Poisson-based algorithms, because of its connection with optimal transport (this will be made precise later). More  about cosmological simulations using \MAG\ and their comparison with standard $N$-body simulations can be found in the recent article \cite{LBM24} by Lévy, Brenier and the second author.  From a theoretical point of view, another advantage is that its tight relation with optimal transport provides us with easy geometrical interpretations.  Moreover, it is shown in the article  \cite{BBM24} by Bonnefous, Brenier and one of us (RM) that \MAG\ can describe a  scalar field that adds a new long-range interaction to general relativity (a \emph{fifth force}) and is often evoked in modified theories of gravity such as Galileons, opening a promising domain to further explore.
\\
However, the main  drawback of \MAG\ is that it is not a fundamental theory of gravitation; it is only effective for solving \EUR.

From now on we drop the size $L$ by setting $L=1$, and denote $ \mathbb{T}^d= \mathbb{T} _{ L=1}^d.$

\subsection*{Zeldovich approximation}
 In \cite{Bre11} Brenier transposed Peebles' \EUR\ problem in the framework of the semi-Newtonian gravitational model of the early Universe where the trajectory $t\mapsto x_t\in \mathbb{T}^3$ of each particle (typically a cluster of galaxies!) satisfies the Newton-like equation of motion 
\begin{align}\label{eq-107}
\frac{2t}{3} \ddot x_t+ \dot x_t=- \nabla \varphi_t(x_t), \quad t\ge 0,
\end{align}
where the potential $ \varphi$ solves
\begin{align}\label{eq-108}
1+t \Delta \varphi_t=  \mu_t
\end{align}
and $ \mu_t$ is the matter density.
This describes classical Newtonian interactions taking place in an Einstein-de\,Sitter space corresponding to a Big Bang scenario.  The couple of equations \eqref{eq-107}\&\eqref{eq-108} is called the semi-Newtonian system ({\small SNS}).
\\
At time $t=0,$ we see that
\begin{align*}
 \mu_0\equiv 1, \quad \dot x_0=-\nabla \varphi_0(x_0),
 \quad \Delta \varphi_0(x)= \lim _{ t\to 0^+} \frac{ \mu_t(x)-1}{t}=: \mu'_0(x).
\end{align*}

\begin{remarks}\ \begin{enumerate}[(a)]
\item
One is  far from the classical $N$-body problem.

\item
Equation \eqref{eq-108} is simply \eqref{eq-03} with $L=1$ and $t \varphi_t$ instead of $ \varphi_t.$

\item
The physical assumption: $ \mu_0\equiv 1,$ is experimentally verified up to a very high precision level,  see Figure \ref{fig-cmb}.\end{enumerate}\end{remarks}

Zeldovich approximation \cite{Zel70} is simply
\begin{align*}
\widetilde x_t=x_0-t \nabla \varphi_0(x_0),\qquad 0\le t\le t_*(x_0),
\end{align*}
where $0<t_*(x_0)\le \infty$ is the first time where a collision with another particle occurs. In this  scenario, the initial fluctuation of the field generates a  path with  constant velocity until its first collision. 

Brenier proved in \cite{Bre11} that keeping \eqref{eq-107}, but  replacing 
 the balanced Poisson equation \eqref{eq-108} by the Monge-Ampère equation
 \begin{align}\label{eq-110}
\det( \mathbb{I}+t\Hess \varphi_t)=  \mu_t,
 \end{align}
one obtains a least action principle admitting the Zeldovich approximations as its \emph{exact} solutions. A similar reasoning will be exposed in a while to arrive at \eqref{eq-16} for the \MAG\ problem \eqref{eq-01}\&\eqref{eq-04}, see Definition \ref{def-08} below. Note that \eqref{eq-108} and \eqref{eq-110} are exactly \eqref{eq-03} and \eqref{eq-04b} with $L=1$ and different notations for the potential.

\subsection*{Optimal transport versus $N$-body simulation}

\begin{figure}
\centering
\begin{subfigure}{.5\textwidth}
  \centering
    \includegraphics[width=5cm]{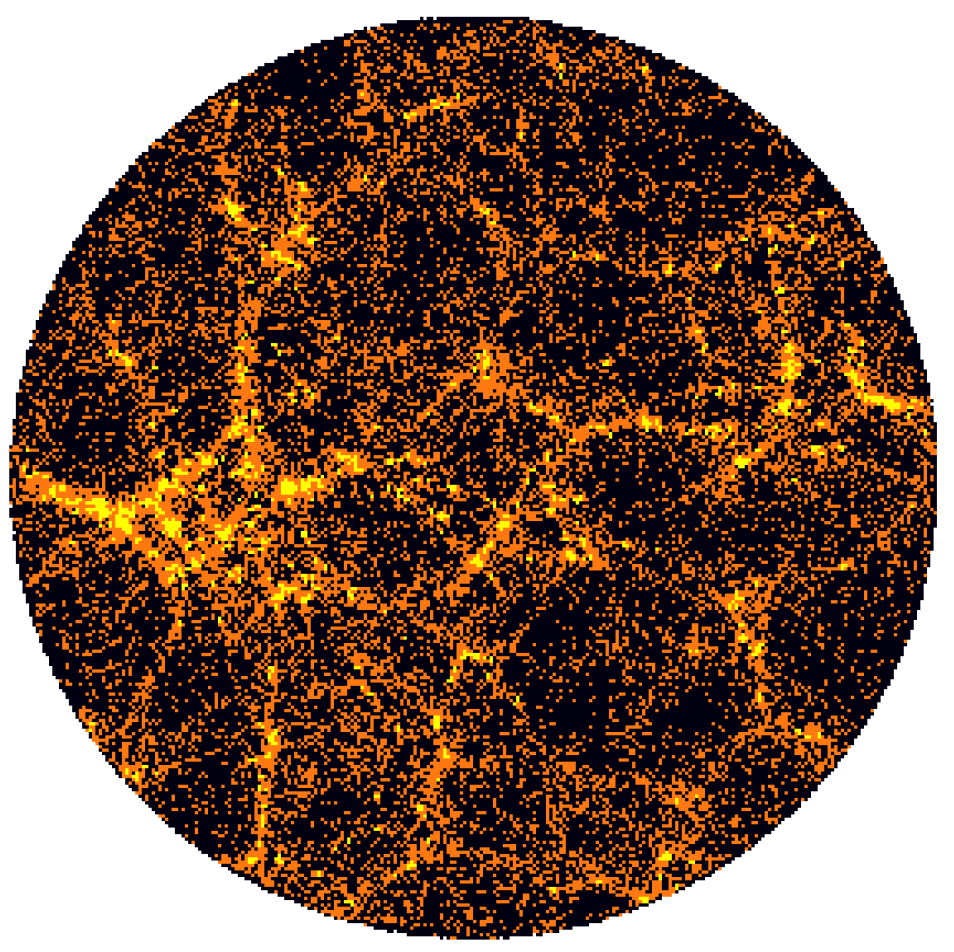}
\end{subfigure}%
\begin{subfigure}{.5\textwidth}
  \centering
\includegraphics[width=6cm]{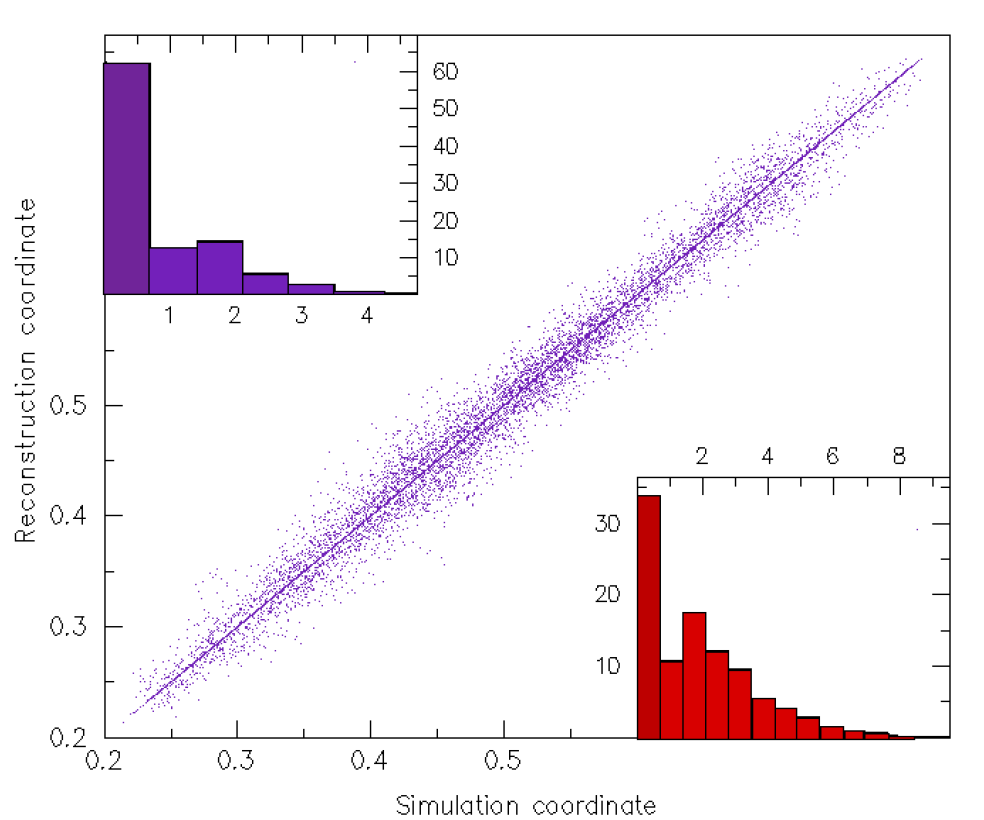}    
\end{subfigure}
  \caption{Optimal transport vs.\ $N$-body simulation, \cite{FMMS02}.}\label{fig-fmms}
\end{figure}

Frisch, Matarrese, Sobolevskii and the second author of the present article have shown  in    \cite{FMMS02}, more than twenty years ago, that with the simplified dynamics of the Zeldovich approximation,   \EUR\  is exactly the  Monge quadratic optimal transport problem between the initial uniform distribution of matter $ \mu_0\equiv 1$ and the  distribution of matter of the present epoch $ \mu_T,$ provided that the Zeldovich map $x_0\mapsto \widetilde x_T(x_0)$ is the gradient of a convex potential.  See also  \cite{BFHLMMS} for  more mathematical details. Numerical simulations in 
\cite{FMMS02} also demonstrate that the solution of  this quadratic optimal transport problem  is  close to the result of   $N$-body simulations as performed following \cite{HBS91,BP97} for instance.  Figure \ref{fig-fmms}  illustrates the comparison between a standard $N$-body simulation and a construction using optimal transport. More precisely, one compares the joint distributions of the couples of initial and final positions of the $N$-body simulation with the optimal transport plan between the initial and final marginal distributions of the $N$-body simulation. The yellow points of the left-hand side picture correspond to  significant errors of matching. The dots near
the diagonal on the right-hand side graphic  are a scatter plot of reconstructed (using optimal transport)  initial points vs.\ simulation initial points. The upper
left inset is a histogram (by percentage) of distances  between such points; 62\%\ are assigned exactly (up to the grid precision). The lower right inset is a similar histogram for reconstruction on a finer grid where 34\%\ are assigned exactly.

\subsection*{\MAG's definition}
The first definition of \MAG\ that appeared in the literature is the following

\begin{definition}[\MAG's approximation of {\small SNS}, \cite{Bre11}]
On the torus $ \mathbb{T}^d= \Rd/[0,1]^d,$ the dynamics defined by the couple of equations \eqref{eq-107}\&\eqref{eq-110}     is the \emph{Monge-Ampère gravitation} (\MAG) system approximating  the semi-Newtonian system ({\small SNS}) \eqref{eq-107}\&\eqref{eq-108}.
\end{definition}

\MAG's approximation of {\small SNS} is effective in the regime of short time and weak gravity.

From now on, following Brenier and his co-authors \cite{Bre16,ABB21,ABB20} we drop the relativistic dynamics \eqref{eq-107} and go back to the usual Newton equation \eqref{eq-01}. With this choice the underlying physics is  less realistic, but the mathematics are clearer and easier to handle. This simplification seems to us to be a reasonable choice for a preliminary approach of a more complete theory to be developed later.

 \begin{definition}[\MAG\ on $ \mathbb{T}^d $, absolutely continuous fluid, \cite{Bre11}] \label{def-08}
   The dynamics defined on the torus $ \mathbb{T}^d= \Rd/[0,1]^d$ by the couple of equations \eqref{eq-01} and
   \begin{align}\label{eq-04}
 \det( \mathbb{I}+\Hess \varphi_t)=\mu_t,
 \quad \textrm{on}\ \mathbb{T}^d ,
 \end{align}
    is called \emph{Monge-Ampère gravitation}, \MAG\ for short.
 \end{definition}

Furthermore, to simplify the mathematics and the presentation, we keep  Brenier's choice   in \cite{Bre11,ABB20} to consider \MAG\      in the whole space $ \Rd$ rather than in the torus $ \mathbb{T}^d,$ see Definition \ref{def-06} below.  However,   a severe  critic of this choice is expressed at Remark \ref{rem-02}-(i).

\section{\MAG. Action functional}
\label{sec-MAG-mapping}

This section is aimed at giving a short presentation of the Monge-Ampère gravitation  theory. It describes step by step the way for justifying  Brenier's Definition \ref{def-07} of \MAG's action functional.

 \subsection*{Optimal transport}

At first sight, it seems to be useless to replace the linear equation \eqref{eq-03} by a nonlinear one. But, be aware that solving the gravitational problem \eqref{eq-01}\&\eqref{eq-03}, or its standard counterpart \eqref{eq-01}\&\eqref{eq-02}, remains inaccessible in most situations. On the other hand, as will be seen in a moment, the connection of the Monge-Ampère equation \eqref{eq-04} with quadratic optimal transport permits a rather simple  geometric interpretation of   \MAG\ which leads to a fast numerical algorithm in \cite{Bre11}. 


Following Brenier  \cite{Bre11}, let us make precise the connection of \eqref{eq-04} with quadratic optimal transport.  
 Take a probability measure  $ \lambda$  on $\Rd$ and transport it by the measurable map $\Tf :\Rd\to\Rd$ to obtain its image 
 \[\mu:=\Tf \pf\lambda,\]
 defined by $ \mu(dy)= \lambda\big( (\Tf )^{ -1}(dy)\big).$
We have in mind that   $ \mu$  plays the role of the actual distribution of matter. The measure $ \lambda$ will be specified at \eqref{eq-07} so that some analogue of \eqref{eq-04} is satisfied, see \eqref{eq-08} below.
\\
  If  $ \mu$ is  absolutely continuous and  $\Tf $ is injective and differentiable almost everywhere, then $ \lambda$ is also absolutely continuous, the inverse map $(\Tf) ^{ -1}=:\Tb$  is differentiable $\mu$-a.e., $\Tb \pf \mu=\lambda$  and  the Monge-Ampère equation
\begin{align}\label{eq-05}
\mu(y)=\lambda(\Tb(y))\ |\det(\nabla\Tb)(y)|, \quad y\in\Rd,\ \textrm{a.e.}
\end{align}
is satisfied.
Here and in the remainder of this article, any absolutely continuous measure and its density with respect to Lebesgue measure are denoted by the same letter.
\\
The  Monge  optimal  transport problem which is relevant for our purpose is
\begin{align}\label{eq-09}
\inf _{ \Tb : \Tb \pf \mu=\lambda} \int_{\Rd} |y-\Tb (y)|^2\, \mu(dy),
\end{align}
where, as for equation \eqref{eq-05}, the unknown is the map $\Tb :\Rd\to\Rd$, while  the  probability measures $ \lambda$ and $\mu$ are prescribed. Brenier's theorem \cite{Bre91} tells us that when  $\int_{\Rd} |x|^2\,\lambda(dx)< \infty,$ $\int_{\Rd} |y|^2\, \mu(dy)< \infty$ and $\mu$ is absolutely continuous,  \eqref{eq-09} admits a unique solution $\Tb$. Moreover, 
\begin{align}\label{eq-10}
\Tb =\nabla \theta,\quad \textrm{a.e.,}
\end{align}
where $\theta:\Rd\to\RR$ is a convex function.  As Aleksandrov's theorem states that any convex function  on $\Rd$ is twice differentiable almost everywhere,  the Jacobian matrix $\nabla \Tb =\Hess \theta$ is  well-defined a.e.\ and positive semi-definite. Introducing the function 
\begin{align}\label{eq-06}
\varphi(y):= \theta(y)- |y|^2/2,
\end{align}
with \eqref{eq-10} equation \eqref{eq-05} reads as
$ 
\mu(y)=\lambda(\Tb(y)) \det( \mathbb{I}+\Hess \varphi)(y), 
$
$ y\in\Rd,\ \textrm{a.e.}
$
Finally, we observe that in the special case where the source distribution $ \lambda$ is uniform, that is: 
\begin{align}\label{eq-07}
\lambda(dx)= \1_D(x) \,dx
\end{align}
for some measurable subset $D\subset \Rd$ with a unit   volume $\Leb(D)=1,$ \eqref{eq-05} simplifies as
$
\mu(y)=\1 _{ \{\Tb(y)\in D\}}\, \det( \mathbb{I}+\Hess \varphi)(y),
$  $ y\in\Rd,\  \textrm{a.e.}$  With \eqref{eq-10} and \eqref{eq-06}, one can recast this equation   for the unknown $ \varphi:$
\begin{align*}
\1 _{ \{y+\nabla \varphi(y)\in D\}}\, \det( \mathbb{I}+\Hess \varphi)(y)=\mu(y), \quad y\in\Rd,\  \textrm{a.e.}
\end{align*}
It is an actualization of the analogue   \eqref{eq-04} of the balanced Poisson equation, once the state space  $ \mathbb{T}^d $ is replaced by $\Rd$.  

\begin{definition}
\label{def-06}
\emph{(Dynamics of \MAG\  pushed by $D$ in $ \mathbb{R}^d$. Absolutely continuous distribution of matter).} 
 Let $D\subset \ZZ$ satisfy $\Leb(D)=1.$ 
The dynamics of the \emph{Monge-Ampère gravitation pushed by the source set $ D$} is defined by
\begin{enumerate}[(i)]
\item
the following system of Newton's equations \eqref{eq-01}:
\begin{align*}
\left\{
\begin{array}{lcll}
\ddot X_t(\xi)&=&-\nabla \varphi_t( X_t(\xi)),\quad &0\le t<t_*,\\
X_0(\xi)&=& \xi, &t=0,
\end{array}
\right.
\quad \xi\in \supp( \mu_0),
\end{align*}
where the initial probability measure $ \mu_0$ is assumed to be absolutely continuous, 
\item
the distribution of matter
\begin{align*}
\mu_t=(X_t)\pf \mu_0,\qquad 0\le t< t_*
\end{align*}
where $t_*\in[0, \infty]$ is the infimum of all times $t$ such that $ \mu_t$ fails to be absolutely continuous, 
\item
and the Monge-Ampère equation 
\begin{align}\label{eq-08}
\1 _{ \{y+\nabla \varphi_t(y)\in D\}}\, \det( \mathbb{I}+\Hess \varphi_t)(y)=\mu_t(y), \quad y\in\Rd,\  \textrm{a.e.},\quad 0\le t<t_*.
\end{align}
\end{enumerate} \end{definition}

\subsection*{\MAG's force field}

Replacing \MAG's  equation \eqref{eq-04} by \eqref{eq-08}, 
the right hand side of Newton's equation \eqref{eq-01} is
\begin{align}\label{eq-11}
-\nabla \varphi(y)=y-\nabla \theta(y)=y-\Tb(y), \quad    y\in\Tf(D),\  \textrm{a.e.}
\end{align}
This is the explicit connection of \MAG\ with quadratic optimal transport. It holds both on $ \mathbb{T}^d $ and $ \mathbb{R}^d$, since on $ \mathbb{T}^d $ one chooses \[D= \mathbb{T}^d .\]

\begin{figure}
\centering
\begin{subfigure}{.5\textwidth}
  \centering
    \includegraphics[width=7cm]{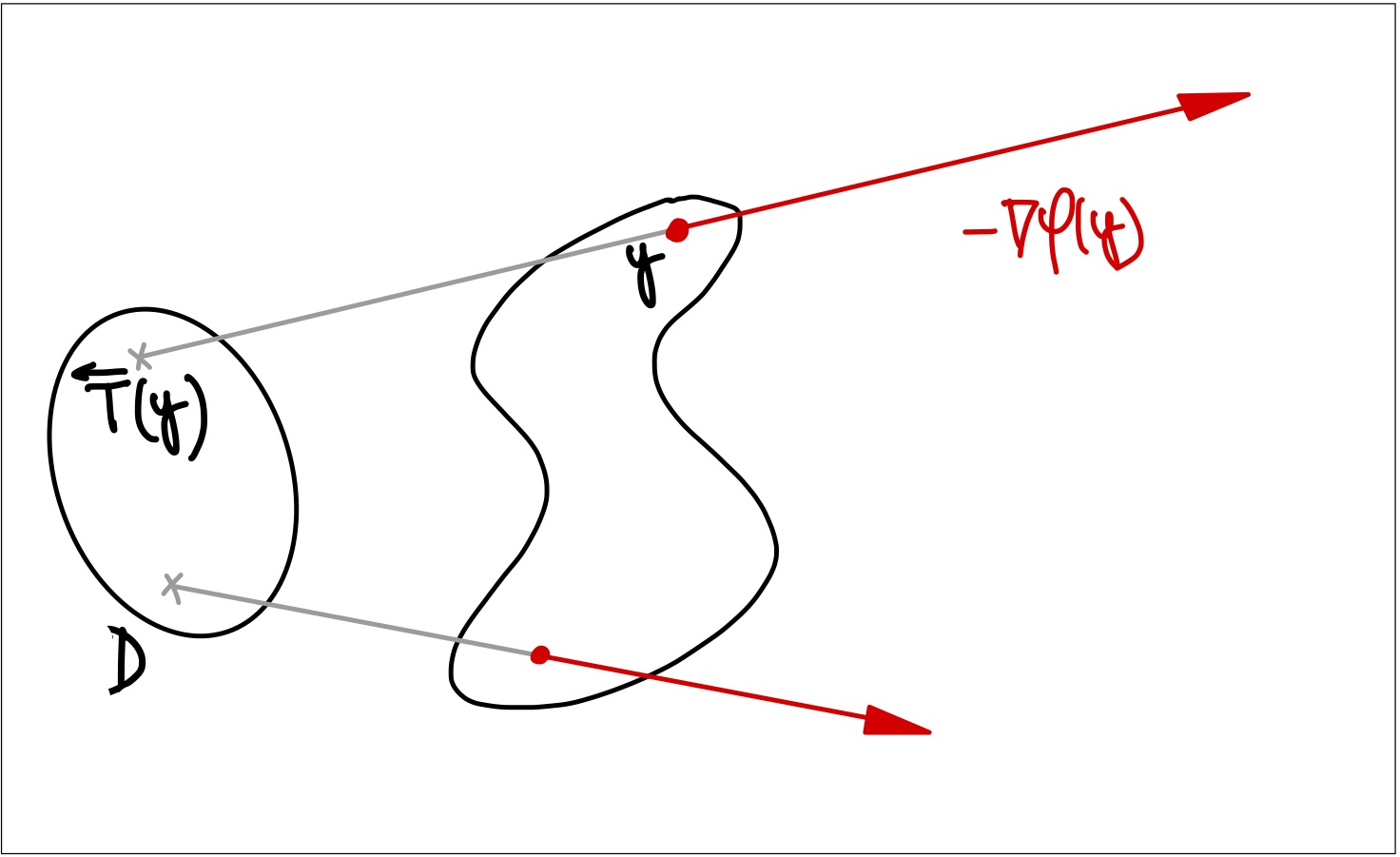}
  \caption{Matter outside $D$}
  \label{fig-pushed-by-D-a}
\end{subfigure}%
\begin{subfigure}{.5\textwidth}
  \centering
\includegraphics[width=7cm]{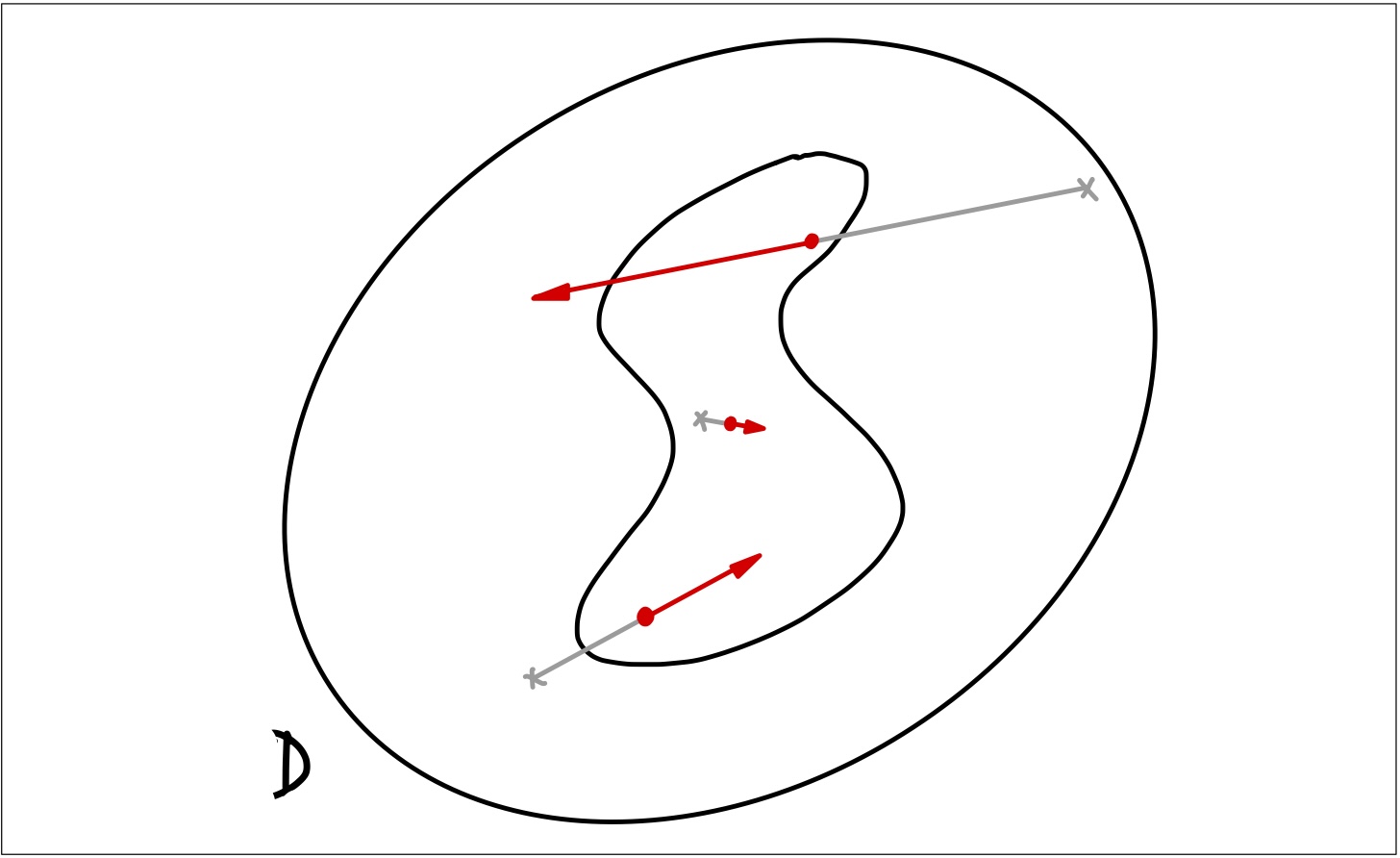}    
  \caption{Matter inside $D$}
  \label{fig-pushed-by-D-b}
\end{subfigure}
  \caption{The force field depends on $D$}\label{fig-pushed-by-D}
\end{figure}

\begin{remarks}\ \label{rem-02}
\begin{enumerate}[(i)]
\item
Since the force field $-\nabla \varphi$ depends on the choice of the source set $D$, it cannot be considered as a proxy for a physical gravitational field in $\Rd.$ This is illustrated at Figure \ref{fig-pushed-by-D}. The relative position of $D$ and the cloud of matter plays an essential role. In particular, one sees that Figure \ref{fig-pushed-by-D-a} illustrates anything but a self-attractive gravitation!
\\
  However, when the state space is $ \mathbb{T}^d $, as is usual in computational cosmology, the source measure $ \lambda$ is chosen to be the uniform probability measure on the whole set $ \mathbb{T}^d .$ Since \eqref{eq-08} becomes \eqref{eq-04} in this special case, this leads to \MAG\ gravitation as defined at Definition \ref{def-08}.
\\
Still, we go on working in $\Rd$ for simplicity.
\item
Note also that $-\nabla \varphi$ is only defined (almost everywhere) on the support of the target measure $ \mu.$ This is not an issue because we are only concerned with the evolution of $ \mu_t$.
\item
When considering such an evolution, it is not necessary to assume that the initial matter density $ \mu_0$ is equal to the source measure $ \lambda$ specified at \eqref{eq-07}. Indeed,  \emph{$ \lambda$ should be regarded as an artefact}.
\end{enumerate}
\end{remarks}

\subsection*{Polar factorization}

The last ingredient to be introduced for a geometric interpretation of \MAG\  is Brenier's polar factorization theorem \cite{Bre91}. During this subsection, $ \lambda$ is any  \emph{absolutely continuous} probability measure   on $\Rd$ satisfying $\int_{\Rd} |x|^2\, \lambda(dx)< \infty.$ The polar factorization theorem states that for  any measurable mapping $\yy :D\subset\Rd\to\Rd$   in 
\[H:=L^2 _{ \Rd}(D, \lambda),\] 
i.e.\ verifying $\|\yy \|_H^2:=\int_D |\yy (x)|^2\, \lambda(dx)< \infty,$ and such that the probability measure  \[\mu:=\yy \pf \lambda\] is also absolutely continuous, we have
\begin{align}\label{eq-12}
\Tb\circ \yy = \proj_S(\yy ),
\end{align}
where $\Tb$ is the optimal transport map from $ \mu$ to $ \lambda,$ and $\proj_S:H\to H$ is the orthogonal projection in the Hilbert space $H$ onto the subset
\begin{align*}
S:= \left\{ \xx  \in H: \xx\pf \lambda = \lambda\right\}
\end{align*}
of all $ \lambda$-preserving maps. 

\begin{center}
		\includegraphics[width=6cm]{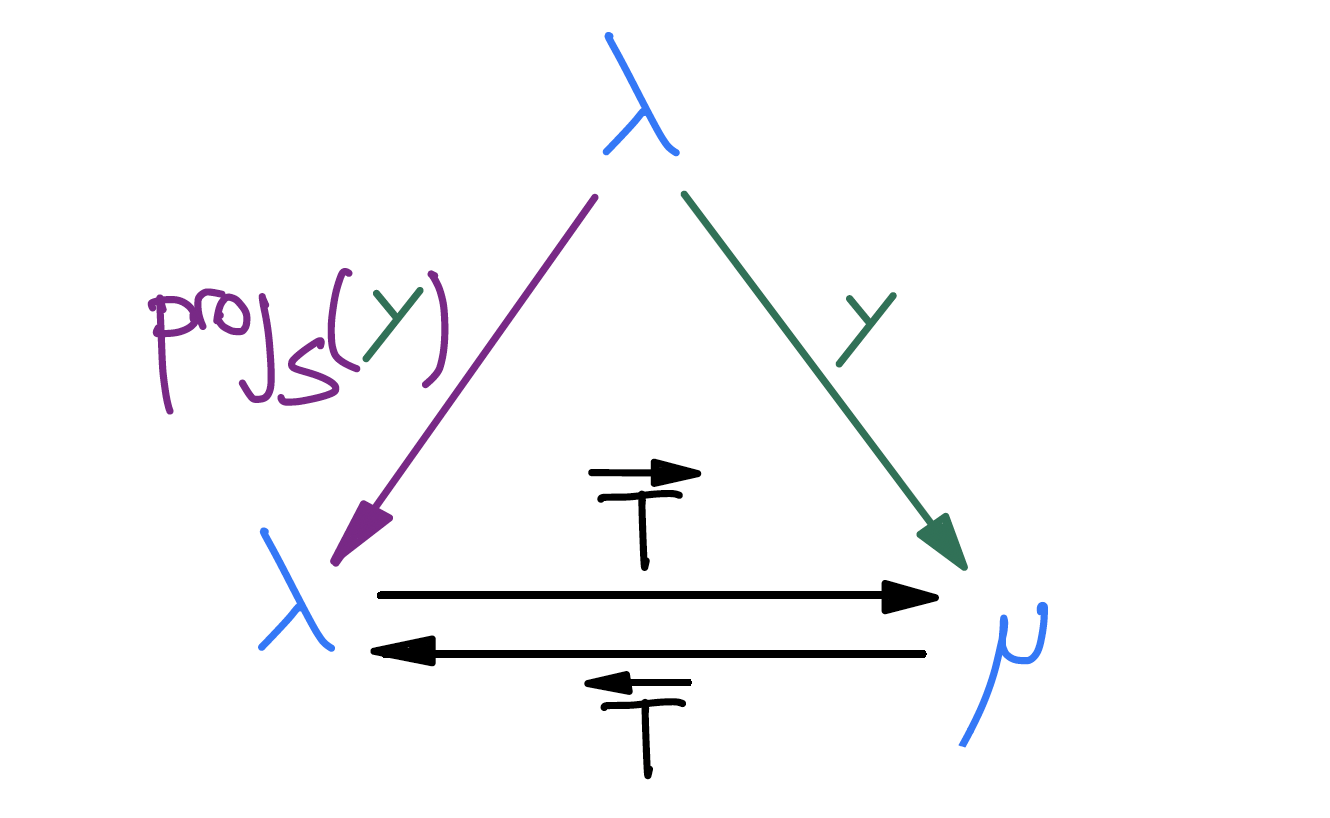}
\end{center}

From now on, elements in $H$ will be written with bold letters. 
Note that 
\begin{align}\label{eq-112}
\quad\forall \xx  \in S, \ 
\|\xx  \|_H^2=\int_D |\xx  (x)|^2\, \lambda(dx)=\int_D |x|^2\, \xx\pf\lambda(dx)=\int_D |x|^2\, \lambda(dx)=:r^2< \infty.
\end{align} 
Hence, $S$ is a subset of a sphere; in particular it is not convex. 
 The uniqueness of the projection of $\yy $ on $S$, which is implied by the assumed absolute continuity of $ \mu,$ is part the statement.  This theorem is often expressed the reverse way:
\begin{align*}
\yy =\Tf\circ \proj_S(\yy ),
\end{align*}
where $\Tf=(\Tb) ^{ -1}$ is the optimal transport map from $ \lambda$ to $ \mu.$ It is worth seeing $ \proj_S(\yy )$ as a type of permutation preserving $ \lambda$ primer to the least cost mapping $\Tf.$

\subsection*{Geometric expression of \MAG}
In view of \eqref{eq-01}, \eqref{eq-11} and \eqref{eq-12}, the collection of mappings $(\yy _t; t\ge 0)$ in $H$ is driven by \MAG\ if
\begin{align}\label{eq-13}
\ddot \yy _t=\yy _t-\proj_S(\yy _t), \quad 0\le t<t_*.
\end{align}
It is a detailed description of the evolution  of the infinite particle system $(\yy _t(x); x\in D) _{ 0\le t<t_*}.$ In particular it gives the evolution of the fluid density $\mu_t:=(\yy _t)\pf\lambda.$

\subsection*{Concentration of matter}

However, for \eqref{eq-13} to be valid it is necessary that $ \mu_t$ is absolutely continuous to ensure that $\yy_t$ admits a \emph{unique} projection  on $S.$ Considering the converse of this implication, one sees that \emph{as soon as $\yy_t$ has several projections on $S$, the matter distribution $ \mu_t$ becomes singular with respect to Lebesgue measure.} This appearance of singularities should be interpreted  as  concentration of matter. 
 As $S$ is  non-convex,  this happens as a rule. 
 \\
 Choosing \emph{one} candidate $\widehat{\proj}_S(\yy _t)$ as a function of the whole {set} $\Proj_S(\yy _t)$ of all the orthogonal projections of $\yy_t$ on $S,$ determines the dynamics when concentration of matter occurs.  
 A natural one is 
\begin{align}\label{eq-19}
\widehat{\proj}_S(\yy ):=\proj _{S(\yy )}(\yy )
\end{align}
where  $S(\yy ):=\cl\cv (\Proj_S(\yy ))$ is the closed convex hull of  $\Proj_S(\yy ).$
For a   justification of this choice, see Appendix \ref{app-c}.

\begin{definition}[\MAG\ dynamics allowing for  matter concentration, \cite{Bre11}]
 \MAG\ dynamics described at Definition \ref{def-06} is only valid for an absolutely continuous distribution of matter $ \mu_t.$ In view of \eqref{eq-13}, an extension allowing for matter concentration is 
\begin{align}\label{eq-14}
\ddot \yy _t=\yy _t-\widehat{\proj}_S(\yy _t)\in H, \quad t\ge 0.
\end{align}
\end{definition}

\begin{remarks}\ \begin{enumerate}[(a)]
\item
Because of the introduction of the extension $\widehat{\proj}_S$ of $\proj_S$ when $\Proj_S$ contains several points,  this definition allows for mass concentration: time $t$ goes beyond $t_*$. 
\item
It depends on the choice of the extension $\widehat{\proj}_S$ of $\proj_S.$ 
\item
  Later, it will be convienient to work with action functionals, rather than  with their Euler-Lagrange solutions.    One must be aware that the connection between \eqref{eq-14} and the natural candidate \eqref{eq-15} below  for a corresponding action functional is not fully established, see Remark \ref{rem-01}-(b).
\end{enumerate}\end{remarks}

\subsection*{Action functional}
A short reminder about action functionals and their connection with the equations of motion is proposed at Appendix \ref{app-b}.

In \cite{Bre11}, it is emphasized that the function
\begin{align*}
 \Phi (\yy ):= \inf _{ \xx  \in S}\| \yy -\xx  \|_H^2/2
 	=\|\yy \|_H^2/2-\Pi_S(\yy )+r^2/2,\quad \yy \in H,
 \end{align*}
 where \[\yy\mapsto\Pi_S(\yy ):=\sup _{ \xx  \in S} \left\langle \xx  ,\yy  \right\rangle _H,\]
is differentiable at any $\yy $ which admits a unique projection  on $S$.
Moreover, for any such "nice" $\yy $,
\begin{align}\label{eq-111}
\nabla \Phi(\yy )=\yy -\nabla\Pi_S(\yy )=\yy -\proj_S(\yy ).
\end{align}
This implies
\begin{align*}
\Phi(\yy )=\|\nabla \Phi(\yy )\|_H^2/2,
\end{align*}
and also that a Lagrangian associated to Newton's equation \eqref{eq-13} is 
\begin{align} \label{eq-56}
\|\dot \yy \|_H^2/2+\Phi(\yy ).
\end{align}
Plugging the last but one identity into the last  one,  the Lagrangian becomes 
$	
\|\dot \yy \|_H^2/2+\Phi(\yy )=\|\dot \yy \|_H^2/2+\|\nabla \Phi(\yy )\|_H^2/2
	=\|\dot \yy -\nabla\Phi(\yy )\|_H^2/2+ \langle\nabla\Phi(\yy ),\dot \yy \rangle_H.
$	
As $\langle\nabla\Phi(\yy ),\dot \yy \rangle_H$ is a null Lagrangian (because $ \displaystyle{\langle\nabla\Phi(\yy _t),\dot \yy _t\rangle_H= \frac{d}{dt}\Phi(\yy _t)}$ is a total derivative), one can choose the alternate Lagrangian
$
\|\dot \yy  -\nabla\Phi(\yy )\|_H^2/2
$
without modifying the  dynamics. Its leads to the action functional
\begin{align}\label{eq-16}
\IT\uud \|\dot \yy _t -\nabla\Phi(\yy _t)\|_H^2\ dt
= \IT\uud\|\dot \yy _ t -\yy _t+\proj_S(\yy _t)\|_H^2\ dt,
\end{align}
which becomes meaningless as soon as $t_1$ is larger than the first time $t^*$  when the set $\Proj_S(\yy _{ t^*})$ contains at least two elements.
Finally, Brenier proposes the following

\begin{definition}[Action functional of \MAG\ pushed by $D$, \cite{Bre11}] \label{def-03}
It is
\begin{align}\label{eq-15}
\IT\uud\|\dot \yy _t -\widehat\nabla\Phi(\yy _t)\|_H^2\ dt
=\IT\uud \|\dot \yy _ t -\yy _t+\widehat\proj_S(\yy _t)\|_H^2\ dt,
\end{align}
  the extended gradient $\widehat\nabla \Phi(\yy )$ being defined  by 
 \begin{align*}
 \widehat\nabla \Phi(\yy )
 	:=\yy -\widehat\nabla\Pi_S(\yy )
 	=\yy  -\widehat\proj_S(\yy )
 \end{align*}
 where $\widehat\nabla\Pi_S(\yy )$ is the (unique) element  with minimal norm in the subdifferential $ \partial\Pi_S(\yy )$ of the convex function $\Pi_S$.
 \end{definition}

\begin{remarks}\ \label{rem-01}\begin{enumerate}[(a)]
\item
In view of \eqref{eq-56}, another candidate for an action functional of \MAG\ is
\begin{align*}
\|\dot \yy \|_H^2/2+ \|\widehat\nabla \Phi(\yy )\|_H^2/2. 
\end{align*}
Unlike this action, the action \eqref{eq-15} looks like a Freidlin-Wentzell large deviation rate function. This   will be exploited at Section \ref{sec-B-k-mapping}, see \eqref{eq-21}.
\item
It is noticed in \cite{Bre11} that it is not clear that these  action functionals are equivalent and that one of them admits \eqref{eq-14} as its Euler-Lagrange equation. 
\item
 It happens that $\widehat\nabla\Pi_S(\yy )=\widehat\proj_S(\yy ),$ as defined at \eqref{eq-19}. See Proposition \ref{res-05}.
\end{enumerate}\end{remarks}

\subsection*{Extension to a discrete source measure}

Let us replace the normalized volume measure $ \lambda$   of some  set $D\in \Rd$, see \eqref{eq-07}, by  the normalized counting measure
\begin{align}\label{eq-25}
\lk:= \frac{1}{k}\sum _{ 1\le j\le k} \delta _{ x_j}
\end{align}
 on some finite set $D ^{ (k)}:=\{x_1,\dots, x_k\}\subset \Rd.$  Any mapping $\yy :D ^{ (k)}\to\Rd$ is encoded by the vector $\yy:=(\yy (x_1),\dots, \yy (x_k))\in (\Rd)^k$ whose squared norm in the Hilbert space $H:=L^2 _{ \Rd}( \lk)$  is 
 \[
 \|\yy \|_H^2=k ^{ -1}\sum _{ 1\le j\le k} |\yy (x_j)|^2=k ^{ -1}\|\yy\|^2
 \]
  where $|\sbt|$ and  $\|\sbt\|$ are respectively the  Euclidean norms on $\Rd$ and $ (\Rd)^k .$ Hence,  
 \begin{align*}
 H\simeq  (\Rd)^k .
 \end{align*}
 \begin{definition}[$k$-mapping] \label{def-km}
 When the source measure is the discrete measure $\lk$,  we call any application $\yy: D ^{ (k)}\to\ZZ$ in $H$ a \emph{$k$-mapping}.
 \end{definition}
Keeping our convention of writing elements of $H$ with bold letters, any $k$-mapping will also be written this way. 
 \\
The set $S$ of all $k$-mappings preserving $ \lk$ is
 \begin{align}\label{eq-20}
 S=\{\xx  ^\sigma; \sigma\in \mathcal{S}\}
 \end{align}
  where $ \mathcal{S}$ is the set of all permutations of $k$ elements and  
  \[\xx  ^\sigma:=(x _{ \sigma(1)},\dots, x _{ \sigma(k)})\in H,
  \qquad  \sigma\in \mathcal{S}.
  \]     
Brenier proposed in \cite{Bre16} to extend Definition \ref{def-03} to this semi-discrete\footnote{The term \emph{semi-discrete }refers to the fact that the source measure is discrete and the target one is diffuse. Of course, this diffuseness is only theoretical (this is necessary to apply Brenier's theorem before    concentration of matter occurs), because any numerical method requires to discretize the target measure.} setting which is natural to implement numerical simulations \cite{LBM24}.

\begin{definition}[Action functional of \MAG\ pushed by $\{x_1,\dots, x_k\}$, \cite{Bre16}] \label{def-07}
It is still \eqref{eq-15} as in Definition \ref{def-03}, but $ H\simeq  (\Rd)^k$ and $S$ is given by \eqref{eq-20}.
 \end{definition}

This definition is the basis for the definitions of new action functionals  we will work with in the rest of the article, see  \eqref{eq-28b} and  Definitions  \ref{def-01}, \ref{def-02} and  \ref{def-05}. At present time, it is only  justified by analogy with Definition \ref{def-03}. Its effectiveness in the numerical simulations suggests that one should recover \eqref{eq-15} by letting $k$ tend to infinity and preparing $ \lambda ^{ (k)}$ such that $\lim _{ k\to \infty} \lambda ^{ (k)}=\1 _{ D}\,\Leb.$

\subsection*{Non-locality} \label{sec-nonlocal}

Newton's gravity described by equations \eqref{eq-01} and \eqref{eq-02}  is non-local: a modification of $ \mu_t$ in a neighborhood of some location  affects instantaneously the force field at any remote location. 
Newton considered action at a distance to be an inadequate model for gravity. He even wrote: it is  \emph{so great an absurdity that I believe no man who has in philosophical matters a competent faculty of thinking can ever fall into it.} As is well known, this defect has been corrected by Einstein's theory of general relativity.
\\
Non-locality is also inherent to \MAG, as can be seen through  its fundamental equation \eqref{eq-14}. This is illustrated at the very end of the Appendix section \ref{app-c} (page \pageref{sec-non-local-2}).

\section{\MAG.  Particle system} \label{sec-B-k-mapping}

This section is aimed at giving a  presentation of  a random particle system, see \eqref{eq-21} below,  introduced by Ambrosio, Baradat and Brenier in \cite{ABB20} whose dynamics is related to the action functional of \MAG\ pushed by $\{x_1,\dots,x_k\}$ introduced at Definition \ref{def-07}, see \eqref{eq-22} below. 


 Let us  introduce the following  stochastic differential equation in the set  $\Rdk $ of $k$-mappings
 \begin{align}\label{eq-23}
 \mathsf{X} ^{ \epsilon} _s= \mathsf{X} _{ 0}+ \sqrt{ \epsilon}\, \mathsf{B} _{ s},\quad s_0\le s\le s_1,
 \end{align}
 where $\mathsf{B}$ is a standard Brownian motion  in $\Rdk $ starting from zero,  $ \epsilon>0$ is a fluctuation parameter which is intended to tend to zero, and the law of the initial position $ \mathsf{X} _{ s_0}$ in $\Rdk $ is
 \begin{align}\label{eq-29}
r _{ s_0}:= \textrm{Law}( \mathsf{X} _{ s_0})= \frac{1}{k!} \sum _{ \sigma\in \mathcal{S}} \delta _{ \xx  ^\sigma}.
 \end{align}
The process $ \mathsf{X}^ \epsilon _s=( \mathsf{X}^{\epsilon,1}_s,\dots, \mathsf{X}^{\epsilon,k}_s)$ describes a cloud of $k$ \emph{indistinguishable} Brownian particles in $\Rd$ starting from $\{x_1,\dots,x_k\}$; indistinguishability being a direct consequence of the choice \eqref{eq-29} of the initial law. It is a random \emph{path of $k$-mappings}. The $j$-th coordinate  $ \mathsf{X} ^{ \epsilon,j}$ is the path of the $j$-th particle starting from  the $j$-th random draw without replacement  from the set $\{x_1,\dots, x_k\}.$
Remark that although   the coordinates $ \mathsf{X} ^{ \epsilon,j},$  $1\le j\le k,$ are correlated, they share the same law with initial distribution $ \lk,$ see \eqref{eq-25}. 
 
 For any $s>s_0$, the law  of $ \Xe _s$ is the following mixture of   Gaussian measures  in $\Rdk $ with means $\xx  ^\sigma$ and covariance matrices $ \epsilon  (s-s_0)\, \mathbb{I}:$
 \begin{align}\label{eq-45}
 r^ \epsilon_s(d\xx)
 	= \frac{1}{k!} \sum _{ \sigma\in \mathcal{S}}  (2\pi \epsilon  (s-s_0)) ^{ -dk/2}
	\exp \left( - \frac{\|\xx -\xx  ^\sigma\|^2 }{2 \epsilon  (s-s_0)}\right) \, d\xx .
 \end{align}
 It solves the heat equation
 \begin{align}\label{eq-46}
  \partial_s r^ \epsilon= \epsilon \Delta r^ \epsilon/2
 \end{align}
 which rewrites as the continuity equation
 \begin{align*}
 \partial_s r^ \epsilon+\nabla\scal ( r^ \epsilon\vce)=0
 \end{align*}
 with the current velocity field $\vce$ on $\Rdk $ of the diffusion process $\Xe$  defined for any $s>s_0$ and  any $\xx\in\Rdk $ by
 \begin{align} \label{eq-17}
 \begin{split}
\vce_s(\xx )&=- \epsilon\nabla\log \sqrt{r^ \epsilon_s}(\xx )\\
 	&=- \frac{ \epsilon\nabla r^ \epsilon_s}{2 r^ \epsilon_s}(\xx )
	= \frac{\sum _{ \sigma\in \mathcal{S}} (\xx -\xx  ^\sigma) \exp \big( - {\|\xx -\xx  ^\sigma\|^2 }/ (2\epsilon (s-s_0))\big)}{ 2 (s-s_0) \sum _{ \sigma\in \mathcal{S}}  \exp \big( - {\|\xx -\xx  ^\sigma\|^2 }/ (2\epsilon (s-s_0))\big)}\\
	&= \frac{1}{2 (s-s_0)} \left( \xx - \frac{\sum _{ \sigma\in \mathcal{S}}  \xx  ^\sigma\exp \big( - {\|\xx -\xx  ^\sigma\|^2 }/ (2\epsilon (s-s_0))\big)}{\sum _{ \sigma\in \mathcal{S}}  \exp \big( - {\|\xx -\xx  ^\sigma\|^2 }/ (2\epsilon (s-s_0))\big)}\right).
 \end{split}
 \end{align}
Letting $ \epsilon$ tend  to zero, we see with the Laplace principle that for any   $s>s_0,$
$\lime\vce_s(\xx )=(2 (s-s_0)) ^{ -1}(\xx - \xx   ^{ \sigma(\xx )})$ where $ \xx   ^{ \sigma(\xx )}$ is the closest point from $\xx $ among all the $\xx  ^\sigma$ in $S,$  \emph{provided that this closest point is unique}. In view of \eqref{eq-20}, under this uniqueness assumption this means that
\begin{align}\label{eq-18}
\lime\vce_s(\xx )= \frac{\xx -\proj_S(\xx )}{2(s-s_0)},
\end{align}
a formula similar to the right hand side of equation \eqref{eq-13}. See Lemma \ref{res-y06} at Appendix \ref{app-d} for a detailed proof.   

Getting rid of the factor $(2 (s-s_0)) ^{ -1}$ will  be a matter of  change of time, see the parameter setting \ref{hyp-02} below.

Since the action functional \eqref{eq-16} can be read as some large deviation rate function, it is proposed in \cite{ABB20} to consider the stochastic differential equation in $\Rdk $
\begin{align}\label{eq-21}
d \mathsf{Z} ^{ \epsilon, \eta}_s
	=\vce_s( \mathsf{Z} ^{ \epsilon,\eta}_s)\,ds+ \sqrt{\eta k \kappa_s ^{ -1}}\, d \mathsf{W}_s ,\quad 0<s_0\le s\le s_1,
\end{align}
where $\vce$ is the current velocity \eqref{eq-17}, $ s\mapsto\kappa_s$ is a positive function, $ \eta>0$ is a parameter intended to  decrease to zero and $ \mathsf{W}$ is a standard Brownian motion on $\Rdk $. The Freidlin-Wentzell  large deviation principle roughly states that, for any fixed $ \epsilon>0$ and any initial state $\zz _o$ in $\Rdk ,$ when $\eta$ tends to zero 
\begin{align*}
\Proba( \mathsf{Z} ^{ \epsilon,\eta}\in \sbt\mid  \mathsf{Z} ^{ \epsilon,\eta} _{ s_0}=\zz _o) \underset{\eta\to 0}\asymp \exp \left( - \eta ^{ -1}\ \inf _{ \zz\in \bullet, \zz(s_0)=\zz _o}\widetilde A^ \epsilon( \zz)\right) 
\end{align*}
with  the large deviation rate function 
\begin{align}\label{eq-30}
\widetilde A^ \epsilon [( \zz_s) _{ s_0\le s\le s_1}] = 
 \IS\uud \|\dot \zz_s-\vce_s( \zz_s)\|_H^2\,  \kappa_s \,ds
\end{align}
where $ ( \zz_s) _{ s_0\le s\le s_1}$ stands for a generic absolutely continuous path taking its values in $\Rdk $. 
Remember that $\|\sbt\|_H=k ^{ -1/2}\|\sbt\| _{ \Rdk }$ and also that the factor $k $ is part of the diffusion coefficient in \eqref{eq-21}.
 It is proved in \cite{ABB20}, see also \cite{ABB21}, that 
\begin{align}\label{eq-24}
\Gamma\textrm{-}\lim_{ \epsilon\to 0}\widetilde A^ \epsilon\big[ ( \zz_s) _{ s_0\le s\le s_1})\big]
	=  \IS \uud \Big\|\dot \zz_s- \frac{\zz_s-\widehat\proj_S( \zz_s)}{2(s-s_0)} \Big\|_H^2  \, \kappa_s \,ds,
\end{align}
with $\widehat\proj_S$  defined at \eqref{eq-19}. 
Compare with \eqref{eq-18}. 
\\
It is easily seen (see \eqref{eq-54} below for details) that applying the parameter setting below, we arrive at \eqref{eq-22b} below.   
\begin{paraset}\  \label{hyp-02}
\begin{itemize}
\item
choose $ \kappa_s=2 (s-s_0)$,  
\item
change  time: $ s=s_0+e ^{ 2t}.$ 
\end{itemize} 
\end{paraset}
\noindent
This change of time maps $[s_0,s_1]$ onto $[- \infty, t_1]$ where $t_1:=\log \sqrt{s_1-s_0}$.
Denoting \[ \yy_t:= \zz _{ s}=\zz _{ s_0+e ^{ 2t}},\quad - \infty\le t\le t_1,\] the right-hand side of \eqref{eq-30} becomes
\begin{align}\label{eq-22b}
A^ \epsilon\big[ ( \yy_t) _{ t\le t_1})\big]:=
\IT \uud \Big\|\dot \yy_t - \yy_t + \frac{\sum _{ \sigma\in \mathcal{S}}  \xx  ^\sigma\exp \big( - {\|\yy_t -\xx  ^\sigma\|^2 }/ (2\epsilon e ^{ 2t})\big)}{\sum _{ \sigma\in \mathcal{S}}  \exp \big( - {\|\yy_t -\xx  ^\sigma\|^2 }/ (2\epsilon e ^{ 2t})\big)}\Big\|_H^2 \,dt,
\end{align}
and  the $ \Gamma$-limit \eqref{eq-24} becomes
 \begin{align}\label{eq-22}
 \Gamma\textrm{-}\lim_{ \epsilon\to 0} A^ \epsilon\big[ ( \yy_t) _{  t\le t_1})\big]=
\IT \uud \|\dot \yy_t-\yy_t+\widehat\proj_S(\yy_t)\|_H^2 \,dt.
\end{align}
 As desired, it   is  \MAG's action \eqref{eq-15}.
Remark that the appearance of the hat upon $\proj_S$ is a consequence of the $\Gamma$-limit.
 
 \begin{definition}[Action of $ \epsilon$-\MAG\ pushed by $\{x_1,\dots,x_k\}$] \label{def-04}\  \\
 The functional $A^ \epsilon$ defined at \eqref{eq-22b} is  the action of  $ \epsilon$-\MAG\ pushed by $\{x_1,\dots,x_k\}$.\\ Its  $s$-version is $\widetilde A^ \epsilon$ defined by \eqref{eq-30}.
 \end{definition}

\subsubsection*{Remarks} 
\begin{enumerate}[(a)]
\item
In \cite{ABB20}, the $k$-mapping valued stochastic process   \eqref{eq-21}   is  interpreted as \emph{surfing the heat wave}, pointing out some analogy with de Broglie's pilot wave theory in quantum physics. 
\item
Let us quote a sentence from \cite{ABB20}: \emph{``Unexpectedly, the action \eqref{eq-22} is exactly the one previously suggested by the
third author \emph{(Brenier)} in \cite{Bre11}  to include dissipative phenomena (such as sticky collisions in one space dimension) in the Monge-Ampère gravitational model!''}  The physical interpretation of the \emph{intriguing} particle system \eqref{eq-21}  is also  unclear to us at first sight, while its connection with \MAG\ seems  tight enough to think that it might not be incidental. 
\end{enumerate}

The significant feature of the stochastic differential equation \eqref{eq-21} is that the  \emph{forward velocity} of $ \mathsf{Z} ^{ \epsilon,\eta}$: the drift field $\vce$, is the \emph{current velocity} of someone else, namely $\Xe$.  \emph{How to give a meaning to this substitution?} 
One purpose of this article is to propose at Section \ref{sec-B-k-fluid} a way to obtain a clearer physical meaning. 

In order to proceed in this direction, we need to extend \MAG\ from mappings $\yy$ to fluids.

\section{\MAG\ for a fluid}
\label{sec-MAG-fluid}

We look at  the evolution of a self-gravitating  \emph{fluid} in $\Rd$ governed by a \MAG\ force field. This section only deals with \MAG: we drop $ \epsilon$-\MAG\ for a while.
In particular the state space is $\Rd,$ not $\Rdk,$ and we are back to the notation of Section  \ref{sec-MAG-mapping}.

Section \ref{sec-MAG-mapping} was dedicated to the flow of mappings $t\mapsto \yy_t=(\yy_t(x)) _{ x\in D},$ keeping track of the source element $x\in D.$ This was necessary to obtain a representation in terms of optimal transport. However,  particles in a fluid being indistinguishable, we do not observe the detail of the mapping $\yy_t,$ but only the profile
\begin{align}\label{eq-51}
\mu_t:=(\yy_t)\pf \lambda \in \mathrm{P}(\Rd)
\end{align}
 of positions at time $t$.
The following term  is part of  the integrand 
 of the action \eqref{eq-15}:
 \begin{align*}
 \|\dot \yy_t- \{\yy_t-\widehat\proj_S( \yy_t)\}\|_H^2 
 	= \int_D |\dot\yy_t(x)- \{\yy_t(x)-[\widehat\proj_S( \yy_t)](x)\}|^2\, \lambda(dx).
 \end{align*}
 
 \begin{properties}\label{hyp-01}
 Suppose that there exist a vector field $v_t(y)$   and a map $\widehat T( \mu,y)$ where $y\in\Rd$ and $\mu\in\PZ$ such that   for any $x$ in $D$,
 \begin{align}
\bullet\quad& \dot\yy_t(x)=v_t\big(\yy_t(x)\big)\in\Rd, \label{eq-35a}\\
\bullet\quad& [\widehat\proj_S( \yy_t)](x)=\widehat T\big({ \mu_t},\yy_t(x)\big)\in \Rd\quad \textrm{where}\quad \mu_t:=(\yy_t)\pf \lambda. \label{eq-35b}
 \end{align}
 \end{properties}
 Then,  \eqref{eq-15} writes as
 \begin{align}\label{eq-36}
 \begin{split}
   \IT\uud\|\dot \yy_t- 
  \{ \yy_t &-\widehat\proj_S( \yy_t)\}\|_H^2  \,dt\\
	&=\int _{ D\times [t_0,t_1]}\uud |v_t\big(\yy_t(x)\big)-\{\yy_t(x)-\widehat T \big({ \mu_t},\yy_t(x)\big)\}|^2\,   \lambda(dx)dt\\
	&=\IRT \uud |v_t(y)-[y-\widehat T ({ \mu_t},y)]|^2\, \mu_t(dy) dt\\
	&=\IT\uud \left\langle |v_t-(\Id-\widehat T_ {\mu_t})|^2, \mu_t \right\rangle \, dt
	=\IT\uud \|v_t-(\Id-\widehat T_ {\mu_t})\|^2 _{\mu_t}\,dt,
 \end{split}
 \end{align}
 where last but one equality is simply a change of notation: $ \left\langle f, \mu \right\rangle :=\IZ f \, d \mu$, and $\|\sbt\| _{ \mu}$ is a shorthand for $\|\sbt\| _{  L^2 _{ \Rd}( \mu)}$.
By \eqref{eq-35a}, $v_t$ is the velocity field of the fluid with density $ \mu_t$. Hence,   the continuity equation 
 \begin{align}\label{eq-26}
 \partial_t \mu_t+\nabla \scal (\mu_t v_t)=0
 \end{align}
 is satisfied in a weak sense. Although  this notation suggests that $ \mu_t$ should be absolutely continuous, no such hypothesis is done at Proposition \ref{res-01} below: the weak formulation \eqref{eq-27}   is valid for any probability measure $ \mu_t$. 
 
Dropping  the requirement \eqref{eq-51} that the density profile is represented by $ \mu_t= (\yy_t)\pf \lambda$ for some flow of mappings $(\yy_t) _{t_0\le  t\le t_1},$ the identity \eqref{eq-36} suggests the following extension to a fluid of the definition of Monge-Ampère gravitation.

 \begin{definition}[\MAG's action for a fluid pushed by $ \lambda$]\label{def-01}
 Let $ \lambda$ be any probability measure on $D\subset\Rd$.
  A  least action principle for a  \MAG\  self-gravitating fluid pushed by $ \lambda$ is
 \begin{align}\label{eq-28}
 \inf _{ ( \mu,v)}  \IT\uud \|v_t-(\Id-\widehat T_ {\mu_t})\|^2 _{ \mu_t}\,dt
 \end{align}
 where
 \begin{enumerate}[(i)]
 \item
  the infimum runs through all   $( \mu,v)$ satisfying the continuity equation \eqref{eq-26} in the weak sense \eqref{eq-27}, and such that the endpoint marginals $ \mu _{t_0}$ and $ \mu _{ t_1}$ are prescribed,
  \item
 $\widehat T_ \mu$ is the optimal map $\Tb_\mu$ transporting $ \mu$ to $ \lambda,$ if the Monge transport problem admits a unique solution, or some extension of it if uniqueness fails, see Remark \ref{rem-03}-(c) below.
 \end{enumerate}
 \end{definition}
 This model depends on the choices of the measure $ \lambda$ and  the extension $\widehat T$ of $\Tb.$ 
 
 \begin{remarks}\ \label{rem-03}  Clearly, its physical adequateness depends on the fulfillment of  Properties \ref{hyp-01}.
 Let us comment on them. 
 \begin{enumerate}[\ (a)]
 \item
 As far as one is interested in the evolution of the density $ \mu_t$ forgetting  the details of the mapping $\yy_t,$ property \eqref{eq-35a} is justified by Proposition \ref{res-01}.
 \item 
 Let us have a look at property \eqref{eq-35b}. 
As long as  $\mu_t$  remains absolutely continuous, for any $x\in D,$   the $x$-th component $[\widehat\proj_S( \yy_t)](x)\in\Rd$  of    $\widehat\proj_S( \yy_t)$  only depends on $\mu_t$ and the position $\yy_t(x)$ of the $x$-th particle in $\Rd.$ Indeed, in this situation $\Proj_S(\yy_t)$ contains the single element $\proj_S(\yy_t)=\Tb_{\mu_t}\circ \yy_t$ where by \eqref{eq-10}
\begin{align}\label{eq-55}
\Tb_{\mu_t}=\nabla \theta_t
\end{align}
is the optimal transport map from the absolutely continuous measure $ \mu_t$ to the  measure $ \lambda,$ in which case $[\widehat\proj_S(\yy_t)](x)=[\proj_S(\yy_t)](x)=\Tb_{\mu_t}( \yy_t(x)).$ 
\item
We think that it is \emph{physically reasonable}  to assume that  property \eqref{eq-35b} still holds when $ \mu_t$ fails to be absolutely continuous. A privileged model  should consist of replacing $\Tb_{\mu_t}(y)$ by the orthogonal projection in $\Rd$ of $y$  on the closed convex hull of the subset   $\cup _{ \pi _{ \mu_t}} \supp \pi _{ \mu_t}(\sbt\mid y)$, where the union runs through the collection $\{\pi_ {\mu_t}\}$ of all solutions of the Monge-Kantorovich  optimal transport problem \emph{from $ \mu_t$ to $ \lambda,$} and $ \supp \pi _{ \mu_t}(\sbt\mid y)$ stands for the support of the optimal plan $\pi_{ \mu_t}$ conditioned by the knowledge of the source location $y.$ Recall Definition \ref{def-03} and see Appendix \ref{app-c}.
 \end{enumerate}
\end{remarks}

\subsection*{From $\inf _{ (\mu,u)}$ to $\inf _{ (\mu)}$}

In view of the Benamou-Brenier formula \eqref{eq-BB} stated at Appendix \ref{app-OW}, for a standard  Lagrangian $  (t, \mu, \dot \mu)\mapsto \alpha_t\|\dot \mu\|^2_ \mu/2- \mathcal{U}_t( \mu),$ where $ \alpha:[t_0,t_1]\to [0, \infty)$ is a nonnegative function, we obtain
\begin{align*}
\inf _{ ( \mu,u)} \IT  \Big(\IZ \uud |u_t|^2 \alpha_t\, d \mu_t - \mathcal{U}_t( \mu_t)\Big)\,dt
	= \inf _{ ( \mu)}\IT \big(\uud \|\dot \mu_t\|^2 _{ \mu_t} \alpha_t-\mathcal{U}_t( \mu_t)\big)\,dt,
\end{align*}
where $\|\dot \mu_t\|_{ \mu_t} $ is the norm of the tangent vector $\dot\mu_t$ in the \OW-manifold. 
We shall also be concerned by Lagrangian of type $  (t, \mu, \dot \mu)\mapsto \uud \|\dot \mu- \nabla w_t\|^2 _{ \mu}\alpha_t$  where, again, $ \alpha$ is a nonnegative function. 

\begin{lemma}\label{res-11}
Suppose that $w:[t_0,t_1]\times\Rd\to\RR$ is a regular function, then
\begin{align*}
 \IT\uud& \|\dot \mu_t- \nabla w_t\|^2 _{ \mu_t}\alpha_t\, dt\\&=
 \IT \left\{\uud\|\dot \mu_t\|^2 _{ \mu_t}\alpha_t+\uud\|\nabla w_t\|^2 _{ \mu_t}\alpha_t +  \IZ \partial_t( \alpha_t  w_t)\, d \mu_t\right\} \, dt\\
&\hskip 5cm -\IZ \alpha _{ t_1} w _{ t_1}\, d \mu _{ t_1}+\IZ \alpha _{ t_0} w _{ t_0}\, d \mu _{ t_0}.
\end{align*}
\end{lemma}

\begin{proof}
Expanding the square gives us
\begin{align*}
\|\dot \mu- \nabla w_t\|^2 _{ \mu}=\|\dot \mu\|^2_ \mu+\|\nabla w_t\|^2_ \mu-2(\nabla w_t,\dot \mu)_ \mu
	=\|\dot \mu\|^2_ \mu+\|\nabla w_t\|^2_ \mu- 2(\grOW\mu  \mathcal{W}_t,\dot\mu)_\mu
\end{align*}
where $\grOW{}$ is the gradient with respect to the \OW-metric derived from 
\eqref{eq-49}, and $ \mathcal{W}_t( \mu):= \IZ w_t\, d \mu$,  see Lemma \ref{res-20}. The result follows from
\begin{align*}
\frac{d}{dt}( \alpha _{ t} \mathcal{W} _{ t}( \mu _{ t}) )= \alpha_t (\grOW{\mu_t}  \mathcal{W}_t,\dot\mu_t)_{\mu_t}
		+  \IZ \partial_t( \alpha_t w_t)\,d \mu_t.
\end{align*}
See \cite{AGS05} for a proof of this chain rule.
\end{proof}
Therefore, the Lagrangian $\uud \|\dot \mu- \nabla w_t\|^2 _{ \mu} \alpha_t$ is equivalent to the modified Lagrangian
$	
\alpha_t\|\dot \mu\|^2_ \mu/2+\alpha_t\|\nabla w_t\|^2_ \mu/2 +  \IZ \partial_t( \alpha_t  w_t)\, d \mu
$	
which has the form $\alpha_t \|\dot\mu\|_\mu^2/2- \mathcal{U}_t( \mu).$ Hence,
\begin{align}\label{eq-50}
\inf _{ ( \mu,u)} \IT  \Big(\IZ \uud |u_t-\nabla w_t|^2\, d \mu_t\Big)\alpha_t\,dt
	= \inf _{ ( \mu)}\IT \uud \|\dot \mu_t-\nabla w_t\|^2 _{ \mu_t}\,\alpha_t\,dt.
\end{align}
In particular, with some { minor additional work}, we obtain the following

\begin{proposition}[\MAG's action for a fluid pushed by $ \lambda$] \label{res-12}
If \eqref{eq-55} extends to 
\begin{align*}
\widehat T_ {\mu_t}\in \partial \theta_t,
\end{align*}
 meaning that the extension   $\widehat T_t$ of $\Tb$ is still a subgradient of a convex function  $ \theta_t$ for almost all $t,$ then the least action principle \eqref{eq-28} of Definition \ref{def-01} reads as
\begin{align}\label{eq-28b}
 \inf _{ ( \mu)}  \IT\uud \|\dot \mu_t-(\Id-\widehat T_ {\mu_t})\|^2 _{ \mu_t}\,dt,
\end{align}
 where
   the infimum runs through all   $( \mu)$  such that the endpoint marginals $ \mu _{ t_0}$ and $ \mu _{ t_1}$ are prescribed.
\end{proposition}

\section{$ \epsilon$-\MAG\ for a $k$-fluid.  Action functional}
\label{sec-fluid}

The main difficulty of  the least action principle \eqref{eq-28b} is to handle  the optimal transport  term $\widehat  T_{ \mu}$. To do so, we keep the idea of  \cite{ABB20} of working with  the $k$-mapping-valued stochastic process  $ \Xe $, see \eqref{eq-23}-\eqref{eq-29}, 
 because the symmetrization operator  $(k!) ^{ -1}\sum _{ \sigma\in \mathcal{S}}$ together with the Laplace principle when $ \epsilon$ tends to zero, is a good way to recover $\proj_S$ and therefore optimal transport, see \eqref{eq-24}.  
But, this is at the price of replacing the diffuse source measure $ \lambda$ defined at  \eqref{eq-07}, by its discrete analogue $\lk$ defined at \eqref{eq-25}.

On the other hand, we leave apart the enigmatic process $ \mathsf{Z} ^{ \epsilon,\eta}$ defined at \eqref{eq-21}. It will be replaced at Section \ref{sec-nbm}  by some empirical process   built on $\Xe$, see \eqref{eq-48}.

From now on, we only  consider  $ \epsilon$-approximations of \MAG, in the sense that we do not let $ \epsilon$ tend to zero, leaving open the problem of the fulfillment of property \eqref{eq-35b}.

\begin{definitions}[$k$-mapping and $k$-fluid]\ 
\begin{itemize}
\item
Any element $\zz$ of $\Rdk$ is interpreted as a \emph{$k$-mapping}, i.e.\ $\zz: \left\{x_1,\dots,x_k\right\} \to\Rd.$
\item
Any probability measure $p\in \PRdk $ on  $\Rdk$  is interpreted as a fluid of $k$-mappings, and it is called a \emph{$k$-fluid}, for short.
\end{itemize}
\end{definitions}

\subsection*{An $ \epsilon$-approximation of \MAG\ for  $k$-fluids}

Arguing as in  Section \ref{sec-MAG-fluid}, the $ \epsilon$-\MAG\ action functional \eqref{eq-30}: $\widetilde A^ \epsilon[\zz]= \IS \ud   \|\dot \zz_s-\vce_s( \zz_s)\|_H^2\ \kappa_s \,ds,$ admits the $k$-fluid  analogue:
\begin{align}\label{eq-33}
\IS \uud \|\dot\pp_s-\vce_s\|^2 _{ p_s}\, \kappa_s\, ds,
\end{align}
where $\vce$ is the current velocity \eqref{eq-17} of $ \Xe $,  $s\mapsto p_s\in \mathrm{P}(\Rdk )$ is  the path of density distributions    of a   $k$-fluid, $\dot\pp$ is its \emph{gradient} velocity field, meaning that the continuity equation
\begin{align}\label{eq-32}
 \partial_s p+\nabla\scal ( p\,\dot\pp)=0
\end{align}
 is satisfied in the weak sense and  $\dot\pp_s\in \gcb(p_s),$ 
 so that it is a tangent vector in the \OW-manifold, see Appendix \ref{app-OW} for the notion of \OW-manifold and Definition \ref{def-H} for the definition of $\gcb(p)$ .     
 The action \eqref{eq-33} is valid because $\vce$ is a gradient field, so that one can apply   Lemma \ref{res-11}, leading to \eqref{eq-50}.
 \\
 Similarly, at the limit $ \epsilon\to 0,$  reasoning as during the derivation of Proposition \ref{res-12}, the $k$-fluid analogue of the $k$-mapping action \eqref{eq-24}  is 
\begin{align}\label{eq-47}
\IS \uud \Big\|\dot\pp_s- \frac{\Id-\widehat T _{ p_s}}{2(s-s_0)} \Big\| _{ p_s}^2  \, \kappa_s \,ds
\end{align}
provided that the extension $\widehat T$ of $\Tb$ is  still a subgradient of some convex function.\\
In view of the $\Gamma$-limit \eqref{eq-24}, we see that \eqref{eq-33} is a reasonable $ \epsilon$-approximation of \eqref{eq-47}.

Applying  to \eqref{eq-47} the parameter setting \ref{hyp-02}: 
 $\kappa_s =2 (s-s_0)$, $ s=s_0+ e ^{ 2t}, $ as was done in order  to go from \eqref{eq-24} to \eqref{eq-22}, and setting 
 \begin{align}\label{eq-82}
 q_t:=p_s= p_{s_0+ e ^{ 2t}},
 \end{align} 
 we obtain its $t$-analogue
 \begin{align}\label{eq-52}
\ITT \uud \|\dot{\qq}_t- (\Id-\widehat T _{ q_t})\| _{ q_t}^2  \, dt, 
\end{align}
which is the \MAG\ action for $k$-fluids, analogous to the \MAG\ action functional for $k$-mappings obtained at Proposition \ref{res-12}.
Defining 
\begin{align}\label{eq-97} 
m^ \epsilon_t(d\yy):= r^ \epsilon_s(d\yy)=r^ \epsilon _{s_0+ e ^{ 2t}}(d\yy)
	\overset{ \eqref{eq-45}}= \frac{1}{k!} \sum _{ \sigma\in \mathcal{S}}  (2\pi \epsilon  e ^{ 2t}) ^{ -dk/2}
	\exp \left( - \frac{\|\yy -\xx  ^\sigma\|^2 }{2 \epsilon  e ^{ 2t}}\right) \, d\yy ,
\end{align}
we obtain
 \begin{align}\label{eq-54}
 \dot\mm^\epsilon_t(\yy)\overset{\checkmark}=2 e ^{ 2t} \vce _{s_0+ e ^{ 2t}}(\yy)
 	\overset{ \eqref{eq-17}}=\yy - \frac{\sum _{ \sigma\in \mathcal{S}}  \xx  ^\sigma\exp \big( - {\|\yy -\xx  ^\sigma\|^2 }/ (2\epsilon e ^{ 2t})\big)}{\sum _{ \sigma\in \mathcal{S}}  \exp \big( - {\|\yy -\xx  ^\sigma\|^2 }/ (2\epsilon e ^{ 2t})\big)},
 \end{align}
 where the marked equality follows from
 \begin{align*}
 \partial_t m_t= \partial_t r _{s_0+ e ^{ 2t}}=2e ^{ 2t} \partial_s r_s
 =-2 e ^{ 2t}\nabla\scal (r_s\dot\rr_s)=-2e ^{ 2t}\nabla\scal (m_t \dot\rr _{s_0+ e ^{ 2t}})
 \end{align*}
 with $s=s_0+e ^{ 2t},$ that is: $\partial_tm_t+\nabla\scal\big(m_t(2e ^{ 2t}\dot\rr _{s_0+ e ^{ 2t}})\big)=0.$
It is natural to propose the following  
  \begin{definition}[$\epsilon$-\MAG\ for a   $k$-fluid] \label{def-02}
  We consider flows $(q_t) _{ t_0\le t\le t_1}$ of $k$-fluids.
  \\
    A  least action principle for an  $\epsilon$-\MAG\  self-gravitating $k$-fluid is
 \begin{align}\label{eq-37}
 \inf _{ ( q)}  \IT \uud\|\dot\qq_t-\dot\mm^\epsilon _t\|^2 _{ q_t}\,dt
 \end{align}
  The   infimum runs through all   $( q)$   with prescribed endpoint marginals $ q _{ t_0}$ and $ q _{ t_1}$.
  \end{definition}
The action 
 \eqref{eq-37} is a reasonable $ \epsilon$-approximation of \eqref{eq-52}.
 

\subsection*{From $k$-fluids  to fluids} \label{sec-fff}

The detail of the evolution of the flow of $k$-mappings $ p_s\in \PRdk $ is necessary for  optimal transport to enter the game. However,  particles  are not coloured by any  rank of trial (the $k$ slots in $\Rd$ of a vector in $H\simeq (\Rd)^k)$. We only see a monochrome cloud. This means that instead of a probability measure $p$ on $(\Rd)^k,$ we have  to consider its projection 
\begin{align*}
\proj _{ dk\to d}(p):= k ^{ -1} \sum _{ 1\le j\le k} p_j  \in \mathrm{P}(\Rd),\qquad p\in \PRdk ,
\end{align*}
on $ \mathrm{P}(\Rd),$ where for any $1\le j\le k,$ the $j$-th marginal $p_j\in \mathrm{P}(\Rd)$ of $p$ is  defined by: $p_j(dy):=p(\Rd\times\cdots\Rd\times dy\times\Rd\times\cdots\times\Rd)$ where $dy\subset\Rd$ occupies the $j$-th slot. The weights $k ^{ -1}$ are those of $ \lk$ because
\[
\lk\overset{ \eqref{eq-25}}= k ^{ -1}\sum _{ 1\le j\le k} \delta _{ x_j}=\proj _{ dk\to d}\Big( (k!) ^{ -1}  \sum _{ \sigma\in \mathcal{S}} \delta _{ \xx  ^\sigma}\Big) \overset{ \eqref{eq-29}}= \proj _{ dk\to d} (\textrm{Law}( \mathsf{X} _{ 0})).
\]
Introducing the $j$-th projection  from $(\Rd)^k$ to $\Rd$ defined for any $(y_1,\dots,y_k)\in(\Rd)^k$ by
\[ \proj _j(y_1,\dots,y_k)=y_j\in\Rd,\]  we see that $p_j=(\proj _j)\pf p$.
On the other hand, the continuity equation  $ \partial_s p+\nabla\cdot(p\vv)=0$ in $ \Rdk ,$ implies the continuity equation
\begin{align}\label{eq-39}
\partial_s \mu+\nabla\scal( \mu v)=0
\end{align}
in $\Rd$, via the transformation  
\begin{align*}
\proj _{ dk\to d}(p,\vv)= ( \mu,v) \quad \textrm{where}\quad
\left\{ 
\begin{array}{l}
\displaystyle{\mu=\proj _{ dk\to d}(p)=k ^{ -1}\sum _{ 1\le j\le k}(\proj _j)\pf p},\\ 
 \displaystyle{ v(y)= k ^{ -1}\sum _{ 1\le j\le k} E_ p(\proj _j(\vv)\mid \proj _j=y)},\ y\in\Rd.
\end{array}
\right.
\end{align*}
Indeed, we see with \eqref{eq-27} that $ \partial_s p+\nabla\cdot(p\vv)=0$ means that for any $s_0\le s_*\le s_1$ and any function $f$ in $C^1_c(\Rdk),$ we have: 
$\int _{ \Rdk} f\, d(p _{ s_*}-p _{ s_0})
	=\int _{ \Rdk\times [s_0,s_*]} \nabla f(\yy)\cdot \vv_s(\yy)\, p_s(d\yy)ds .$
Applying this identity with  $f(y_1,\dots,y_k)=k ^{ -1}\sum _{ 1\le j\le k} g(y_j)$ for any function $g$ in $C^1(\Rd),$ gives the announced result.

Remark that if $\proj _{ dk\to d}(p)= \mu,$ then $ \proj _{ dk\to d}(p,\dot\pp)=( \mu,\dot \mu).$ This holds because $k ^{ -1}\sum _{ 1\le j\le k} E_ p(\proj_j(\dot\pp)\mid \proj_j=\sbt)$ is a gradient field on $\Rd$ since $\dot\pp$ is a gradient field on $\Rdk.$
As a consequence, the least action problem  in $ \PRdk $  based on \eqref{eq-33}:
\begin{align}\label{eq-38}
\inf _{ (p)} \IS \uud   \| \dot\pp_s-\vce_s\|^2_{p_s}  \,\kappa_s
	\, ds
\end{align}
where the infimum runs through all the $(p)$  with prescribed marginal measures  $p _{ s_0} $ and $ p _{ s_1},$ leads to the following main definition of this article.
\\
Let us denote 
\begin{align*}
\OP:=C([t_0,t_1], \PZ) \quad \textrm{or}\quad \OP:=C([s_0,s_1], \PZ)
\end{align*}
(depending on the $s$ or $t$ context), the set of $\PZ$-valued trajectories, and similarly
\begin{align*}
\OkP:=C([t_0,t_1], \PRdk) \quad \textrm{or}\quad \OkP:=C([s_0,s_1], \PRdk),
\end{align*}
the set of $\PRdk$-valued trajectories.

\begin{definition}[Least action principle for a fluid  driven by $ \epsilon$-\MAG\ pushed by $\{x_1,\dots,x_k\}$] \label{def-05}\ \\
It is the least action principle in $ \mathrm{P}(\Rd)$
\begin{align}\label{eq-40-t}
\inf _{ ( \nu)} \inf \left\{  \IT \ \uud   \| \dot\qq_t-\dot\mm^\epsilon _t\|^2_{q_t}\, dt; (q)\in \OkP : \proj _{ kd\to d}(q_t)= \nu_t, \forall t\right\}
\end{align}
where the leftmost infimum runs through all the $( \nu)\in \OP $    such that $ \nu _{ t_0}$ and $\nu _{ t_1} $ are prescribed.
Its $s$-version is obtained replacing \eqref{eq-40-t} by 
\begin{align}\label{eq-40}
\inf _{ ( \mu)} \inf \left\{  \IS \ \uud   \| \dot\pp_s-\vce_s\|^2_{p_s}\,\kappa_s\, ds; (p)\in \OkP : \proj _{ kd\to d}(p_s)= \mu_s, \forall s\right\}.
\end{align}
\end{definition}

\section{$ \epsilon$-\MAG\ for a $k$-fluid. Newton equation}
\label{sec-Neweq}

In this section, we partly stay at an  informal level, applying Otto's heuristics, but we  also prove rigorous results.  Otto's heuristics  means that, while investigating  the Euler-Lagrange  equation of a least action principle in the \OW-manifold, we only consider a finite dimensional analogy, see Appendix \ref{app-OW}. This type of equation is usually referred to as a Newton equation.

\subsection*{Statement of the main result of this section}

Recall the notation of Section \ref{sec-MAG-mapping}. In particular, the set $S=\{\xx  ^\sigma; \sigma\in \mathcal{S}\}\subset \Rdk$ introduced at  \eqref{eq-20} will play some role.
Let us introduce the field on $\Rdk$ of probability measures on $S$, 
 \begin{align*} 
\pi_t^ \epsilon(\yy)
	:= \sum _{ \sigma} \pi_t^{ \epsilon, \sigma}(\yy)\, \delta _{ \xx^ \sigma}\in \mathrm{P}(S),\qquad\yy\in\Rdk,
\end{align*}
where
 \begin{align*}
\pi_t ^{ \epsilon,\sigma}(\yy)
	:= \frac{w_t ^{ \epsilon, \sigma}(\yy)}{\sum _{ \sigma'}w _t^{ \epsilon,\sigma'}(\yy)}
\qquad \textrm{with}\qquad
w_t ^{ \epsilon, \sigma}(\yy)
	:=\exp \left(- \frac{\|\yy-\xx^ \sigma\|^2}{2 \epsilon e ^{ 2t}}\right) .
\end{align*}
The main reason for introducing $ \pi_t^ \epsilon$ is the following expression of \eqref{eq-54} 
\begin{align}\label{eq-63}
 \dot\mm_t^ \epsilon(\yy)
 	= \yy - \sum _{ \sigma}\pi_t ^{ \epsilon,\sigma}(\yy)\, \xx ^{ \sigma} 
	=\yy-\int _{ S} \xx\ [\pi_t^ \epsilon(\yy)](d\xx).
\end{align} 
For any $\yy\in \Rdk$, $ \epsilon>0$ and $t$, we write
\begin{align*}
\tilde \xx ^{ \epsilon}_t(\xx,\yy)&:= \xx-\int _{ S} \xx'\ [\pi_t^ \epsilon(\yy)](d\xx'),\quad  \xx\in S\\
F ^{ \epsilon}_t(\yy)&:=( \epsilon e ^{ 2t}) ^{ -1}\, \int _{ S} \big((\yy- \xx)\scal\tilde \xx ^{ \epsilon}_t(\xx,\yy)\big)\,\tilde \xx ^{ \epsilon}_t(\xx,\yy)\, \pi^ \epsilon_t(\yy)(d\xx).
\end{align*}
Let 
\[
S_*(\yy):= \left\{ \xx_*^1(\yy),\dots,\xx_* ^{ n_*(\yy)}(\yy)\right\} =\argmin_S \|\yy-\sbt\|\subset S
\] 
be the subset of all the $n_*(\yy)$ closest points to $\yy$  in  $S$, i.e.\ the orthogonal projection in $\Rdk$ of $\yy$ on $S.$  We also introduce
\begin{align*}
\pi_*(\yy):= n_*(\yy)^{ -1} \sum _{ 1\le n\le n_*(\yy)} \delta _{ \xx_*^n(\yy)}
\end{align*}
the uniform probability measure on $S_*(\yy),$   and
\begin{align*}
A_*(\yy):=[\mathrm{Cov}(\xx^1_*\dots,\xx ^{ n_*})\,\overline\xx_*](\yy)
	=n_* ^{ -1}\sum _{ 1\le n\le n_*} [\overline\xx_*\scal(\xx_*^n-\overline\xx_*)\  (\xx_*^n-\overline\xx_*)](\yy)
\end{align*}
where $\overline\xx_*(\yy)=n_*(\yy)^{ -1} \sum _{ 1\le n\le n_*(\yy)}   \xx_*^n(\yy).$

The acceleration  $\nabla ^{ \mathrm{OW}} _{ \dot\qq_t}\dot\qq_t$ in the \OW-manifold of the path $t\mapsto q_t$ is introduced at \eqref{eq-OWacc} in  Appendix \ref{app-OW}.

\begin{theorem}[Newton equation for a $k$-fluid driven by $ \epsilon$-\MAG] \label{res-25}
Any solution $(q)$ of the least action principle \eqref{eq-37} solves the Newton equation with the acceleration field $\nabla ^{ \mathrm{OW}} _{ \dot\qq_t}\dot\qq_t$ given, for any $\yy\in\Rdk,$ by
\begin{align}\label{eq-62}
\nabla ^{ \mathrm{OW}} _{ \dot\qq_t}\dot\qq_t(\yy)
	= \dot\mm^ \epsilon_t(\yy)+ F ^{ \epsilon}_t(\yy).
\end{align}
There exist $C,a>0$  such that for all $t$ and all $\yy\in\Rdk$, 
\begin{align}\label{eq-65}
\left| F ^{ \epsilon}_t(\yy)- \epsilon ^{ -1} e ^{ -2t} A_*(\yy) \right|
 \le C (\|\yy\|+1) e ^{ -2t} \epsilon ^{ -1}\exp(-a e ^{ -2t}\|\yy\|/ \epsilon)
\underset{ \epsilon\to 0}\longrightarrow 0.
\end{align}
Moreover, there exists a negligible subset $ \mathcal{N}$ of $\Rdk$ which is 
 a finite union of  vector subspaces with codimension at least 2, such that 
 \begin{align}
 A_*(\yy)=0, \quad \textrm{for all } \yy\not\in \mathcal{N}.
 \end{align}
Therefore,
\begin{align}
\textrm{if }& \yy\not\in \mathcal{N},\quad  
&\nabla ^{ \mathrm{OW}} _{ \dot\qq_t}\dot\qq_t(\yy)
	&= \dot\mm^ \epsilon_t(\yy)&+&\ o _{ \epsilon\to 0}(1), \label{eq-64a}
\\
\textrm{if }& \yy\in \mathcal{N},\quad 
&\nabla ^{ \mathrm{OW}} _{ \dot\qq_t}\dot\qq_t(\yy)
	&= \dot\mm^ \epsilon_t(\yy) + \epsilon ^{ -1} e ^{ -2t} A_*(\yy)\hskip -1.2cm &+&\   o _{ \epsilon\to 0}(1) \qquad
	\label{eq-64b}
\end{align}
with $\sup _{ \yy\in\Rdk}\|A_*(\yy)\|\le 2r^3< \infty,$ where $r$ is the common norm of the elements of $S$.
Furthermore,
\begin{align}\label{eq-58}
\nabla ^{ \mathrm{OW}} _{ \dot\qq_t}\dot\qq_t
	=\nabla ^{ \mathrm{OW}} _{ \dot\mm^ \epsilon_t}\dot\mm^ \epsilon_t.
\end{align}
\end{theorem}

\begin{remark}\label{rem-x2}
The right hand side of the Newton equation \eqref{eq-62} does not depend on $q_t.$ It is a constant force field which is equal to $\nabla ^{ \mathrm{OW}} _{ \dot\mm^ \epsilon_t}\dot\mm^ \epsilon_t$, see \eqref{eq-58}.
\end{remark}

\begin{proof}["Proof"]
It is a direct consequence of    Lemmas \ref{res-13}, \ref{res-15}, \ref{res-14} and \ref{res-16} below.
\\
We only propose some heuristics for a proof of Lemma \ref{res-13}. On the other hand, the proofs of Lemmas \ref{res-15}, \ref{res-14} and \ref{res-16} are complete. 
\end{proof}

\begin{remarks}\ \begin{enumerate}[(a)]
\item
With the notation of  Lemma \ref{res-16}, the thin set is $ \mathcal{N}= \left\{\yy:n_*(\yy)\ge 3\right\} .$
\item
Statement \eqref{eq-64a} with \eqref{eq-63} is reminiscent to \eqref{eq-14}.
\item
Statement \eqref{eq-64b} expresses  the divergence as $ \epsilon$ tends to zero of the force field on $ \mathcal{N}.$ This strong force field on $ \mathcal{N}$ is responsible for \emph{concentration of matter}.
\item
Identity \eqref{eq-58}  is the precise meaning of \cite{ABB20}'s expression: \emph{surfing the heat wave}.
\end{enumerate}\end{remarks}

\subsection*{Proofs of the lemmas}

Let us  prove these lemmas.

\begin{lemma}\label{res-13}
Any solution $(q)$ of the least action principle \eqref{eq-37} solves the Newton equation 
\begin{align}\label{eq-57}
\nabla ^{ \mathrm{OW}} _{ \dot\qq_t}\dot\qq_t
	=2 \dot\mm^ \epsilon_t+ 4 \epsilon^2e ^{ 4t} \nabla[\mathcal{Q}(m^ \epsilon_t|\Leb)],
\end{align}
where  we introduced the quantum potential
\begin{align}\label{eq-103}
\mathcal{Q}(m|\Leb):=-\frac{\Delta \sqrt{m}}{2\sqrt{m}}.
\end{align}
\end{lemma}

\begin{proof}[Heuristics for a proof]
Defining for any $q\in \PRdk $ 
\begin{align*}
\mathcal{F}_t(q)= - 2\epsilon e ^{ 2t} \int _{ \Rdk} \log \sqrt{m^ \epsilon_t}\,dq,
\end{align*}
we see with Lemma \ref{res-20}-(a) that $\dot\mm^ \epsilon_t=\grOW q \mathcal{F}_t$. Hence
the Lagrangian $ \uud\|\dot\qq-\dot\mm^\epsilon _t\|^2 _{ q}$  of \eqref{eq-37} writes as
\begin{align*}
\uud\|\dot\qq-\grOW{}\mathcal{F}_t(q)\|^2_q,
\end{align*}
the analogue of which in  an $n$-dimensional Riemannian manifold  is 
\begin{align*}
L(t,q,v)=\uud | v- \grad f_t(q)|^2_q
\end{align*}
for some function $f$. Beware, during this analogy, $q$ belongs to $\RR^n,$ it is not a measure. 
More precisely, in  a  coordinate system  
\begin{align*}
L(t,q,v)= \uud g _{ ij} (q)\big(v^i- g ^{ ik} (q)\partial_k  f_t(q)\big)\big(v^j- g ^{ jl}(q) \partial_l f_t(q)\big)
\end{align*}
where  the metric tensor is $g=(g _{ ij}),$ its inverse is $g ^{ -1}=( g ^{ ij})$ and we use Einstein's summation convention. 
Expending the square 
\begin{align*}
L(t,q,v)=
	\uud |v|^2_q+\uud|\grad f_t(q)|^2_q-(\grad f_t(q),v)_q.
\end{align*}
As
\begin{align*}
\ddt [f_t( \omega_t)]= \partial_t f_t( \omega_t)+ df_t( \omega_t)\scal \dot \omega_t
	= \partial_t f_t( \omega_t)+ (\grad f_t( \omega_t),\dot \omega_t) _{ \omega_t},
\end{align*}
the Lagrangian $L$ is equivalent to
\begin{align*}
\uud |v|^2_q+\uud|\grad f_t(q)|^2_q+ \partial_tf_t(q),
\end{align*}
in the sense that both least action principles attached to these Lagrangians and  sharing the same prescribed endpoints admit the same solutions. Let us introduce the scalar potential 
\[U_t(q):=-\uud|\grad f_t(q)|^2_q- \partial_tf_t(q),\] 
so  that this Lagrangian has the standard form: $\uud|v|^2_q-U_t(q).$   Therefore, the Euler-Lagrange equation $\ddt [ \partial_v L(t, \omega_t,\dot \omega_t)]- \partial_q L(t, \omega_t,\dot \omega_t)=0$ writes as the Newton equation
\begin{align*}
\nabla _{\dot \omega_t}\dot \omega_t=-\grad U_t( \omega_t),
\end{align*}
where $ \nabla$ is the Levi-Civita connection of the metric tensor $g$, see \cite{JLJ98}.
\\
By analogy, the Newton equation in the \OW-manifold is
\begin{align*}
\nabla ^{ \mathrm{OW}} _{ \dot\qq_t}\dot\qq_t=-\grOW{\qq_t} \mathcal{U}_t
\end{align*}
where the analogue of $U_t$  is 
\begin{align*}
\mathcal{U}_t(q)=-\uud \int _{ \Rdk} |\grOW q \mathcal{F}_t|^2\, dq
	- \partial_t \mathcal{F}_t(q),\qquad q\in \PRdk .
\end{align*}
Let us compute $ \mathcal{U}_t(q).$ 
 As $\grOW q \mathcal{F}_t=-2 \epsilon e ^{ 2t}\nabla\log \sqrt{m_t^ \epsilon},$ we obtain
\begin{align*}
\uud \int _{ \Rdk} |\grOW q \mathcal{F}_t|^2\, dq
 =2 \epsilon^2 e ^{ 4t}\int _{ \Rdk} |\nabla \log \sqrt{m ^ \epsilon_t} |^2\, dq.
\end{align*}
Let us look at  $ \partial_t \mathcal{F}_t(q).$ We have $m^ \epsilon_t= r^ \epsilon _{ e ^{ 2t}}$ and $ \partial_s r^ \epsilon_s= \epsilon \Delta r^ \epsilon_s/2$. Hence $ \partial_t m^ \epsilon_t=2 e ^{ 2t}{ \partial_sr^ \epsilon _{ s}}_{|s=e ^{ 2t}}= \epsilon e ^{ 2t} {\Delta r^ \epsilon_s}_{|s=e ^{ 2t}}= \epsilon e ^{ 2t} \Delta m^ \epsilon _{ t}$ and 
\begin{align*}
-\partial_t \mathcal{F}_t(q)
	= 4\epsilon e ^{ 2t} \int _{ \Rdk} \log \sqrt{m^ \epsilon_t}\,dq
	+ \epsilon e ^{ 2t} \int _{ \Rdk} \partial_t \log {m^ \epsilon_t}\,dq.
\end{align*}
But
\begin{align*}
 \partial_t \log {m^ \epsilon_t}
 	&=\frac{\partial_t  m^ \epsilon_t}{ m^ \epsilon_t}
	=\epsilon e ^{ 2t} \frac{ \Delta m^ \epsilon_t}{ m^ \epsilon_t}
	=\epsilon e ^{ 2t} \big( \Delta\log m^ \epsilon_t+ |\nabla \log m^ \epsilon_t|^2\big)\\
	&= \epsilon e ^{ 2t} \big( 2\Delta\log \sqrt{m^ \epsilon_t}+4 |\nabla \log \sqrt{m^ \epsilon_t}|^2\big),
\end{align*}
where we used $ { { \Delta u}/{u}= \Delta\log u +|\nabla\log u|^2}.$ 
This implies that
\begin{align*}
\mathcal{U}_t(q)
	&= 4\epsilon e ^{ 2t} \int _{ \Rdk} \log \sqrt{m^ \epsilon_t}\,dq
	+2\epsilon^2 e ^{ 4t} \int _{ \Rdk}\big(\Delta\log \sqrt{m^ \epsilon_t}+ |\nabla \log \sqrt{m^ \epsilon_t}|^2  \big)\,dq\\
	&= 4\epsilon e ^{ 2t} \int _{ \Rdk} \log \sqrt{m^ \epsilon_t}\,dq
	+2\epsilon^2 e ^{ 4t} \int _{ \Rdk} \frac{\Delta \sqrt{m^ \epsilon_t}}{\sqrt{m^ \epsilon_t}}\,dq.
\end{align*}
Its gradient is the \emph{constant} vector field
\begin{align*}
\grOW q \mathcal{U}_t
	=  4\epsilon e ^{ 2t} \nabla \log \sqrt{m^ \epsilon_t}
	+4\epsilon^2 e ^{ 4t} \nabla \left(\frac{\Delta \sqrt{m^ \epsilon_t}}{2\sqrt{m^ \epsilon_t}}\right) 
	=-2 \dot\mm^ \epsilon_t
	-4\epsilon^2 e ^{ 4t} \nabla [\mathcal{Q}(m^ \epsilon_t|\Leb)]
\end{align*}
Finally, the Newton equation we are after is \eqref{eq-57}.
\end{proof}

A rigorous proof of this lemma could be based on a theory of Hamiltonian  equations on the \OW-manifold as in  \cite{AG08}.

\begin{lemma}\label{res-15}
Identity \eqref{eq-58} holds.
\end{lemma}

\begin{proof}
We also obtain \eqref{eq-58} 
because the trajectory $(m ^ \epsilon_t)$ trivially minimizes  the Lagrangian $ \uud\|\dot\qq-\dot\mm^\epsilon _t\|^2 _{ q}$, implying  that it solves the least action principle (with well chosen endpoints).   Remark \ref{rem-x2} tells us that its acceleration is the constant force field \eqref{eq-62}. This leads us to the announced identity.
\end{proof}

\subsection*{Notation} 
\begin{itemize}
\item
For any  $\yy\in\Rdk,$\quad  $\yy=(y_i^l) _{ 1\le i\le k, 1\le l\le d}=(\underbrace{y_1^1,\dots,y_1^d} _{ \yy_1\in\Rd};\dots;\underbrace{y_k^1,\dots,y_k^d} _{ \yy_k\in\Rd})\in(\Rd)^k.$ 

\item
For any  $1\le i\le k,$

\begin{itemize}
\item
for any regular  $f: \Rdk\to\RR,$ \quad
 $ \partial_i f= \partial _{ \yy_i}f=( \partial _{ y_i^l}f) _{ 1\le l\le d}\in\Rd;$

 \item
 for any regular  $u: \Rdk\to\Rd,$ \quad
 $ \partial_i u= \partial _{ \yy_i}u=( \partial _{ y_i^l}u^n) _{ 1\le l,n\le d}\in \textrm{M} _{ d\times d}.$
\end{itemize}

\item
For any $ \pi\in \mathrm{P}(S)$,
\begin{itemize}
\item
 for any vector valued function $u$ on $S$ (as $S$ is a finite set, $u$ is a vector),
$ \langle u, \pi\rangle := \int _{ S} u(\xx)\, \pi(d\xx)=: \langle u(\xx), \pi\rangle$ where last identity is a practical  abuse of notation which permits us to write for instance $ \langle \xx,\pi\rangle =\int _{ S} \xx\, \pi(d\xx);$ 

\item
once $\pi$ is clear from the context, for any $\xx\in S,$  we write $\tilde \xx:= \xx- \langle \xx,\pi\rangle=\xx-\int _{ S} \xx'\, \pi(d\xx')$ or more specifically $\tilde \xx(\yy)=\tilde\xx(\xx,\yy):=\xx-\langle \xx',\pi(\yy)\rangle.$ 
\end{itemize}
 
 \item
 For any $a=(a^l) _{ 1\le l\le d}, b=(b^n) _{ 1\le n\le d}\in\Rd,$ $a\otimes b$ is the  $d\times d$-matrix defined by  $a\otimes b:= (a^lb^n) _{ 1\le l,n\le d}.$ \end{itemize}

\begin{lemma}\label{res-14}
\qquad
$	
4\epsilon^2 e ^{ 4t} \nabla [\mathcal{Q}(m^ \epsilon_t|\Leb)](\yy)
=-\dot\mm_t^ \epsilon(\yy)+F ^{ \epsilon}_t(\yy).
$	
\end{lemma}
	
\begin{proof}
Since $ \epsilon$ and $t$ are fixed, we do not write them as indices. 
The leftmost equality in \eqref{eq-54} is
\begin{align*}
\dot\mm=-2 \epsilon e ^{ 2t}\nabla \log \sqrt{m}.
\end{align*}
It implies 
\begin{align*}
\mathcal{Q}(m|\Leb)
	&= -\uud \Delta\log \sqrt{m}-\uud|\nabla\log \sqrt{m}|^2
	= -\uud \nabla\scal \nabla\log \sqrt{m}-\uud|\nabla\log \sqrt{m}|^2\\
	&= \frac{1}{4 \epsilon e ^{ 2t}}\nabla\scal\dot\mm
		- \frac{1}{8 \epsilon^2e ^{ 4t}}|\dot\mm|^2,
\end{align*} 
leading to
\begin{align}\label{eq-61}
4\epsilon^2 e ^{ 4t} \nabla [\mathcal{Q}(m^ \epsilon_t|\Leb)]
	= \epsilon e ^{ 2t}\nabla(\nabla\scal \dot\mm)
		-\nabla (\|\dot\mm\|^2)/2.
\end{align}
This shows that $\yy\mapsto\dot\mm(\yy)$ is an  ingredient that we have to work with.
It will be  convenient to use its representation  \eqref{eq-63}:
\begin{align*}
 \dot\mm(\yy)
	=\yy- \langle \xx,\pi(\yy)\rangle.
\end{align*} 
 The basic block of our calculation is, for any $1\le i\le k$ and any $ \sigma\in \mathcal{S},$
 \begin{align}\label{eq-59}
  \partial_i \pi^ \sigma(\yy)=( \epsilon e ^{ 2t}) ^{ -1} \pi^ \sigma(\yy)\tilde \xx_i ^ \sigma(\yy),\quad \forall \yy\in\Rdk.
 \end{align}
Let us show it. First of all
\begin{align*}
 \partial_i w^ \sigma(\yy)
 	=-( \epsilon e ^{ 2t}) ^{ -1}(\yy_i-\xx_i^ \sigma)w^ \sigma(\yy).
\end{align*}
Hence,
\begin{align*}
  \partial_i \pi^ \sigma
  &= \frac{ \partial_i w^ \sigma}{\sum _{ \sigma'}w ^{ \sigma'}}
  	- \frac{w^ \sigma \sum _{ \sigma'} \partial_iw ^{ \sigma'}}{(\sum _{ \sigma'}w ^{ \sigma'})^2}
	= - ( \epsilon e ^{ 2t}) ^{ -1} \pi^ \sigma \Big[ \yy_i-\xx_i^ \sigma-\sum _{ \sigma'}(\yy_i-\xx_i ^{ \sigma'}) \pi ^{ \sigma'} \Big]\\
	&=  ( \epsilon e ^{ 2t}) ^{ -1} \pi^ \sigma (\xx_i^ \sigma - \langle \xx_i, \pi\rangle)
	=  ( \epsilon e ^{ 2t}) ^{ -1} \pi^ \sigma\tilde\xx_i^ \sigma,
\end{align*}
which is \eqref{eq-59}.
This implies that for any $1\le i,j\le k,$
\begin{align}\label{eq-60}
\partial_j\dot\mm_i= \delta _{ ij}\Id _{ \Rd}
	- ( \epsilon e ^{ 2t}) ^{ -1} \langle \tilde \xx_i\otimes\tilde \xx_j, \pi\rangle ,
\end{align}
because
\begin{align*}
 \partial_j\dot\mm_i(\yy)
 	&= \partial_j(\yy_i- \langle \xx_i, \pi(\yy)\rangle )
	= \delta _{ ij}\Id _{ \Rd}- \langle \xx_i\otimes \partial_j \pi(\yy)\rangle \\
	& \overset{ \eqref{eq-59}}=  \delta _{ ij}\Id _{ \Rd}-( \epsilon e ^{ 2t}) ^{ -1} \langle \xx_i\otimes \tilde\xx_j(\yy), \pi(\yy)\rangle 
	=  \delta _{ ij}\Id _{ \Rd}-( \epsilon e ^{ 2t}) ^{ -1} \langle \tilde\xx_i (\yy)\otimes \tilde\xx_j (\yy), \pi (\yy)\rangle
\end{align*}
since $ \langle \tilde \xx_j, \pi\rangle =0.$ Consequently,  
\begin{align*}
\nabla\scal\dot\mm
	=\sum _{ i,l} \partial _{ y_i^l}\dot m_i^l
	=\sum _{ i,l} \left( 1-( \epsilon e ^{ 2t}) ^{ -1} \langle (\tilde x_i^l)^2, \pi\rangle \right) 
		=kd-( \epsilon e ^{ 2t}) ^{ -1} \langle \|\tilde\xx\|^2, \pi\rangle 
\end{align*}
and for any $1\le i\le k,$
\begin{align*}
 \epsilon e ^{ 2t} \partial_i( \nabla\scal\dot\mm)
 	&=- \partial_i \langle \|\tilde\xx\|^2, \pi \rangle 
	=- \langle \partial_i \|\tilde\xx\|^2, \pi \rangle -  \langle \|\tilde\xx\|^2, \partial_i\pi \rangle \\
	&\overset{ \eqref{eq-59}}= - \langle \partial_i \|\tilde\xx\|^2 + ( \epsilon e ^{ 2t}) ^{ -1} \tilde\xx_i, \pi \rangle.
\end{align*}
As
\begin{align*}
 \partial_i\|\tilde\xx\|^2
 	= \partial_i \sum _{ j}|\tilde\xx_j|^2
	=2 \sum_j [\partial_i\tilde\xx_j ]\, \tilde\xx_j
	\overset{ \eqref{eq-59}}=-2( \epsilon e ^{2t}) ^{ -1}\sum_j \langle \tilde\xx_i\otimes\tilde\xx_j, \pi\rangle  \,\tilde\xx_j,
\end{align*}
we have
\begin{align*}
\langle  \partial_i\|\tilde\xx\|^2, \pi\rangle =0,
\end{align*}
because $ \langle\tilde \xx_j, \pi\rangle =0.$
We finally obtain
\begin{align*}
 \epsilon e ^{ 2t} \partial_i( \nabla\scal\dot\mm)
 	=-( \epsilon e ^{ 2t}) ^{ -1} \langle \|\tilde\xx\|^2\, \tilde\xx_i, \pi\rangle .
\end{align*}
On the other hand,
\begin{align*}
\partial_i(\|\dot \mm\|^2/2)
	&= \partial_i\sum_j |\dot\mm_j|^2/2
	=\sum_j (\partial_i\dot\mm_j)\, \dot\mm_j
	\overset{ \eqref{eq-60}}= \sum_j \delta _{ ij}\dot\mm_j
		- ( \epsilon e ^{ 2t}) ^{ -1}\sum_j \langle \tilde\xx_i\otimes\tilde\xx_j, \pi\rangle \, \dot\mm_j\\
	&=\dot\mm_i- ( \epsilon e ^{ 2t}) ^{ -1}\sum_j \langle \tilde\xx_i\otimes\tilde\xx_j, \pi\rangle \, \dot\mm_j.
\end{align*}
These last two identities, together with \eqref{eq-61}, lead us to
\begin{align*}
4\epsilon^2 e ^{ 4t} \nabla [\mathcal{Q}(m|\Leb)]
=-\dot\mm ^ \epsilon_t
+( \epsilon e ^{ 2t}) ^{ -1} \Big( \langle \tilde\xx\otimes\tilde\xx, \pi\rangle \,\dot\mm
			- \langle \|\tilde\xx\|^2\,\tilde\xx, \pi\rangle \Big).
\end{align*}
As $\dot\mm(\yy)=\yy- \langle \xx, \pi\rangle $ does not depend on the variable $\xx,$ which is integrated, we see that
\begin{align*}
\langle \tilde\xx\otimes\tilde\xx, \pi\rangle \,\dot\mm
			- \langle \|\tilde\xx\|^2\,\tilde\xx, \pi\rangle
	&= \langle [\tilde\xx\otimes\tilde\xx] \,\dot\mm-  \|\tilde\xx\|^2\,\tilde\xx, \pi\rangle
	= \langle [\tilde\xx\otimes\tilde\xx] \,(\dot\mm- \tilde\xx), \pi\rangle\\
	&= \langle [\tilde\xx\otimes\tilde\xx] \,(\yy- \xx), \pi\rangle
	= \big\langle \big((\yy- \xx)\scal\tilde\xx\big)\,\tilde\xx, \pi\big\rangle.
\end{align*}
This completes the proof of the lemma.
\end{proof}

Since   $( \epsilon e ^{ 2t} )^{ -1}\big\langle \{(\yy- \xx)\scal\tilde\xx\}\,\tilde\xx, \pi^ \epsilon\big\rangle$ might diverge as $ \epsilon$ tends to zero, we have a closer look at it. 
Let us denote the energy gap
\begin{align*}
c(\yy):= \min _{ \xx\in S\setminus S_*(\yy)}\|\yy-\xx\|^2/2-\min _{ \xx\in S}\|\yy-\xx\|^2/2> 0.
\end{align*}

\begin{lemma}\label{res-16}
Let $r:=\|\xx\|,$ $\xx\in S,$ be the common norm of the elements of $S$. There is some $a>0$   such that for any nonzero  $\yy\in\Rdk,$ 
\begin{enumerate}[(a)]
\item \quad
$0< a\|\yy\|\le  c(\yy) \le 2r\|\yy\|$\quad and
\item \quad
$
|F ^{ \epsilon}_t(\yy)- \epsilon ^{ -1} e ^{ -2t} A_*(\yy)|
\le  8k!r^2(\|\yy\|+r)  e ^{ -2t}  \epsilon ^{ -1} \exp(-e ^{ -2t} a \|\yy\|/ \epsilon).$
\item\quad
$\sup_\yy|A_*(\yy)|\le 2 r^3.$
\item\quad
Furthermore, 
if $n_*(\yy)=1$ or $n_*(\yy)=2,$ then $A_*(\yy)=0$.\\\null \quad But this generally fails when $n_*(\yy)\ge 3.$
\end{enumerate}
\end{lemma}

\begin{proof}
Let us prove (a). For $\yy=0,$ $S_*(0)=S.$ Hence, $c(0)=+ \infty$ because a minimum on an empty set can be set as infinite (this is coherent with the bounds to appear below). For any nonzero $\yy,$ denote $\xx_*(\yy)$ an element of $S_*(\yy)$ and $\hat\xx(\yy)$ a minimizer of $\min _{ \xx\in S\setminus S_*(\yy)}\|\yy-\xx\|^2/2$, so that 
\begin{multline*}
c(\yy)= \|\yy-\hat\xx(\yy)\|^2/2 - \|\yy-\xx_*(\yy)\|^2/2\\
=(\xx_*-\hat\xx)(\yy)\scal\yy+\|\hat\xx(\yy)\|^2/2-\|\xx_*(\yy)\|^2/2
=(\xx_*-\hat\xx)(\yy)\scal\yy
\end{multline*}
because $\|\hat\xx(\yy)\|=\|\xx_*(\yy)\|=r.$ The same computation  shows that $\xx_*$ is an orthogonal projection of $\yy$ on $S$ if and only if $\xx_*\scal \yy\ge \xx\scal \yy,$ for all $\xx\in S.$ It follows that $\yy\mapsto\xx_*(\yy)$ and $\yy\mapsto\hat\xx(\yy)$  are functions of the unit vector $\uu_\yy:=\yy/\|\yy\|.$ Therefore, for any $\yy\neq 0,$ $
0<c(\yy)=\|\yy\|\, \uu_\yy\scal (\hat\xx-\xx_*)(\uu_\yy).$ As $S$ is a finite set,  $a:=\inf _{ \uu:\|\uu\|=1} \uu\scal (\hat\xx-\xx_*)(\uu)>0.$ Finally,
\begin{align*}
0<a\|\yy\|\le c(\yy)=\uu_\yy\scal (\hat\xx-\xx_*)(\uu_\yy)\,\|\yy\|\le 2r\|\yy\|.
\end{align*}
Let us prove (b). 
 We easily see that the total variation between $\pi^ \epsilon_t(\yy)$ and $\pi_*(\yy)$ is upper bounded by
\begin{align*}
\|\pi^ \epsilon_t(\yy)-\pi_*(\yy)\| _{ \mathrm{TV}}
	\le 2k!n_*(\yy) ^{ -1}  \exp( -e ^{ -2t}c(\yy)/ \epsilon)
	\le 2k!  \exp( -e ^{ -2t}c(\yy)/ \epsilon).
\end{align*}
This implies that
\begin{align*}
\big|\big\langle \{(\yy- \xx)\scal\tilde\xx\}\,\tilde\xx, \pi^ \epsilon_t(\yy)-\pi_*(\yy)\big\rangle\big|
	\le  8k!(\|\yy\|+r)r^2 \,\exp(-e ^{ -2t}c(\yy)/ \epsilon).
\end{align*}
It remains to evaluate and bound 
\begin{align*}
A:=
\big\langle \{(\yy- \xx)\scal\tilde\xx\}\,\tilde\xx, \pi_*(\yy)\big\rangle
	= n_* ^{ -1}\sum _{ 1\le n\le n_*} \{(\yy-\xx_*^n)\scal (\xx_*^n-\overline\xx_*)\}\, (\xx_*^n-\overline\xx_*),
\end{align*}
where we  drop the explicit dependence on $\yy$.
We are going to take advantage of both invariances: $ \|\xx_*^n\|=r$ and $\|\yy-\xx_*^n\|=\ell$ for all $1\le n\le n_*.$ As for any $n,$ $\ell^2:=\|\yy-\xx_*^n\|^2=\|\yy\|^2+\|\xx_*^n\|^2-2\yy\cdot\xx_*^n=\|\yy\|^2+r^2-2\yy\cdot\xx_*^n,$ we see that $\yy\cdot\xx_*^n$ does not depend on $n$. Hence, $\yy\scal(\xx_*^n-\overline\xx_*)=0$ for all $n$, and 
\begin{multline*}
A=-n_* ^{ -1}\sum _{ 1\le n\le n_*} \{\xx_*^n\scal (\xx_*^n-\overline\xx_*)\}\, (\xx_*^n-\overline\xx_*)\\
	=-n_* ^{ -1}\sum _{ 1\le n\le n_*} \{\|\xx_*^n\|^2-\overline\xx_*\scal\xx_*^n\}\, (\xx_*^n-\overline\xx_*)
	=-n_* ^{ -1}\sum _{ 1\le n\le n_*} \{r^2-\overline\xx_*\scal\xx_*^n\}\, (\xx_*^n-\overline\xx_*)\\
	=n_* ^{ -1}\sum _{ 1\le n\le n_*} \{\overline\xx_*\scal\xx_*^n\}\, (\xx_*^n-\overline\xx_*)
	= \mathrm{Cov}(\xx^1_*\dots,\xx ^{ n_*})\,\overline\xx_*
\end{multline*}
where the last two equalities follow from $\sum_n (\xx_*^n-\overline\xx_*)=0.$ 

Let us prove (c): $A=n_* ^{ -1}\sum _{ 1\le n\le n_*} \{\overline\xx_*\scal\xx_*^n\}\, (\xx_*^n-\overline\xx_*)$ implies $\sup_\yy|A_*(\yy)|\le 2 r^3.$

Let us prove (d).

\begin{enumerate}[(i)]
\item
If $n_*=1,$ then  $A=0$ because $\xx_*=\overline\xx_*.$
\item
If $n_*=2,$ we also have $A=0$ because, denoting $\xx_*^1=a$ and $\xx_*^2=b,$ 
\begin{align*}
A&= \ud \Big( \frac{a+b}{2}\scal a\Big) \frac{a-b}{2}
+\ud \Big( \frac{a+b}{2}\scal b\Big) \frac{b-a}{2}
	\\
	&\hskip 2cm 
	= \Big((a+b)\scal (b-a)\Big) (b-a)/8= (\|b\|^2-\|a\|^2)\, (b-a)/8 =0.
\end{align*}
Last equality holds because $\|a\|=\|b\|,$ since $a$ and $b$ belong to $S.$
\item
If $n_*\ge 3,$ this is not true anymore: $A$ does not vanish in general. 
As an example, take the three vectors $a=(1,0,0,\dots,0)$, $b=(0,1,0,\dots,0)$ and $c=(-1,0,0,\dots,0)=-a$, lying on a sphere centered at zero.  Their barycenter $\overline\xx_*$ is $b/3$, and
\begin{align*}
A= 1/3\ [(a\scal b/3)\,(a-b/3)+(b\scal b/3)&\, (b-b/3)+(c\scal b/3)\,(c-b/3)]\\
	&= 1/3\ [(b\scal b/3)\, (b-b/3)]= 2b/27\neq 0,
\end{align*}
because $a\scal b=c\scal b=0$ and $\|b\|^2=1.$
\end{enumerate}
This completes the proof of the lemma.
\end{proof}

\section{Schrödinger problem}
\label{sec-nbm}

This section is dedicated to already well known  results  about the  large deviation principle for the empirical process of a collection of independent copies of diffusion processes by Dawson and Gärtner  \cite{DG87} and Föllmer \cite{Foe85}, its connection with the Schrödinger problem \cite{Sch31,Sch32,Foe85,Leo12e} and an expression of its large deviation  rate function as a Lagrangian action in the \OW-manifold derived in \cite{CCGL20}. This is a preliminary step for the construction at Section \ref{sec-B-k-fluid} of an interacting Brownian particle system whose empirical process satisfies the Gibbs conditioning principle of Statement \ref{res-10} below.

\subsection*{Our  goal}

Our ultimate goal would be to provide some physical representation differing from the enigmatic model \eqref{eq-21}. We are in search for some collection $(  \XtN) _{ N\ge 1}$    of   random elements in the set $\OkP :=C([s_0,s_1], \mathrm{P}(\Rdk ))$ of all continuous paths on the set  $ \mathrm{P}(\Rdk )$ of $k$-fluids which 
satisfies the following 
\begin{conjecture}\label{res-10} \emph{(Gibbs conditioning principle).} 
There is exists some $\XtN$ such that 
 for any probability measures $ \alpha$ and $ \beta$ on $\Rdk,$ conditionally on  $ \XtN(s_0)\simeq \alpha$ and $ \XtN(s_1)\simeq \beta,$  the most likely trajectory $p \in\OkP $ of $\XtN$ as $N$ tends to infinity solves the least action principle \eqref{eq-38}.
 \end{conjecture}
 
 \begin{remarks}\ \begin{enumerate}[(a)]
\item
The fluctuation parameter $ \epsilon>0$ is fixed once for all. We are only concerned by limits as $N$ tends to infinity.
 \item
  This statement is fuzzy: a rigorous one  should consider a $ \Gamma$-limit along decreasing neighborhoods of $ \alpha$ and $ \beta$.
 \end{enumerate}\end{remarks}
 
This article does not succeed in providing such a result. However, we shall obtain at  Section \ref{sec-kappa} a Gibbs conditioning principle which is half way in this direction, see Corollary \ref{res-10b} and Proposition \ref{res-10c}.

 \subsubsection*{Time $s$ versus time $t$}
 
In the present section, we shall stick to time $s$  and consider the least action principle \eqref{eq-38} rather than its $t$-version \eqref{eq-37}.  Once a particle system corresponds to the  least action principle  \eqref{eq-38}, it simply remains to apply the parameter setting \ref{hyp-02} to arrive at \eqref{eq-37}.

  \subsection*{A Brownian cloud  related to $\epsilon$-\MAG\ for a $k$-fluid}
  
  As a first step, we recall the large deviation principle satisfied by the empirical process
\begin{align}\label{eq-48}
X^N: s\in [s_0,s_1]\mapsto \frac 1N   \sum _{ 1\le i\le N} \delta _{ \mathsf{X} ^{ \epsilon}_i(s)}
\in \mathrm{P}\big(\Rdk \big)
\end{align}
of  a  sequence  $( \mathsf{X} ^{ \epsilon}_i) _{ i\ge 1}$ of independent copies  of $ \mathsf{X} ^{ \epsilon},$ see \eqref{eq-23} and \eqref{eq-29}, that is 
\begin{align*}
\mathrm{Law} ( \mathsf{X} ^{ \epsilon}_i) _{ i\ge 1} 
=(R^ \epsilon) ^{ \otimes \infty}
\end{align*}
where $R^ \epsilon$ is the law of the process $ \Xe $. 
\\
  The random process $X^N$ describes a  Brownian cloud of $N$ particles evolving in $\Rdk .$
 Its large deviations  in $\OkP$:
\begin{align*}
\Proba( X^N\in \sbt )
\underset{N\to \infty}\asymp \exp \Big( - N\ \inf _{ p\in \bullet} J(p)\Big),
\end{align*}
 are well known since the pioneering article \cite{DG87} by Dawson and Gärtner, and its presentation by Föllmer in \cite{Foe85}. The rate function $p\in \OkP \mapsto J(p)\in [0, \infty]$  is expressed below at Proposition \ref{res-08}.  This implies that, for any two prescribed time marginals $ \alpha$ and $ \beta$ in $\mathrm{P}(\Rdk ),$
 \begin{multline*}
 \Proba(X^N\in\sbt \mid X^N(s_0)\simeq \alpha , X^N(s_1)\simeq \beta ) \\\underset{N\to \infty}\asymp \exp \Big( - N\ \Big[\inf _{ p\in \bullet, p _{ s_0}=\alpha , p _{ s_1}=\beta } J(p)-\inf _{ p: p _{ s_0}=\alpha , p _{ s_1}=\beta } J(p)\Big]\Big) ,
 \end{multline*}
 which in turns implies the following  
\begin{statement}\label{res-09} \emph{(Gibbs conditioning principle).} 
For any probability measures $ \alpha$ and $ \beta$ on $\Rdk,$ conditionally
 on  $ X^N(s_0)\simeq \alpha$ and $ X^N(s_1)\simeq \beta,$  the most likely trajectory $p\in\OkP $ of $X^N$ as $N$ tends to infinity solves the least action principle
 \begin{align}\label{eq-31}
 \inf  J(p),\quad p\in\OkP : p _{ s_0}=\alpha , p _{ s_1}=\beta .
 \end{align}
\end{statement}
We  focus on $X^N$ 
 because   the minimization problem  \eqref{eq-31}  happens to be close to the least action principle \eqref{eq-38} we are after, see Proposition \ref{res-02} below.   
Let us give some indications about the computation of $J.$ 
The random paths $ \Xe _i$ take their values in the space 
\[\Ok:=C([s_0,s_1],\Rdk )\] 
of all continuous paths on $\Rdk .$
   The large deviation principle as $N$ tends to infinity of  their empirical measures
\[
\Xh:= \frac 1N   \sum _{ 1\le i\le N} \delta _{ \mathsf{X} ^{ \epsilon}_i}\in  \POk,
\] 
which take their values in the set $\POk$ of all probability measures on $\Ok,$ 
is given by  Sanov's theorem, see \cite{DZ} for instance, which roughly states that
\begin{align*}
\Proba( \Xh\in \sbt )
\underset{N\to \infty}\asymp \exp \Big( - N\ \inf _{ P\in \bullet}H(P|R^ \epsilon)\Big) 
\end{align*}
where 
\[
H(P|R^ \epsilon):=E_P \log(dP/dR^ \epsilon)\in[0, \infty],
\quad P\in  \POk,
\]
is the relative entropy of $P$ with respect to $R^ \epsilon.$ 
By the contraction principle, the large deviation rate function  for $X^N$ is
\begin{align}\label{eq-100}
J(p)=\inf \{ H(P|R^ \epsilon); P\in\POk: P_s=p_s, \forall s\in[s_0,s_1] \},\quad p=(p_s) _{ s_0\le s\le s_1}\in \OkP 
\end{align}
where $P_s\in \mathrm{P}(\Rdk )$ is the $s$-th marginal of $P$.

\subsection*{Schrödinger problem}
The minimization problem \eqref{eq-31} is also called the Schrödinger problem. It was  addressed and solved (at least formally)  in 1931 by Schrödinger in \cite{Sch31,Sch32}. In view of \eqref{eq-100}, its solution $(p_s)$ is the time marginal flow  of the solution $P$, i.e.
\begin{align}\label{eq-102}
p_s=P_s, \quad s_0\le s\le s_1,
\end{align}
of the following entropy minimization problem
\begin{align}\label{eq-101}
\inf H(P|R^ \epsilon), \quad P\in\POk : P _{ s_0}=\alpha , P _{ s_1}=\beta.
\end{align}
This  formulation   in terms of the relative entropy   is due to Föllmer \cite{Foe85}.
Problem \eqref{eq-101} admits at most one solution because $H(\sbt|R^ \epsilon)$ is strictly convex  and the  constraint set $ \left\{P\in\POk : P _{ s_0}=\alpha , P _{ s_1}=\beta \right\} $ is convex.

\begin{definitions}\label{def-x2}
 If it exists, it is called the \emph{Schrödinger bridge} with respect to $R^ \epsilon$ between $ \alpha$ and $ \beta.$ Its marginal flow \eqref{eq-102} is called the \emph{entropic interpolation}  with respect to $R^ \epsilon$ between $ \alpha$ and $ \beta.$ The entropy minimization problem \eqref{eq-101} is  called the \emph{Schrödinger  problem}.
\end{definitions}

\subsection*{The  function $J$} Finding the minimizer $P\in\POk$ of \eqref{eq-100} is called the Dawson-Gärtner problem \cite{DG87}.  Next proposition gives its solution together with  an expression of the large deviation rate function $J.$

\begin{proposition}\label{res-04}
If $J(p)$ is finite, this infimum is uniquely attained at some $P(p)\in\POk$ which is Markov with a gradient drift field  $\vvf s\in \gcb(p_s)$ for almost every $s,$ 
i.e.\ $P(p)$ solves the martingale problem with Markov generator $ \partial_s+\vvf{s}\scal\nabla+ \epsilon \Delta/2$. Moreover, 
\begin{align}\label{eq-70}
J(p)=H(P(p)|R^ \epsilon)
	=H(p_{ s_0}|r^ \epsilon _{ s_0})
		+\epsilon ^{ -1} \IS \uud  \| \vvf{s}\|^2_{p_s } \,ds.
\end{align}
\end{proposition}

\begin{proof}
See \cite{Foe86,Foe85} for the proof of this result, and also \cite{DG87} for the original  proof leading to an alternate equivalent expression of $J(p)$.
\end{proof}

We prefer departing from the standard representation \eqref{eq-70} of $J$ which was  put forward in \cite{DG87,Foe85}, to exploit the following alternate representation.

\begin{proposition}\label{res-08}
The large deviation rate function $J$ is given for any $p\in \OkP $ by
\begin{align}\label{eq-34}
J(p)=\uud H(p _{ s_0}|r^ \epsilon _{ s_0})+&\uud H(p _{ s_1}|r^ \epsilon _{ s_1})
	+\epsilon ^{ -1} \IS\uud    \| \dot\pp_s-\vce_s\|^2_{p_s}  \,ds
	+ \epsilon \IS I(p_s| r^ \epsilon_s)\, ds.
\end{align}
 \end{proposition}
 \begin{proof}
 For a proof of the identity \eqref{eq-34} which relies on time reversal, see \cite[\S 6]{CCGL20}.
\end{proof}

In  formula \eqref{eq-34},
\begin{align*}
I(p|r):= \uud  \Big\|\nabla\log \sqrt{ \frac{dp}{dr}}\Big\|^2_p
\end{align*}
is the Fisher information of $p$ with respect to $r,$ $r ^{ \epsilon}_s$ is the $s$-th marginal of $R^ \epsilon$, 
 and 
 \begin{align*}
  \dot\pp_s=\vvf{s}- \epsilon\nabla\log \sqrt{p_s}=:\vvc{P(p)}_s
 \end{align*} is the current velocity of $P(p)$. In particular, $\dot\pp_s$ belongs to  $\gcb(p_s)$ for almost every $s$ and it satisfies the continuity equation \eqref{eq-32}    in the weak sense. It is part of this result that $dp/dr$ is regular enough for its derivative to exist in the sense of distributions as soon as $J(p)$ is finite. \\
As a consequence of this proposition, we obtain

\begin{proposition}\label{res-02}
 The least action principle \eqref{eq-31} is equivalent to
\begin{align}\label{eq-41}
\inf _{ p} \IS \Big(\uud   \| \dot\pp_s-\vce_s\|^2_{p_s}  
	+ \epsilon^2 I(p_s| r^ \epsilon_s)\Big)\, ds
\end{align}
where the infimum runs through all the $p\in\OkP$ satisfying $p _{ s_0}=\alpha $ and $ p _{ s_1}=\beta .$ 
\end{proposition}

Hence, a restatement of Statement \ref{res-09} is 
\begin{statement}\emph{(Gibbs conditioning principle).}  
For any probability measures $ \alpha$ and $ \beta$ on $\Rdk,$ conditionally on  $ X^N(s_0)\simeq \alpha$ and $ X^N(s_1)\simeq \beta,$  the most likely trajectory $p\in\OkP $ of $X^N$ as $N$ tends to infinity solves the least action principle \eqref{eq-41}.
 \end{statement}

Comparing \eqref{eq-41} with \eqref{eq-38}: $ \displaystyle{
\inf _{ (p)} \IS \uud   \| \dot\pp_s-\vce_s\|^2_{p_s}  \,\kappa_s
	\, ds}$, we see that it is necessary 
\begin{itemize}
\item
to introduce the coefficient $\kappa_s$ into \eqref{eq-41} and 
\item
to remove from \eqref{eq-41} the Fisher information term.
\end{itemize}

\section{Plugging  $\kappa_s$ in}
\label{sec-kappa}

Our aim is to introduce the coefficient $\kappa_s$ into \eqref{eq-41}. The main result of this section is 

\begin{theorem}\label{res-17}
Let $Z$ be a continuous  Markov process taking its values in $\Rd$ on the time interval $[s_0,s_1].$  Let $(Z''_i; i\ge 1)$ be  an iid sequence of copies of $Z$ and denote $Z^N_s:=N ^{ -1}\sum _{ i=1}^N \delta _{ Z''_i(s)},$ $s_0\le s\le s_1,$  the corresponding empirical process. We assume that   $(Z^N) _{ N\ge 1}$ obeys the large deviation principle:
$	
\Proba(Z^N \in \sbt)
\underset{N\to \infty}\asymp \exp \Big( - N\ \inf _{ \zeta\in \bullet}I(\zeta)\Big),
$	
in    $ C([s_0,s_1],\mathrm{P}(\Rd))$ with rate function    
\begin{align}\label{eq-77}
I(\zeta)=I _{ s_0}(\zeta _{ s_0})+\int _{ s_0} ^{ s_1} \mathcal{L}(s,\zeta_s,\dot\zeta_s)\,ds,\quad \zeta=(\zeta_s) _{ s_0\le s\le s_1},
\end{align}
for some Lagrangian function $ \mathcal{L}$ on the tangent bundle of the \OW-manifold, see  \eqref{eq-78}. 
\\
Let $s\in[s_0,s_1]\mapsto \kappa_s\in (0, \infty)$ be a continuously differentiable positive function.
\\
Let  $(\Zb^N) _{ N\ge 1}$ be the sequence of modified empirical processes described at page \pageref{eq-69}, see \eqref{eq-69}, \eqref{eq-76}, \eqref{eq-68b} and \eqref{eq-69b}. 
\\
Then,  $(\Zb^N) _{ N\ge 1}$ 
obeys the large deviation principle: 
\begin{align*}
\Proba(\Zb^N \in \sbt)
\underset{N\to \infty}\asymp \exp \Big( - N\ \inf _{ \zeta\in \bullet}I^\kappa(\zeta)\Big),
\end{align*}
in    $ C([s_0,s_1],\mathrm{P}(\Rd))$ with rate function
\begin{align}\label{eq-66}
I^\kappa(\zeta)=\kappa _{ s_0}I _{ s_0}(\zeta _{ s_0})+\int _{ s_0} ^{ s_1} \kappa_s \mathcal{L}(s,\zeta_s,\dot\zeta_s)\,ds,\quad \zeta=(\zeta_s) _{ s_0\le s\le s_1}.
\end{align}
\end{theorem}

\begin{proof}[Sketch of proof]
Approximate $\Zb^N$ by a sequence  $\big((\Zt ^{ R,L,N} )_{ N\ge 1}\big) _{ R,L\ge 1}$ that is constructed as in the first model to be described below at page \pageref{sec-fm}, see \eqref{eq-68} and \eqref{eq-80}. Then, relying on \eqref{eq-79}, apply the abstract theorem  \cite[Thm.\,4.2.16]{DZ} on exponentially good approximations.
\end{proof}
We only give a  {  sketch of the proof} to save time. However, the main arguments are   exposed  in the next pages.

\subsection*{Gibbs conditioning principle}

Let us denote $\Xb^N$ the process $\Zb^N$ where $Z=\Xe$ defined at \eqref{eq-23}.
At each time $s,$ $\Xb^N_s$ is the empirical measure of a cloud of $\lfloor \kappa_s N\rfloor$  $k$-mappings. Each of these particles ($k$-mappings) is   a copy of $\Xe.$  But these copies are not independent. Assuming that  $\kappa$ is an increasing function, during any small time interval $[s,s+h],$ a fraction $\kappa'_s\,h+o _{ h\to 0}(h)$ of the particles splits: each of them gives birth to a new particle starting at the same place as its genitor and evolving in the future  according to the kinematics of $\Xe$ and independently of the other particles. 
\\
Remark that although the number $\lfloor \kappa_s N\rfloor$ of particles in the cloud increases with time, $\Xb^N$ is normalized so that its total mass remains constant: $\Xb^N_s(\Rdk)=1$ for all $s.$ As a consequence, the random fluctuation of the whole cloud decreases.

Theorem \ref{res-17} and Proposition \ref{res-08} tell us that $(\Xb^N) _{ N\ge 1}$ obeys the large deviation principle
\begin{align*}
\Proba(\Xb^N \in \sbt)
\underset{N\to \infty}\asymp \exp \Big( - N\ \inf _{ p\in \bullet}J^\kappa(p)\Big),
\end{align*}
in    $ \OkP$ with rate function given for any $p\in\OkP$ by
\begin{equation} \label{eq-81}
\begin{split}
J^\kappa(p)=\uud \kappa _{ s_0} H(p _{ s_0}|r^ \epsilon _{ s_0})&+\uud\kappa _{ s_1} H(p _{ s_1}|r^ \epsilon _{ s_1})\\
	&+\epsilon ^{ -1} \IS\uud    \| \dot\pp_s-\vce_s\|^2_{p_s}\kappa_s  \,ds
	+ \epsilon \IS I(p_s| r^ \epsilon_s)\kappa_s\, ds.
\end{split}
\end{equation}
Note that, because of the $s_1$-term,  the rate function $J$ of Proposition \ref{res-08} has not exactly the form of $I$ at Theorem \ref{res-17}. To obtain the announced result, apply Theorem \ref{res-17} with Proposition \ref{res-04}, and also with its time reversed analogue which shares the same rate function since time reversal is one-one. Then take the half sum of these two expressions to arrive at \eqref{eq-81}, see \cite[\S\,6]{CCGL20} for details. 
It follows that an action functional attached to $(\Xb^N) _{ N\ge 1}$ is
\begin{align*}
 p\in\OkP\mapsto \IS\uud    \| \dot\pp_s-\vce_s\|^2_{p_s}\kappa_s  \,ds
	+ \epsilon^2 \IS I(p_s| r^ \epsilon_s)\kappa_s\, ds.
\end{align*}

\begin{corollary}\label{res-10b}   For any probability measures $ \alpha$ and $ \beta$ on $\Rdk,$ conditionally on  $ \Xb^N(s_0)\simeq \alpha$ and $ \Xb^N(s_1)\simeq \beta,$  the most likely trajectory $p \in\OkP $ of $\Xb^N$ as $N$ tends to infinity solves the least action principle 
\begin{align*}
\inf _{ p} \IS \Big(\uud   \| \dot\pp_s-\vce_s\|^2_{p_s}  
	+ \epsilon^2 I(p_s| r^ \epsilon_s)\Big)\,\kappa_s\, ds
\end{align*}
where the infimum runs through all the $p\in\OkP$ satisfying $p _{ s_0}=\alpha $ and $ p _{ s_1}=\beta .$ 
 \end{corollary}

Let us switch to time $t$ by means of the Parameter setting  \ref{hyp-02}: $\kappa_s=2(s-s_0)$ and $s=s_0+e ^{ 2t},$ and  \eqref{eq-82}: $q_t:=p_s= p_{s_0+ e ^{ 2t}},$ to obtain the $t$-action
\begin{align}\label{eq-83}
 q\in\OkP\mapsto \IT\uud    \| \dot\qq_t-\dot\mm^\epsilon _t\|^2_{q_t}\,dt
	+   \IT 4 e ^{ 4t} \epsilon^2 I(q_t| m^ \epsilon_t)\,  dt,
\end{align}
whose first term is \eqref{eq-37} as desired.

Defining $\Yb^N(t):=\Xb^N(s)=\Xb^N(s_0+e ^{ 2t}),$ Corollary \ref{res-10b} becomes

\begin{proposition}
\label{res-10c} \emph{(Gibbs conditioning principle).} For any probability measures $ \alpha$ and $ \beta$ on $\Rdk,$ conditionally on  $ \Yb^N(t_0)\simeq \alpha$ and $ \Yb^N(t_1)\simeq \beta,$  the most likely trajectory $q \in\OkP $ of $\Yb^N$ as $N$ tends to infinity solves the least action principle 
\begin{align}\label{eq-411}
\inf _{ q} \IT\uud    \| \dot\qq_t-\dot\mm^\epsilon _t\|^2_{q_t}\,dt
	+   \IT 4 e ^{ 4t} \epsilon^2 I(q_t| m^ \epsilon_t)\,  dt
\end{align}
where the infimum runs through all the $q\in\OkP$ satisfying $q _{ t_0}=\alpha $ and $ q _{ t_1}=\beta .$ 
\end{proposition}

\subsection*{Preliminary considerations}

Let $(Z_i; i\ge 1)$ be a sequence of independent copies of some Markov process $(Z_s) _{ s_0\le s\le s_1}$  taking its values in $\Rd$ and such that the empirical process \[Z^N_s:= N ^{ -1}\sum _{ i=1}^N \delta _{ Z_i(s)},\quad s_0\le s\le s_1,\] obeys the large deviation principle
\begin{align*}
\Proba( Z^N \in \sbt)
\underset{N\to \infty}\asymp \exp \Big( - N\ \inf _{ \zeta\in \bullet}I(\zeta)\Big) 
\end{align*}
 with the rate function  $I$  given at \eqref{eq-77}.
In the general setting of the present section, we want to  find some modification $Z ^{ \kappa,N}$ of $X^N$ which obeys the  large deviation principle
\begin{align*}
\Proba( Z ^{ \kappa,N} \in \sbt)
\underset{N\to \infty}\asymp \exp \Big( - N\ \inf _{ \zeta\in \bullet}I^\kappa(\zeta)\Big) 
\end{align*}
 with  the modified rate function $I^\kappa$ given at \eqref{eq-66}.

The key idea for this purpose is the following easy remark.   Sanov's theorem states that the empirical measure $\widehat Z^N:= N ^{ -1} \sum _{ i=1}^N \delta _{ Z_i}$ obeys the large deviation principle
\begin{align*}
\Proba(\widehat Z^N \in \sbt)
\underset{N\to \infty}\asymp \exp \Big( - N\ \inf _{ P\in \bullet}H(P|R)\Big) 
\end{align*}
where $H(P|R)=\int \log(dP/dR)\,dP$ is the relative entropy of $P$ with respect to the law $R$ of the Markov process $Z.$ 
It immediately follows that for any $ \kappa >0,$  the sequence of modified empirical measures 
\[ \widehat Z ^{ \kappa,N}:= \lfloor \kappa N\rfloor ^{ -1} \sum _{ i=1} ^{ \lfloor \kappa N\rfloor} \delta _{ Z_i}\] where $\lfloor a\rfloor$ is the integer value of $a$, obeys the large deviation principle
\begin{align*}
\Proba(\widehat Z ^{ \kappa ,N} \in \sbt)
\underset{N\to \infty}\asymp \exp \Big( -\lfloor \kappa N\rfloor\ \inf _{ P\in \bullet}H(P|R)\Big)
\simeq  \exp \Big( -N\ \inf _{ P\in \bullet}\kappa\,H(P|R)\Big)
\end{align*}
with rate function $\kappa H(\sbt|R)$ instead of $H(\sbt|R).$ It also follows with the  contraction principle that for any continuous mapping $\Phi$, we obtain
\begin{align*}
\Proba(\Phi(\widehat Z^N) \in \sbt)
\underset{N\to \infty}\asymp \exp \Big( - N\ \inf _{ \zeta\in \bullet}I(\zeta)\Big) 
\quad \textrm{where}\quad I(\zeta)=\inf _{ P:\Phi(P)=\zeta} H(P|R)
\end{align*}
and similarly
\begin{align*}
\Proba(\Phi(\widehat Z ^{ \kappa,N}) \in \sbt)
\underset{N\to \infty}\asymp \exp \Big( - N\ \inf _{ \zeta\in \bullet}I^\kappa(\zeta)\Big) 
\quad \textrm{where}\quad I^\kappa(\zeta)=\inf _{ P:\Phi(P)=\zeta}\kappa H(P|R)
\end{align*}
with rate function\begin{align}\label{eq-67}
I^\kappa=\kappa I.
\end{align}
 Since $Z^N=\Phi(\widehat Z^N)$ where $\Phi$ is the continuous application mapping a path measure $P$ to its flow of marginal measures $\Phi(P)=(P_s) _{ s_0\le s\le s_1},$ we see that if $I$ denotes the large deviation rate function as $N$ tends to infinity of the empirical process $Z^N,$ then the modified empirical process
\begin{align*}
 \lfloor \kappa N\rfloor ^{ -1} \sum _{ i=1} ^{ \lfloor \kappa N\rfloor} \delta _{ Z_i(s)},\quad s_0\le s\le s_1,
\end{align*}
obeys the large deviation principle with rate function $I^\kappa =\kappa I.$
Assume in addition that $Z$ is Markov so that the rate function $I$ is additive in time, and more precisely that it  has the form of the action functional \eqref{eq-77}.
Of course, the rate function $I^\kappa=\kappa I$ 
 of the modified empirical process is associated with the Lagrangian $\kappa \mathcal{L}$ instead of $\mathcal{L}.$ 

Considering a varying coefficient $s\mapsto\kappa_s,$ it is tempting to guess that the modified empirical process defined by
\begin{align}\label{eq-53}
 \lfloor \kappa_s N\rfloor ^{ -1} \sum _{ i=1} ^{ \lfloor \kappa_s N\rfloor} \delta _{ Z_i(s)},\quad s_0\le s\le s_1,
\end{align}
obeys the large deviation principle with  the rate function $I^\kappa$ defined at  \eqref{eq-66}.

But this does \emph{not} hold. 
 The reason for this is that during a small time interval $[s, s+ h],$ the "new" particles $j\in \{\lfloor\kappa_s N\rfloor+1,\dots, \lfloor\kappa _{ s+ h}N\rfloor\}$ to be (a): added if $\kappa_s< \kappa _{ s+ h}$, or (b): removed if $\kappa_s> \kappa _{ s+ h},$ are  sampled from the {random} collection (a):  $(Z_i) _{ i\ge\lfloor\kappa_s N\rfloor }$, or (b): $(Z_i) _{ i\le\lfloor\kappa_s N\rfloor }.$  This \emph{random} sampling is subject to large deviations so that some additional cost must be added to the  Lagrangian cost. 
\\
To fix the idea, we assume up to page \pageref{subsec-kappa} that \emph{$s\mapsto \kappa_s$ is increasing}. This is satisfied for our model ($\kappa_s=2(s-s_0)$, see the parameter setting \ref{hyp-02}). The general case will be considered later at \eqref{eq-68b} and \eqref{eq-69b}.
\\
This additional cost is the large deviation cost for observing the empirical measure of the newcomers
\[ 
(\lfloor\kappa _{ s+ h}N\rfloor-\lfloor\kappa _{ s}N\rfloor ) ^{ -1} \sum _{j\in \{\lfloor\kappa_s N\rfloor +1,\dots, \lfloor\kappa _{ s+ h}N\rfloor\} } \delta _{ Z'_j(s)}
\]
 close to the actual state $ \lfloor \kappa_s N\rfloor ^{ -1} \sum _{ i=1} ^{ \lfloor \kappa_s N\rfloor} \delta _{ Z_i(s)}\simeq \zeta_s$ which might be far from the most likely state. Therefore, the large deviations of the empirical process \eqref{eq-53} are not governed by the desired  rate function \eqref{eq-66}. 
 \\
 This remark indicates the way to a good candidate: one must replace the random sampling of the newcomers by an \emph{almost deterministic} one to control the distance between the actual empirical measure and the empirical measure of the newcomers.

 \subsection*{First model} \label{sec-fm}
 
  To do so,  we  approximate the desired empirical process by introducing the newcomers at each time 
  \[
  \sigma_\ell :=s_0+ \ell  h,\quad 1\le \ell \le L-1,
  \]
 where $ h= (s_1-s_0)/L$ for some arbitrarily large integer $L.$  Let $\Zt ^{ L,N}$ denote this approximate process. At time $ \sigma_\ell ^-,$ meaning the left limit $\lim _{ s\to \sigma_\ell , s< \sigma_\ell },$  its value is
 \begin{align*}
 \Zt ^{ L,N} _{ \sigma_\ell ^-}=(N _{ \ell -1}) ^{ -1} \sum _{ i=1} ^{N _{ \ell -1} } \delta _{ Z_i( \sigma_\ell ^-)},\quad \textrm{where}\quad N_\ell :=\lfloor \kappa _{ \sigma _{ \ell }}N\rfloor.
 \end{align*}
 \emph{Beware}: although we use the same notation $Z_i$ for our random paths, they are not independent anymore as previously.
 \\
 The  newcomers $\{Z'_j( \sigma_\ell );N _{ \ell -1}< j\le N_\ell  \}$ at time $ \sigma_\ell $ are sampled from the support of $ Z ^{ L,N} _{ \sigma_\ell ^-}$ in such a way that their empirical measure
 $(N_\ell -N _{ \ell -1})^{ -1} \sum _{ j=1} ^{N_\ell -N _{ \ell -1}} \delta _{ Z'_j( \sigma_\ell )}$ is arbitrarily close, as $N$ tends to infinity, to $ Z ^{ L,N} _{ \sigma_\ell ^-}.$   This is possible, as explained at next subsection,  see \eqref{eq-71} and \eqref{eq-72}.
 At time $ \sigma_\ell ,$ we state
\begin{align}\label{eq-68}
 \Zt ^{ L,N} _{ \sigma_\ell }= N_\ell  ^{ -1}\Big( N _{ \ell -1} \Zt ^{ L,N} _{ \sigma_\ell ^-}
 	+  \sum _{ j=1} ^{N_\ell -N _{ \ell -1}} \delta _{ Z'_j( \sigma_\ell )}\Big)
	=N _{ \ell }^{ -1} \sum _{ i=1} ^{N _{ \ell } } \delta _{ Z_i( \sigma_\ell )}
\end{align}
where last equality results from some relabeling. 
Keeping track of the history of these relabeling procedures, one sees that \emph{at each time $ \sigma_\ell $ some particles split in such a way that for each $\ell $, $ \Zt ^{ L,N} _{ \sigma_\ell }$ is arbitrarily close, as $N$ tends to infinity, to $ \Zt ^{ L,N} _{ \sigma_\ell ^-},$ while the total number of particles jumps from $N _{ \ell -1}$ at time $ \sigma_\ell ^-$ to $N_\ell $ at  time $ \sigma_\ell $ and the total (unit) mass of the empirical process remains constant.}
\\
During the time interval $[ \sigma_\ell , \sigma _{ \ell +1})$, we set
\begin{align}\label{eq-80}
\Zt ^{ L,N}_s=N_\ell  ^{ -1}\sum _{ i=1} ^{ N_\ell } \delta _{ Z_i(s)},\quad \sigma_\ell \le s< \sigma _{ \ell +1}
\end{align}
where  the random paths $(Z_i(s); \sigma_\ell \le s< \sigma _{ \ell +1}) _{ 1\le i\le N_\ell }$ are independent from each other and   each particle $1\le i\le N_\ell $  starts from $Z_i( \sigma_\ell )$ at time $s= \sigma_\ell ,$ and follows the Markovian evolution of the reference process $Z.$ Since this evolution specifies the Lagrangian $\mathcal{L}(s, \zeta_s,\dot\zeta_s)$, with \eqref{eq-67} we obtain  that the empirical process $\Zt ^{ L,N}$ obeys the large deviation principle with rate function
\begin{align*} 
I ^{ \kappa,L}(\zeta) =
	I_{\kappa _{ s_0}}(\zeta _{ \kappa _{ s_0}})
	+\sum _{ \ell =0} ^{ L-1} \int _{ \sigma_\ell } ^{ \sigma _{ \ell +1}} \kappa _{ \sigma_\ell } \mathcal{L}(s,\zeta_s,\dot\zeta_s)\,ds,\quad \zeta=(\zeta_s) _{ s_0\le s\le s_1}.
\end{align*}
The sum $\sum_\ell $ comes from the Markov property of the sample trajectories. It remains to let $L$ tend to infinity and remark that the following $\Gamma$-limit
\begin{align*}
\Glim L I ^{ \kappa, L}=I^\kappa,
\end{align*}
holds  as soon as $s\mapsto\kappa_s$ is continuous, with $I^\kappa$  defined at \eqref{eq-66}.  This implies that any limit point as $L$ tends to infinity  of the sequence of laws of $(\Zt ^{ L,N}) _{ L\ge 1}$ obeys the large deviation principle as $N$ tends to infinity with rate function $I^\kappa,$ see \cite{Ma12} for this argument. If the law of the  random path $Z$ is the unique solution to its martingale problem, for instance if $Z$ is a diffusion process with Lipschitz drift and diffusion fields as in the present setting, then for each $N\ge 1,$ there exists a unique limit point: $ \mathrm{Law}(\Zt ^{ \infty,N}),$ as $L$ tends to infinity. See \cite{Leo11a} for this argument.

\subsection*{Choice of the newcomers}
During this sketch of proof, we admitted that one can sample the newcomers from the support of $ \Zt ^{ L,N} _{ \sigma_\ell ^-}$ in such a way that their empirical measure
   is arbitrarily close, as $N$ tends to infinity, to $ \Zt ^{ L,N} _{ \sigma_\ell ^-}.$ Let us show this.\\
We want to pick, in an almost deterministic way, $N_\ell -N _{ \ell -1}$ distinct points $y'_j,$ $1\le j\le N_\ell -N _{ \ell -1},$ from the the set $\{ Z_i( \sigma_\ell ^-); 1\le i\le N _{ \ell -1}\}$.
\\
 Let us simplify notation.   We wish to choose $1\le n'\le n$   distinct points $y'_1,\dots,y' _{ n'}$ among a subset  $\{y_1,\dots,y_n\}$ of $\Rd$ in such a manner that $\tilde y ^{ n'}:=(n') ^{ -1}\sum _{ j=1}^{n'} \delta _{ y'_j}$ is arbitrarily close, in the sense of narrow convergence, to $\hat y^n:=n^{ -1}\sum _{ i=1}^{n} \delta _{ y_i}$  as  $n$ tends to infinity: 
 \begin{align}\label{eq-74}
 \lim _{ n, n'\to \infty} W_p(\tilde y ^{ n'},\hat y^n)=0
 \end{align}
 where $W_p$ stands for the Wasserstein distance of order $p\ge 1$. Suppose for the moment that  the supports of the empirical measures $\hat y^n$ are included in some fixed bounded box: 
\begin{align}\label{eq-75}
\cup_n \supp(\hat y^n)\subset \Lambda_R:=[-R,R] ^d,
\end{align}
for some large enough $R>0.$ 
 Take an integer $m\ge1$ (which is intended to tend to infinity as $n$ increases) and cover $ \Lambda_R$ by the $(2m)^d$ cubic boxes  $B_{ k _1,\dots,k _d}:=a _{ k _1,\dots,k _d}+ [0,R/m]^d$ where $a _{ k _1,\dots,k _d}:= m ^{ -1}R\ (k _1,\dots,k _d)$ with $k_1,\dots,k_d\in \{-m, -m+1,\dots, m-2,m-1\}.$ In each box $B,$ pick arbitrarily $\lfloor \hat y^n(B)\,n'\rfloor$ distinct points from $B\cap \supp(\hat y^n).$ It is because of this controlled arbitrariness that we describe these trials as \emph{almost deterministic}.  The collection of these  picked points is the main part of the set of newcomers. Their number $n''=\sum _{ k_1,\dots,k_d} \lfloor \hat y^n(B _{ k_1,\dots,k_d})\,n'\rfloor $ statisfies $n'-(2m)^d\le n''\le n'$ because taking the integer part may induce a lack of at most one point in each box. To complete the set of newcomers, simply add $n'-n''$ not already picked points from $\supp(\hat y^n)\cap \Lambda_R.$ It is immediate to see that for any $p\ge 1,$
\begin{align}\label{eq-71}
W_p(\tilde y ^{ n'},\hat y^n)\le D_R\,( 1/m + (2m)^d / n')
\end{align}
where $D_R=2 \sqrt{d} R$ is the diameter of the large box $ \Lambda_R$ and $D_R/m$ is the common diameter of the small boxes $B _{ k_1,\dots,k_d}$. At \eqref{eq-76} below, we shall tune    the coefficients $m$ and $n'$ as functions of $n$ in order to obtain \eqref{eq-74}.

\subsection*{An exponentially good approximation}

The assumption \eqref{eq-75} that the supports of all the measures $\hat y^n$ are included in some bounded set is crucial for this construction to be valid. But in general the supports of the random measures $\Zt ^{ L,N}_s$ are not uniformly bounded. This is why it is needed to rely on a sequence $\big((\Zt ^{ R,L,N} )_{ N\ge 1}\big) _{ R\ge 1}$ of approximations of  $(\Zt ^{ L,N} )_{ N\ge 1}$ living on larger and larger bounded sets $ \Lambda_R$ which is an   exponentially good approximation of $(\Zt ^{ L,N}) _{ N\ge 1}$ in the sense of the definition 
\cite[Def.\,4.2.14]{DZ} in order to apply the theorem \cite[Thm.\,4.2.16]{DZ}. The proxy $\Zt ^{ R,L,N}$ is  defined as $\Zt ^{ L,N}$, but  the reference Markov random path $(Z_s) _{ s_0\le s\le s_1}$ is replaced by its stopped version: 
\begin{align}\label{eq-73}
Z^R_s:= Z _{ s\wedge \tau_R},\quad s_0\le s\le s_1,
\end{align} 
where $\tau_R:=\inf\{s; s_0\le s\le s_1: Z_s\not\in \Lambda_ R\}\in [s_0,s_1]\cup \{\infty\}$ is the first exit time from the box $ \Lambda_R.$ Assuming that the sample paths of $Z$ are continuous, we have: $\sup_s|Z^R_s|\le R.$ 
The main argument for proving that this is indeed an exponentially good approximation is the following estimate
\begin{align}\label{eq-72}
\lim _{ R\to \infty}\sup_L  \limsup _{ N\to \infty} N ^{ -1} \log \PP\Big(\widehat Z^{L,N}\big( \{ \omega=( \omega_s) _{ s_0\le s\le s_1}; \sup _{s_0\le s\le s_1} | \omega_s|\ge R\}\big)> \delta\Big)= - \infty
\end{align}
which holds for any $ \delta>0.$ This result simply follows from Cramér's theorem applied to an iid sequence of Bernoulli random variables with parameter 
\[
\epsilon(R):=\PP( Z^*\ge R)
\le \EE Z^*/R.
\]
where we set $Z^*:= \sup _{s_0\le s\le s_1} | Z_s|.$
Indeed, Cramér's theorem states that, for any $L\ge 1,$
\begin{align*}
\limsup _{ N\to \infty} N ^{ -1} \log \PP\Big(\widehat Z^{L,N}\big( \{ \omega=( \omega_s) _{ s_0\le s\le s_1}; \sup _{s_0\le s\le s_1} | \omega_s|\ge R\}\big)> \delta\Big)
\le -h _{ \epsilon(R)}( \delta)
\end{align*}
where 
\begin{align*}
h _{ \epsilon(R)}(\delta)&= \delta \log(\delta/\epsilon(R))+(1-\delta)\log\big((1-\delta)/(1-\epsilon(R))\big)\\
&\ge \delta \log(\delta/\epsilon(R))+(1-\delta)\log(1-\delta)
\ge \delta\log( \delta R/\EE Z^*)+(1-\delta)\log(1-\delta),
\end{align*}
which implies \eqref{eq-72} as soon as
\begin{align*}
\EE Z^*:=\EE  \sup _{s_0\le s\le s_1} | Z_s|< \infty.
\end{align*}

\subsection*{Second model}
Previous construction is a little bit frustrating, since it does not provide an explicit representation of the limiting process $\Zt ^{ \infty,N},$ as $L$ and $R$ tend to infinity. Based on these previous considerations, we propose \emph{another} sequence of empirical processes $(\Zb^N) _{ N\ge 1}$ and \emph{some heuristics} for proving that its  large deviation rate function as $N$ tends to infinity is also $I^\kappa.$

The evolution of the empirical process $(\Zb^N_s; s_0\le s\le s_1)$ is a concatenation of very short  periods (as $N$ tends to infinity) separated by the increasing sequence of times
\begin{align}\label{eq-69}
\sigma_0:=s_0,\qquad  \sigma _{ \ell +1}:= \sigma_\ell + (\kappa' _{ \sigma_\ell }  N ^{1- a}) ^{ -1},\quad \ell \ge 0,
\end{align}
where  $0 <a<1$ and $\kappa'_s>0$ is the derivative of  $s\mapsto \kappa_s$ which is assumed to be differentiable.  Note that, although $ \sigma_\ell $ depends on $N,$ for a better readability we do not write explicitly this dependence.
At time $ \sigma_\ell $, the empirical process is 
 \begin{align*}
 \Zb ^{ N} _{ \sigma_\ell }=N _{ \ell } ^{ -1} \sum _{ i=1} ^{N _{ \ell } } \delta _{ Z_i( \sigma_\ell )},\quad \textrm{where}\quad N_\ell :=\lfloor \kappa _{ \sigma _{ \ell }}N\rfloor,
 \end{align*}
 and during the time interval $[ \sigma_\ell , \sigma _{ \ell +1})$ all the random paths $(Z_i; 1\le i\le N_\ell )$ evolve independently and follow the Markovian evolution of $Z^R$ defined at \eqref{eq-73}, where the large box parameter is chosen such that
 \begin{align*}
 R=O _{ N\to \infty}(N^b)
 \end{align*}
 for some constant $b>0.$  
   At time $ \sigma_\ell ^-$, the empirical process was 
 \begin{align*}
 \Zb ^{ N} _{ \sigma_\ell ^-}=(N _{ \ell -1}) ^{ -1} \sum _{ i=1} ^{N _{ \ell -1} } \delta _{ Z_i( \sigma_\ell ^-)}.
 \end{align*}
 \emph{Beware}: although we keep the notation $Z_i$, these random paths are not the same as in the first model.

  Let us describe how to choose the  newcomers $(Z'_j( \sigma_\ell ); 1\le j\le N_\ell -N _{ \ell -1})$ at time $ \sigma_\ell $. 
Remark that the increment ratio $(N _{ \ell +1}-N_\ell )/( \sigma _{ \ell +1}- \sigma_\ell )$ must be close to $N \kappa' _{ \sigma_\ell }.$ This implies that $N _{ \ell +1}-N_\ell \approx N\kappa' _{ \sigma_\ell } (\kappa' _{ \sigma_\ell } N ^{ 1- a}) ^{ -1}= N ^{ a}.$ Therefore, 
\begin{align*}
N _{ \ell +1}-N_\ell =  O _{ N\to \infty}(N^a).
\end{align*}
 We  choose the newcomers $Z'_j( \sigma_\ell )$  as in previous subsection, with $R=O _{ N\to \infty}(N^b)$ as above and the small box parameter
 \begin{align*}
 m=O _{ N\to \infty}(N^c),
 \end{align*}
 for some $c>0.$
As we did at \eqref{eq-68}, we set 
\begin{align}\label{eq-68c}
\Zb^N _{ \sigma_\ell }:=N_\ell  ^{ -1} \Big(N _{ \ell -1} \Zb^N _{ \sigma_\ell ^-}+ \sum _{ j=1} ^{N_\ell -N _{ \ell -1}} \delta _{ Z'_j( \sigma_\ell )}\Big).
\end{align}
We want that the cumulated error vanishes as $N$ tends to infinity:
\begin{align*}
\sum_\ell  W_p(\Zb^N _{ \sigma_\ell ^-},\Zb^N _{ \sigma_\ell }) \underset{N\to \infty}\longrightarrow 0.
\end{align*}
Considering  \eqref{eq-71} with $n=O _{ N\to \infty}(N)$, $n'=O _{ N\to \infty}(N^a)$, $R=O _{ N\to \infty} (N^ b)$  and $m=O _{ N\to \infty} (N^ c)$,  one sees that 
$	
W_\ell :=W_p(\Zb^N _{ \sigma_\ell ^-},\Zb^N _{ \sigma_\ell })
	=O _{ N\to \infty}(N ^{ b-c}+N^{-a+ b +cd}).
$	
On the other hand, \eqref{eq-69} implies that $ \Delta \sigma_\ell := \sigma _{ \ell +1}- \sigma_\ell =O _{ N\to \infty}(N ^{a-1}).$ Therefore
\begin{align*}
\sum_\ell  W_p(\Zb^N _{ \sigma_\ell ^-},\Zb^N _{ \sigma_\ell })\le |s_1-s_0| \sup_\ell  (W_\ell / \Delta \sigma_\ell )
	=O _{ N\to \infty}(N  ^{ \max(1- a+b-c, 1-2a+ b+cd)})
\end{align*}
 vanishes as $N$ tends to infinity if  one chooses $0<a<1,$ $b>0$ and $c>0$ such that 
$1- a+b-c<0$ and  $1-2a+ b+cd<0.$ One can take for instance:
\begin{align}\label{eq-76}
a= \frac{d+2}{d+3},\quad 
b= \frac{1}{3d(d+3)}
\quad \textrm{and}\quad
c= \frac{2d+1}{2d(d+3)}.
\end{align}
We have shown that
\begin{align}\label{eq-79}
\sum_\ell  W_p(\Zb^N _{ \sigma_\ell ^-},\Zb^N _{ \sigma_\ell })	=O _{ N\to \infty}(N  ^{- \gamma}),\quad \textrm{for some}\quad \gamma>0,\quad \textrm{almost surely}.
\end{align}

\subsection*{General case for $\kappa$} \label{subsec-kappa}
Up to now we decided, to make things easier, to assume that $s\mapsto \kappa_s$ is increasing. The general case where $\kappa$ is any positive continuously differentiable function is simply obtained by replacing 
 \eqref{eq-68c}  by
\begin{align}\label{eq-68b}
 \Zb ^{ N} _{ \sigma_\ell }= N_\ell  ^{ -1}\Big( N _{ \ell -1} \Zb ^{ N} _{ \sigma_\ell ^-}
 	+ \mathrm{sign}(\kappa _{ \sigma_\ell }-\kappa _{ \sigma _{ \ell -1}}) \sum _{ j=1} ^{|N _{ \ell }-N _{ \ell -1}|} \delta _{ Z'_j( \sigma_\ell )}\Big)
	=N _{ \ell }^{ -1} \sum _{ i=1} ^{N _{ \ell } } \delta _{ Z_i( \sigma_\ell )}
\end{align}
  and replacing \eqref{eq-69} by
\begin{align}\label{eq-69b}
\sigma_0:=s_0,\qquad  \sigma _{ \ell +1}:= \sigma_\ell + (|\kappa' _{ \sigma_\ell }|  N ^{1- a}) ^{ -1}+N ^{ -1},\quad \ell \ge 0.
\end{align}
Last equality in \eqref{eq-68b} results from some relabeling, and  last term $N ^{ -1}$ in \eqref{eq-69b} is added to make sure that $ \sigma _{ \ell +1}> \sigma_\ell $ even when $\kappa' _{ \sigma_\ell }=0.$

Looking at \eqref{eq-68b}, we see that when $\kappa _{ \sigma_\ell }< \kappa _{ \sigma _{ \ell -1}},$ the newcomers are removed from the support of $ \Zb ^{ N} _{ \sigma_\ell ^-}$ and  when $\kappa _{ \sigma_\ell }>\kappa _{ \sigma _{ \ell -1}},$ the newcomers are added to  $ \Zb ^{ N} _{ \sigma_\ell ^-}$, while the total (unit) mass of  $ \Zb ^{ N} _s$  is conserved.

\subsection*{Non-locality} 

A collection of a number $N$ of $k$-mappings is a cloud of $Nk$ particles. As soon as this cloud is subject to some physical mechanism, its very structure implies non-locality. For instance, suppose that some $k$-mapping gives birth to a new $k$-mapping, as was suggested before,  or that it dies or is pushed by some force field, then  the $k$ particles corresponding to the $k$ coordinates  of the $k$-mapping, branch or die or respond to the force field \emph{simultaneously}, although they might be very far apart. There is an action at a distance. This is no surprise since we already noticed at the very end of Section \ref{sec-MAG-mapping} (page  \pageref{sec-nonlocal})  that, similarly to Newton's theory of gravity, \MAG\ is a non-local theory. 

However, the splitting mechanism reveals this non-locality in a dramatic way since, unlike the effect of the displacement of remote particles which decreases when their distance  increases (as is the case for \MAG, see page \pageref{sec-non-local-2} for an illustration), there is no notion of tiny splitting: it occurs or not.  We do not know if there is some way to circumvent this difficulty by replacing this splitting mechanism by something more flexible.

\subsection*{Is this splitting hypothesis pertinent?}

Besides this non-locality issue,
 this splitting hypothesis\footnote{In physics, a hypothesis is not an assumption; it is waiting to be confirmed or dispelled. Here and below we use this term in a colloquial sense, not in the mathematical sense.}  is strange: how to incorporate it in a gravitation model?  Another physical hypothesis would be that no such splitting occurs but the temperature decreases with time in such a way that  $ \kappa_s$ (which is an inverse temperature) satisfies the Parameter setting \ref{hyp-02}.  This would work fine.

 Nevertheless, as soon as we start with a cloud of Brownian particles and once we do not assume such a decreasing temperature, the splitting  appears when doing the maths. In any case, we  decided to propose this splitting picture because it might be of some use in other types of model.

\section{$\epsilon$-\MAG\ for a $k$-fluid. Quantum force?}
\label{sec-B-k-fluid}

It remains to remove the rightmost term from the action \eqref{eq-83} to arrive at $ \epsilon$-\MAG's action \eqref{eq-37}. This amounts to subtract the potential energy 
\begin{align*}
 \mathcal{I}^\epsilon _t(q):= 4 e ^{ 4t} \epsilon^2  I(q| m^ \epsilon_t)
\end{align*}
 from the Lagrangian $(t,q,\dot\qq)\mapsto \uud    \| \dot\qq-\dot\mm^\epsilon _t\|^2_{q}
	+   \mathcal{I}^\epsilon (t,q)$ to arrive at the desired Lagrangian
\begin{align*}
(t,q,\dot\qq) \mapsto \uud    \| \dot\qq-\dot\mm^\epsilon _t\|^2_{q}.
\end{align*}

\begin{statement}
 In terms of a Newton equation,  this means that the additional force field
 \begin{align}\label{eq-qff}
 - \gOW_q  \mathcal{I}^\epsilon _t
 \end{align}
  is applied to the $k$-fluid that minimizes the action \eqref{eq-83}. 
 \end{statement}

\begin{phyp}\label{physhyp}
The force field \eqref{eq-qff} is a \emph{quantum} force field.
\end{phyp}
 
In the remainder of this section  we present two known results (Theorems \ref{res-21} and \ref{res-21b}) that support this  hypothesis by analogy. However, the hypothesis is rather vague and the way this  force emanates from the fluid itself or its surrounding  remains unknown to us.

\subsection*{Reference path measure}

For the remainder of this section, we drop \MAG's specific setting and we  consider an abstract setting where the  state  space is $\Rn$, the measure
 \[
 m(dx):=m(x)\, dx
 \]
on $\Rn$  is the equilibrium measure of the Markov process with generator
\begin{align*}
A^m:= \nabla \log \sqrt{m} \scal \nabla+   \Delta/2,
\end{align*}
where the  function  $m:\Rn\to(0, \infty)$ is positive everywhere and differentiable. 
 It is assumed that there is a unique path measure $R^m$ that solves the martingale problem with generator $A^m$ and initial measure $m.$  Denoting $(X_t) _{ t_0\le t\le t_1}$ the canonical process, this is equivalent to
\begin{align}\label{eq-89}
\begin{split}
&dX_t=  \nabla \log \sqrt{m} (X_t)\, dt +  dB ^{ R^m}_t,\quad t_0\le t\le t_1,\quad R^m\ae,\\
&X _{ t_0}\sim m,
\end{split}
\end{align}
where $B ^{ R^m}$ is an $R^m$-Brownian motion. \emph{Not only $R^m$ is $m$-stationary, but also  it is reversible.} 

\begin{definition}[Reference path measures]
Let 
$\Rem$ be the $m$-reversible Markov measure with generator $ \epsilon A^m$ and initial distribution $m$, that is 
\begin{align}\label{eq-y1}
\Rem(d \omega)= R^m( X^ {(\epsilon)}( d\omega))
\end{align}
 with $X^{ (\epsilon)}_t:= X _{t/ \epsilon },$ meaning that $\Rem$ is $R^m$ slowed down by a factor $ \epsilon$. We choose  as a reference path measure
\begin{align}\label{eq-91}
\Re:=&\exp \left( \epsilon ^{ -1} \int _{ t_0} ^{ t_1} V_t(X_t)\, dt\right) \,\Rem,
\end{align}
where  $V$ is some scalar potential.
\end{definition}

 \subsection*{Schrödinger bridges and entropic interpolations} 
 
We  recall some results from  \cite{Zam86} and \cite{Co18} about the Schrödinger problem that we already met at \eqref{eq-101}.  This will help us to understand the force field   $+ \gOW_q  \mathcal{I}^\epsilon _t$, with a plus  instead of a minus sign, which is easier to interpret than its opposite.

\subsubsection*{The Schrödinger problem and its solution}

It is known that, in a generic situation, the unique solution of the \emph{Schrödinger problem}
\begin{align*}
\inf _{ Q: Q _{ t_0}= q _{ t_0}, Q _{ t_1}= q _{ t_1}}  H(Q|R)
\end{align*}
where $Q$ is a path measure with prescribed initial and final marginals: $q _{ t_0}$ and $q _{ t_1}$, writes as
\begin{align} \label{eq-92}
Q= f _{ t_0}(X _{ t_0})g _{ t_1}(X _{ t_1})\, R,
\end{align}
for some nonnegative measurable functions $f _{ t_0}$ and $g _{ t_1}$ on $\Rn$, see \cite{Leo12e} and the references therein for an overview of the Schrödinger problem.  Let us recall the Definitions \ref{def-x2} in this broader context.

\begin{definitions}[Schrödinger bridge, entropic interpolation]
The entropic minimizer $Q$ is called the \emph{Schrödinger bridge} with respect to  $R$  between $q _{t_0}$ and $q _{ t_1}$. \\
The entropic interpolation with respect to  $R$ between $q _{t_0}$ and $q _{ t_1}$  is the flow of time-marginals $(q_t:=(X_t)\pf Q; t_0\le t\le t_1)$   of the Schrödinger bridge $Q.$ 
\end{definitions}
If $R$ is Markov, then the Schrödinger bridge $Q$ inherits the Markov property from it.

\subsection*{Entropy versus quantum}

This subsection is dedicated to the presentation of two known results that exhibit a striking analogy between the entropic and quantum evolutions. This is a modern way, in terms of \OW-geometry, of recasting Schrödinger's original insight \cite{Sch31,Sch32}.

 Let us denote the average potentials by
\begin{align*}
 \mathcal{V}_t(q):= \int _{ \Rn} V_t\, dq
 \qquad \textrm{and}\qquad
 \mathcal{V}_t ^{ m, \epsilon}(q)
	:= \int V_t ^{ m, \epsilon} \, dq
\end{align*}
with
\begin{align*}
 V_t ^{ m, \epsilon}:=V_t+ \epsilon^2 \mathcal{Q}(m|\Leb)
\end{align*}
where the quantum potential $ \displaystyle{ \mathcal{Q}(m|\Leb):= - \frac{A^m \sqrt{m}}{ \sqrt{m}}}$ already appeared at \eqref{eq-103}.

\begin{theorem}[Conforti \cite{Co18}] \label{res-21}
Any entropic interpolation $(q_t)$ with respect to $\Re$ solves the Newton equation
\begin{align}\label{eq-y3}
\ddot\qq_t= -\gOW _{ q_t}( \mathcal{V}_t -\epsilon^2   I(\sbt|m))
\end{align}
in the \OW-manifold. This is equivalent to
\begin{align}\label{eq-y5}
\ddot\qq_t= -\gOW _{ q_t}( \mathcal{V} ^{ m, \epsilon}_t -\epsilon^2   I(\sbt|\Leb)).
\end{align}
\end{theorem}

\begin{remark}
Lebesgue measure plays a privileged role in quantum mechanics. This is the reason why we state at \eqref{eq-y5} this Newton equation in terms of the standard Fisher information $I(\sbt|\Leb).$ This will be convenient for a comparison with its quantum analogue \eqref{eq-y4}.
\end{remark}

Let $\Psi$ be a complex function solving the  \emph{Schrödinger equation}
\begin{align}\label{eq-y2}
\big(- i\hbar \partial_t  -   \hbar^2   \Delta  /2 +V\big)\Psi =0.
\end{align}
 Born's formula:
\begin{align*}
q_t=\Psi_t\overline\Psi_t=|\Psi_t|^2
\end{align*}
states that the measure $q_t(dx):=q_t(x)\,dx$ can be interpreted as the density of presence of the state of the system at time $t$. Assuming that $ \Psi _{ t_*}$ is normalized at some instant $t_*$, that is $q _{ t_*}(\RR^n)=1,$  it is also normalized at any time, so that $(q_t)$ is a flow of probability measures.

\begin{theorem}[von Renesse \cite{vR11}] \label{res-21b}
Any normalized solution $\Psi$ of the  Schrödinger equation \eqref{eq-y2} is such that $q:=|\Psi|^2$ solves the Newton equation
\begin{align}\label{eq-y4}
\ddot\qq_t= -\gOW _{ q_t}( \mathcal{V}_t +\hbar^2   I(\sbt|\Leb))
\end{align}
in the \OW-manifold.
\end{theorem}

 Apart from the minor transformation: $ \mathcal{V}\leftrightarrow \mathcal{V} ^{ m, \epsilon}$ (both potentials are linear and therefore responsible for  uniform force fields in the \OW-manifold),  we see that we pass from the thermal evolution \eqref{eq-y5} to the quantum evolution \eqref{eq-y4} and vice-versa by  applying the rule
\begin{align*}
-\epsilon^2 I(q|\Leb)\leftrightarrow +\hbar^2 I(q|\Leb).
\end{align*}
This makes a huge difference between the dynamics. Thinking of $I(\sbt|\Leb)$ as a convex function (it is convex along the usual affine combinations of probability measures), the concave  potential $-\epsilon^2 I(q|\Leb)$ tends to generate relaxation trajectories (typically thermal relaxation to equilibrium)  or exponentially divergent ones, while the convex potential $\hbar^2 I(q|\Leb)$ tends to generate oscillating trajectories, a typical quantum behavior.

\subsection*{Main ideas of the "proofs" of Theorems \ref{res-21} and \ref{res-21b}}

For the comfort of the reader, we give some details of the ``proofs'' of these theorems. They rely on Propositions \ref{res-19}, \ref{res-22} and Lemma \ref{res-20} below.

To our knowledge, there are no proofs of these theorems in a wide enough generality. Theorem \ref{res-21}  in \cite{Co18} holds under strong regularity assumptions and  \cite{vR11}'s proof of Theorem \ref{res-21b} is only formal, very much in the spirit of what follows.

\begin{proposition}[Entropic version of Born's formula]\label{res-19}
Let $Q= f _{ t_0}(X _{ t_0})g _{ t_1}(X _{ t_1})\, \Re$ be a Schrödinger bridge with respect to $\Re,$ see  \eqref{eq-92}, and let  $q_t=(X_t)\pf Q$ denote the corresponding entropic interpolation.  
For any $t,$   $q_t$ is absolutely continuous with respect to $m$ and 
\begin{align}\label{eq-85}
d q_t/dm= f_tg_t, \quad t_0\le t\le t_1,
\end{align}
where for any $t_0\le t\le t_1$ and $x\in\Rn,$ 
\begin{align}\label{eq-84}
\begin{split}
f_t(x)&:=E _{ \Rem}\left(f _{ t_0}(X _{ t_0}) \exp \Big( \epsilon ^{ -1} \int _{ t_0} ^{ t} V_s(X_s)\, ds\Big) \mid X_t=x\right),\\
g_t(x)&:=E _{ \Rem}\left(\exp \Big( \epsilon ^{ -1} \int _{ t} ^{ t_1} V_s(X_s)\, ds\Big)g _{ t_1}(X _{ t_1})\mid X_t=x\right).
\end{split}
\end{align}
\end{proposition}

\begin{proof}
By a general result of integration theory 
\[
\frac{dq_t}{dm}(x)
	=\frac{d(X_t)\pf Q}{d(X_t)\pf \Rem }(x)
	=E _{ \Rem }\Big( \frac{dQ}{d\Rem }\mid X_t=x\Big).
\]
With \eqref{eq-91} and \eqref{eq-92}, we see that for any $t_0\le t\le t_1,$
\begin{align*}
\frac{dq_t}{dm}(X_t)
	&=E _{ \Rem } \left( \frac{dQ}{d\Re}\, \frac{d\Re}{d\Rem }\mid X_t \right) \\
	&=E _{ \Rem }\Big(f _{ t_0}(X _{ t_0}) g _{ t_1}(X _{ t_1})\exp \Big( \epsilon ^{ -1} \int _{ t_0} ^{ t_1} V_t(X_t)\, dt\Big)\mid X_t \Big)\\
	&=E _{ \Rem }\Big(f _{ t_0}(X _{ t_0})\exp \Big( \epsilon ^{ -1} \int _{ t_0} ^{ t} V_s(X_s)\, ds\Big)\\
	&\hskip 5cm \times \exp \Big( \epsilon ^{ -1} \int _{ t} ^{ t_1} V_s(X_s)\, ds\Big) g _{ t_1}(X _{ t_1})\mid X_t \Big)\\
	&=E _{ \Rem }\Big(f _{ t_0}(X _{ t_0})\exp \Big( \epsilon ^{ -1} \int _{ t_0} ^{ t} V_s(X_s)\, ds\Big)\mid X_t \Big)\\
	&\hskip 5cm \times E _{ \Rem }\Big(\exp \Big( \epsilon ^{ -1} \int _{ t} ^{ t_1} V_s(X_s)\, ds\Big) g _{ t_1}(X _{ t_1})\mid X_t \Big)\\
	&=:f_t g_t(X_t),
\end{align*}
where we used the Markov property of $\Rem $ at last but one equality.
\end{proof}

The functions $f(t,x)$ and $g(t,x)$ are solutions of the  parabolic equations
\begin{align}\label{eq-93}
\left\{
\begin{array}{ll}
&(- \epsilon \partial_t+ \epsilon^2 A^m+V)f=0,\\
&f(t_0)=f _{ t_0},
\end{array}
\right.
\qquad
\left\{
\begin{array}{ll}
&( \epsilon\partial_t+ \epsilon^2 A^m+V)g=0,\\
&g(t_1)=g_{ t_1},
\end{array}
\right.
\end{align}
and \eqref{eq-84} are their Feynman-Kac representations. The Born-like formula \eqref{eq-85} was first established by Schrödinger in \cite{Sch31,Sch32} in the Brownian case and extended later by Zambrini in \cite{Zam86}. This "Markovian" proof first  appeared in \cite{Leo12e}. Introducing 
\begin{align*}
\theta_t:= \epsilon\log \sqrt{g_t/ f_t},
\qquad \rho_t:=d q_t/dm,
\end{align*}
we see  with \eqref{eq-85} that
\begin{align}\label{eq-94}
f= \sqrt{ \rho} e ^{ - \theta/ \epsilon},\qquad g= \sqrt{ \rho} e ^{  \theta/ \epsilon}.
\end{align}
Let us introduce the relative quantum potential
\[\mathcal{Q}(q|m):= - \frac{A^m \sqrt{\rho}}{ \sqrt{\rho}}
\quad \textrm{where}\quad \rho= dq/dm. \]
Several basic properties of this  potential are gathered at Appendix \ref{app-e}.

\begin{proposition}\label{res-22}
Let $Q$ be any Schrödinger bridge with respect to $\Re$ and let  $q_t=(X_t)\pf Q$ denote the corresponding entropic interpolation.  We have
 \begin{align}\label{eq-86}
\begin{split}
&\partial_tq+\nabla\cdot(q\, \nabla \theta)=0,\\
&\partial_t\theta+  |\nabla \theta|^2/2 + V-  \epsilon^2    \mathcal{Q}(q|m)=0.
\end{split}
\end{align}
\end{proposition}

\begin{proof}["Proof"] This "proof" is incomplete: we pass over regularity issues. 
\\
Besides $ \theta,$ let us introduce the functions 
\[\varphi:= \epsilon\log f,\quad \psi:= \epsilon\log g \quad \textrm{and}\quad \eta:= ( \varphi+ \psi)/2.\] We also see that $ \theta= ( \psi- \varphi)/2.$ The  parabolic equations \eqref{eq-93} are equivalent to 
$  e ^{ -\varphi/ \epsilon}(- \epsilon\partial_t+ \epsilon^2 A^m+V) e ^{ \varphi/ \epsilon}=0$ and 
$  e ^{ - \psi/ \epsilon}(\epsilon\partial_t+\epsilon^2 A^m+V) e ^{ \psi/ \epsilon}=0,$ that is
\begin{align*}
(a)\quad&- \partial_t \varphi+|\nabla \varphi|^2/2+V+ \epsilon A^m \varphi=0,\\
(b)\quad&\hskip 0.5cm \partial_t \psi+|\nabla \psi|^2/2+V+ \epsilon A^m \psi=0.
\end{align*}
Taking the half difference $[(b)-(a)]/2$ and  half sum $[(a)+(b)]/2$ of these equations, we obtain
\begin{align*}
(i)\quad&\partial_t \eta +\nabla \theta\cdot\nabla \eta+ \epsilon A^m \theta =0,\\
(ii)\quad&\partial_t \theta+|\nabla \theta|^2/2+V+ \epsilon A^m\eta+|\nabla \eta|^2/2 =0.
\end{align*}
By \eqref{eq-85}:  $ \eta= \epsilon\log \sqrt{q}- \epsilon\log \sqrt{m},$ and because $ \partial_t m=0,$ equation $(i)$ writes as : 
\begin{align*}
  \partial_t \log \sqrt{q}
 	+ \nabla\log \sqrt{q} \cdot \nabla \theta
	+  \Delta \theta/2
	=0.
\end{align*}
Multiplying by $2q,$ this amounts to
\begin{align*}
0= \partial_t q +  \nabla q\cdot \nabla \theta  + q \Delta \theta
	= \partial_t q+\nabla \cdot(q\nabla \theta),
\end{align*}
which is the continuity equation in \eqref{eq-86}.\\ 
The Hamilton-Jacobi equation in \eqref{eq-86} is simply equation $(ii)$, once one notices that  
\begin{align*}
\epsilon A^m\eta+|\nabla \eta|^2/2
	&=  \epsilon\nabla\log \sqrt{m}\cdot \nabla\eta + \epsilon \Delta\eta/2+|\nabla \eta|^2/2\\
	&= \epsilon^2 \Big[\nabla\log \sqrt{m}\cdot \nabla\log \sqrt{ \rho} 
		+ \uud \Delta\log \sqrt{ \rho}+\uud |\nabla \log \sqrt{ \rho}|^2\Big]\\
	&\overset{\eqref{eq-99}}=\epsilon^2\Big[ \nabla\log \sqrt{m}\cdot  \frac{\nabla \sqrt{ \rho}}{ \sqrt{ \rho}} 
		+  \frac{ \Delta \sqrt{ \rho}}{2 \sqrt{ \rho}}\Big]
	= \epsilon^2\frac{ A^m \sqrt{ \rho}}{ \sqrt{ \rho}} 
	=- \epsilon^2 \mathcal{Q}(q|m),
\end{align*}
where we used $ \eta= \epsilon\log \sqrt{ \rho}  $ at second equality.
\end{proof}

\begin{remark}
The system of equations \eqref{eq-86} is the entropic analogue of Madelung's equations in quantum mechanics \cite{Ma26}, see \eqref{eq-86c} below. One passes from the entropic to the quantum versions of these equations  applying the quantization rule: $\epsilon \leftrightarrow i\hbar$, that is simply replacing $- \epsilon^2  \mathcal{Q}(q|m)$ by $+ \hbar^2  \mathcal{Q}(q|m).$
\end{remark}

\begin{lemma}[\cite{vR11}] \label{res-20}
The Otto-Wasserstein gradients of $ \mathcal{V}_t: p\mapsto \int V_t\, dp$ and  $ I(\sbt|m_t)$ are
\begin{align*}
(\textrm{a})\quad
 \gOW_p \mathcal{V}_t=\nabla V_t\qquad \textrm{and} \qquad
(\textrm{b})\quad \gOW_p I(\sbt|m_t)=   \nabla \mathcal{Q}(p|m_t).
\end{align*}
\end{lemma}

\begin{proof}["Proof"] This "proof" relies on Otto's heuristics, see Appendix \ref{app-OW}.  We skip the regularity problems.
\\
For any regular path $(p_\tau ),$ 
\begin{align} \label{eq-88}
\frac{d}{d\tau} \mathcal{F}(p_\tau )= ( \gOW _{ p_\tau } \mathcal{F},\dot\pp_\tau )^{ \mathrm{OW}} _{ p_\tau }
=\int _{ \Rn}  \gOW _{ p_\tau }  \mathcal{F}\cdot\dot\pp_\tau \, dp_\tau 
\end{align}
where $(\nabla u,\nabla v)^{ \mathrm{OW}}_p=\int _{ \Rn} \nabla u\cdot\nabla v\, dp$ is the  inner product of the tangent space at $p$ of the \OW-manifold.

The dependence in $t$ of $ \mathcal{V}_t$ or $I(\sbt|m_t)$ does not play any role in the computation of the "spatial" $p$-gradients which must be considered as partial derivatives in $p$ with a fixed time parameter $t$. Therefore, we shall write $V, \mathcal{V},m$  instead of $V_t, \mathcal{V}_t,m_t$ without any loss of generality.
\begin{enumerate}[(a)]
\item
Taking $ \mathcal{F}= \mathcal{V},$ we see that
\begin{align*}
\frac{d}{d\tau} \mathcal{V}(p_\tau )
	=\int _{ \Rn} V \partial_\tau  p_\tau \,d\Leb
	=-\int _{ \Rn} V\,\nabla\scal(p_\tau \dot\pp_\tau )\,d\Leb
	=\int _{ \Rn}\nabla V\scal \dot\pp_\tau \, dp_\tau .
\end{align*}
We used the continuity equation $ \partial_\tau  p+\nabla \scal(p\dot\pp)=0$ at second identity and a standard integration by parts at last identity. Comparing with \eqref{eq-88} leads to the announced result.
\item
For any small perturbation $x\mapsto h(x)$ with a small gradient, we have
\begin{align*}
\begin{split}
|\nabla \sqrt{ \rho+ h}|^2
= |\nabla (\sqrt \rho&+  h/(2 \sqrt{ \rho})+o(h))|^2
=|\nabla \sqrt{ \rho}|^2+\nabla \sqrt{ \rho}\cdot\nabla( h/ \sqrt{ \rho})
	+ o( h, \nabla h)\\
&=|\nabla \sqrt{ \rho}|^2
	+\nabla\log \sqrt{ \rho}\cdot\nabla h 
	-|\nabla\log \sqrt{ \rho}|^2 h
	+ o( h, \nabla h).
\end{split}
\end{align*}
Hence,  denoting $\rho_\tau = dp_\tau /dm$ and taking $ \mathcal{F}=I(\sbt|m)$ in \eqref{eq-88}, we see that
\begin{align*}
\frac{d}{d\tau}I( p_\tau |m)
	&
=\ud \int _{ \Rn} \partial_\tau |\nabla \sqrt{ \rho_\tau }|^2\, d m
	=\ud \int _{ \Rn} ( \nabla\log \sqrt{ \rho_\tau }\cdot\nabla \partial_\tau  \rho_\tau 
	-|\nabla\log \sqrt{ \rho_\tau }|^2 \partial_\tau  \rho_\tau )\, dm\\
	&\overset{(i)}=-\int _{ \Rn} ( A^m\log \sqrt{ \rho_\tau } 
	+ |\nabla\log \sqrt{ \rho_\tau }|^2/2 )\, \partial_\tau  p_\tau \, d \Leb \\
	&\overset{(ii)}= \int _{ \Rn} \mathcal{Q}(p_\tau |m)\, \partial_\tau  p_\tau \, d \Leb
	\overset{(iii)}=  \int _{ \Rn} \nabla  \mathcal{Q}(p_\tau |m)\cdot \dot \pp_\tau  \, dp_\tau .
\end{align*}
At $\overset{(i)}=$, we used $ \partial_\tau  \rho_\tau \,dm= \partial_\tau  p_\tau \, d \Leb$  and  the integration by parts formula \eqref{eq-87}.
 Equality  $\overset{(ii)}=$ follows from $ \Delta u/u= \Delta\log u+|\nabla\log u|^2.$
At $\overset{(iii)}=$, we used the continuity equation $ \partial_\tau  p+\nabla \scal(p\dot\pp)=0$ and a standard integration by parts. We conclude comparing with \eqref{eq-88}.
\end{enumerate}
This completes the "proof" of the lemma.
\end{proof}

\begin{proof}[``Proof'' of Theorem \ref{res-21}] \label{pf-res-21} This proof is incomplete because, although 
the proof of Proposition \ref{res-19} was complete,  the proofs of Proposition \ref{res-22} and Lemma \ref{res-20} were only sketched: we passed over  regularity issues and relied upon  Otto's heuristics.

The first identity in  \eqref{eq-86} is a continuity equation: $\nabla \theta$ is the current velocity of $(q_t)$. Since in addition it is a gradient, we have $\dot\qq_t= \nabla \theta_t.$  The second identity   in \eqref{eq-86} is a Hamilton-Jacobi equation. Taking its gradient (again we skip a  regularity problem)  leads us to Newton's equation: $ (\partial_t+\nabla \theta_t\scal\nabla) \nabla \theta_t+\nabla V_t- \epsilon^2 \nabla  \mathcal{Q}(q_t|m)=0.$ But 
$(\partial_t+\nabla \theta_t\scal\nabla) \nabla \theta_t
=\nabla ^{ \mathrm{OW}} _{ \nabla \theta_t }\nabla \theta_t 
=\nabla ^{ \mathrm{OW}} _{ \dot\qq_t}\dot\qq_t=:\ddot\qq_t.$ Hence,   $\ddot\qq_t=-\nabla V_t + \epsilon^2 \nabla\mathcal{Q}(q_t|m).$ We obtain \eqref{eq-y3} with  Lemma \ref{res-20}.

With  \eqref{eq-95b}, one transforms \eqref{eq-y3} into \eqref{eq-y5}.
\end{proof}

\begin{proof}[``Proof'' of Theorem \ref{res-21b}] 
In \cite{Sch31,Sch32,Zam86}, the standard case $m=\Leb$ is considered. In this case  \eqref{eq-93} writes as
\begin{align*}
\left\{
\begin{array}{rcl}
- \epsilon \partial_t  f+   \epsilon^2   \Delta  f/2 +V  f&=&0\\
+ \epsilon \partial_t  g+   \epsilon^2   \Delta  g/2+V g&=&0.
\end{array}
\right.
\end{align*}
Applying the correspondences
\begin{align}\label{eq-96}
\epsilon \leftrightarrow i\hbar,
\qquad
( f, g)\leftrightarrow (\Psi,\overline\Psi)
\end{align}
where $\Psi$ is a complex function and $\overline\Psi$ is its complex conjugate, we obtain
\begin{align*}
\big(- i\hbar \partial_t  -   \hbar^2   \Delta  /2 +V\big)\Psi =0,
\end{align*}
the second equation:
$+ i\hbar \partial_t \overline\Psi-  \hbar^2   \Delta\overline\Psi/2+V \overline\Psi=0,$ being its complex conjugate. This is Schrödinger equation for the wave function $\Psi$, while \eqref{eq-85} becomes Born's formula:
\begin{align*}
q=\Psi\overline\Psi=|\Psi|^2.
\end{align*}
Formula \eqref{eq-94} becomes the polar decomposition $\Psi= \sqrt{q} e ^{ i \theta/\hbar}$  and \eqref{eq-86} with $m=\Leb$ becomes the couple of Madelung equations \cite{Ma26}
\begin{align}\label{eq-86c}
\begin{split}
&\partial_tq+\nabla\scal(q\, \nabla \theta)=0,\\
&\partial_t\theta+  |\nabla \theta|^2/2 + V+\hbar^2 \mathcal{Q}(q|\Leb) =0,
\end{split}
\end{align}
leading, as in the "proof" of Theorem \ref{res-21}, to  Newton's equation\begin{align*}
\ddot\qq_t=-\gOW _{ q_t}  \Big( \mathcal{V}+ \hbar^2 I(\sbt|\Leb)\Big).
\end{align*}
This result is Theorem \ref{res-21b}. It was  discovered by von Renesse   \cite{vR11}.
\end{proof}

\begin{remark}
Interpreting the equivalent correspondences \eqref{eq-96b} or \eqref{eq-96} is far from common sense. They are a \emph{quantization rule}. Although Schrödinger did not state them explicitly, it is clearly \eqref{eq-96} that he had in mind when writing  \cite{Sch31,Sch32}. 
\end{remark}

\subsection*{Which quantum force?} 
Let $\Psi$ be a solution of  
the \emph{nonlinear} Schrödinger equation 
\begin{align}\label{eq-nls}
- i\hbar \partial_t  \Psi-   \hbar^2   \Delta \Psi  /2 +V(|\Psi|^2)\, \Psi =0,
\end{align}
where $V(x,q)$ is a scalar potential depending on the  measure  $q$. 

\begin{enumerate}
\item
Is there a corresponding 
 Newton equation \[\ddot \qq=-\gOW_q \mathcal{U}\] in the \OW-manifold  satisfied by $q_t=|\Psi_t|^2\,\Leb$?
 \item
  In case this holds, what is the relation between $ \mathcal{U}(q)$ and $V(q,\sbt)?$
  \item
   What in particular if one takes $ \mathcal{U}(q)=I(q|\Leb)$ (which is a rather wild function)?
\end{enumerate}
 
 Answering these questions is necessary to understand the  physical mechanism that is responsible     for  the quantum force    $-\gOW I(\sbt|\Leb)$ because the relevant information is carried by the potential $V(q,x).$ 

Since we gave some evidence that  quantum forces decrease temperature, the variation of the coefficient in front of the Laplacian must play some role. Here we simply wrote $\hbar^2$ which, unlike our $ \epsilon^2,$ is a universal constant that cannot vary. But a complete Schrödinger equation involves $\hbar^2/m$ instead of $\hbar^2.$ Some variation of mass  might play some role. 
 
 It is also important to have in mind that for a usual fluid, $q$ is a measure on the state space (for instance $\RR^3$), while a wave function $\Psi$,  solution of a Schrödinger equation,   is a function on the configuration space (for instance $\RR ^{ 3N})$. One must project the abstract measure  $|\Psi|^2\, \Leb$ from the configuration  space onto the state space to obtain the evolution of the fluid. One can proceed as at  page \pageref{sec-fff}, introducing some projection similar to $\proj _{ dk\to d}$,  in the spirit of the Density Functional Theory \cite{CF23,Fri24}.

\appendix

\section{What remains to be done}
\label{app-a}

To derive a tractable least action principle for a fluid subject to $ \epsilon$-\MAG, we need to
\begin{enumerate}[\quad(1)]
\item
give a more explicit formula for the least action principle \eqref{eq-40};
 \item
as \eqref{eq-40} depends on $k$, look at its limit as $k$ tends to infinity when the limiting source measure $\lim _{ k\to \infty} \lk=\1 _{ D}\,\Leb$ is the normalized volume of some set $D$, see \eqref{eq-07}.
\end{enumerate}
A challenging problem is to
\begin{enumerate}[\quad(3.i)] 
	\item 
	 answer  the questions addressed below \eqref{eq-nls};
	\item
	find a model accounting for the quantum force $ -\gOW_q  \mathcal{I}^\epsilon _t$ as a limit as $N$ tends to infinity of some system of $N$ interacting particles.
\end{enumerate}
Of course,
\begin{enumerate}[\quad (4)]
\item
let $\epsilon$ tend to zero.
\end{enumerate}
Finally, with Remark \ref{rem-02}-(i) in mind,
\begin{enumerate}[\quad(5)]
\item
translate these results from $\Rd$ to the torus $ \mathbb{T}^d$.
\end{enumerate}
This list only pertains to mathematics. Translating these mathematical results -- or more likely some modified versions of them -- into meaningful physics is an open challenge.

\[
    \begin{tikzpicture}[
node distance = 5mm and 7mm,
   arr/.style = {-Triangle},
   box/.style = {rectangle, draw, semithick,
                 minimum height=9mm, minimum width=17mm,
                 fill=white, drop shadow},
                        ]
\node (n1) [box] {Brownian $k$-mappings};
\node (n2) [ below left=of n1] {$\eta\to 0$};
\node (n3) [below right=of n1] {$N\to \infty$ \& splitting};
\node (n5) [below=of n3] {quantum force};
\node (n7) [box, below=of n5] {$ \epsilon$-\MAG\ for $k$-fluids};
\node (n4) [box, below=of n2] {$ \epsilon$-\MAG\ for $k$-mappings};
\node (n6) [below=of n4] {$ \epsilon\to 0$};
\node (n8) [box, below=of n6] {\MAG\ for $k$-mappings};
\node (n9) [below=of n7] {$ \epsilon\to 0$};
\node (n11) [box, below=of n9] {\MAG\ for $k$-fluids};
\node (n13) [below=of n11] {marginal projection \eqref{eq-40-t} \& $k\to \infty$};
\node (n15) [box, below=of n13] {\MAG\ for fluids};

\draw   (n1) -| (n2);
\draw   (n1) -| (n3);
\draw[arr]   (n2) -- (n4);
\draw   (n3) -- (n5);
\draw[arr]   (n5) -- (n7);
\draw   (n4) -- (n6);
\draw[arr]   (n6) -- (n8);
\draw   (n7) -- (n9);
\draw[arr]   (n9) -- (n11);
\draw   (n11) -- (n13);
\draw[arr]   (n13) -- (n15);
\end{tikzpicture}
\]

The left-hand side of this flow chart corresponds to the Ambrosio-Baradat-Brenier particle system \eqref{eq-21}-\eqref{eq-22}. Its right-hand side corresponds to the approach of the present article. We stopped at $ \epsilon$-\MAG\ for $k$-fluids, using mathematics that need to be consolidated later at several points.  The mathematics of the quantum force arrow are still to be discovered.  The two lowest arrows on the right also remain to be explored.

\section{Least action principle}
\label{app-b}

This short note is about basic calculus of variations with no emphasis on rigorous derivations. A good reference  is the textbook by Gelfand and Fomin \cite{GF63} which doesn't seek mathematical rigor either.
The action functional is defined by 
\begin{equation*}
\mathcal{A}(\omega)= \int _{ t_0}^{ t_1} L(t, \omega_t, \dot \omega_t)\,dt
\end{equation*}
where $ \omega=( \omega_t) _{ t_0\le t\le t_1}$ is a regular $\Rn$-valued path with time derivative $\dot \omega$. The function 
$$
L: (t,q,v)\in [t_0,t_1]\times\ZZZ\mapsto L(t,q,v)\in\RR
$$ 
is called the Lagrangian. It is  assumed to be sufficiently differentiable. 
\\
The first variation  of $ \mathcal{A}$ at $ \omega$ in the direction $\eta$ is
\begin{align*}
\mathrm{d} \mathcal{A} _{ \omega}(\eta)
	:= \lim _{ h\to 0} h ^{ -1} [ \mathcal{A}( \omega+h\eta)- \mathcal{A}( \omega)]
\end{align*}
and, as a definition, a critical path $ \omega$ of $ \mathcal{A}$ satisfies $\mathrm{d} \mathcal{A} _{ \omega}(\eta)=0$ for all $\eta.$

\begin{theorem}
Any critical path $ \omega$ of $ \mathcal{A}$ solves the  Euler-Lagrange equation
\begin{align*}
\partial_q L(t, \omega_t,\dot \omega_t)- \frac{d}{dt}\{\partial_v L(t, \omega_t, \dot \omega_t)\}=0,\quad t_0\le t\le t_1.
\end{align*}
\end{theorem}

\begin{proof}
 The Taylor expansion of $L$ leads us to
\begin{eqnarray*}
\mathrm{d} \mathcal{A} _{ \omega}(\eta)
	&=& \lim _{ h\to 0} h ^{ -1} \int _{  t_0}^{ t_1} [L(t, \omega_t+ h\eta_t,\dot \omega_t+ h\dot \eta_t)-L(t, \omega_t, \dot \omega_t)]\,dt\\
	&=&  \int _{ t_0} ^{ t_1} [\partial_q L(t, \omega_t,\dot \omega_t)\cdot  \eta_t +\partial_v L(t, \omega_t, \dot \omega_t)\cdot \dot \eta_t]\, dt\\
		&=& \partial_v L(t_1, \omega _{ t_1},\dot \omega _{ t_1})\cdot \eta(t_1)
		-\partial_v L(t_0, \omega _{ t_0},\dot \omega _{ t_0})\cdot \eta(t_0)\\
		&& \hskip 4cm +  \int _{ t_0} ^{ t_1} \big[\partial_q L(t, \omega_t,\dot \omega_t)- \frac{d}{dt}\{\partial_v L(t, \omega_t, \dot \omega_t)\}\big]\cdot  \eta_t \, dt
\end{eqnarray*}
where  last identity is obtained by integrating by parts.
Since  $\eta$ is arbitrary, it is necessary that along any critical path of  $ \mathcal{A}$  the integrand $\partial_q L(t)- \frac{d}{dt}\{\partial_v L(t)\}$ vanishes.
\end{proof}

The least action problem is
\begin{align*}
\inf\{ \mathcal{A}( \omega); \omega: \omega _{ t_0}=a, \omega _{t_1 }=b\}
\end{align*}
where the endpoint positions $a$ and $b$ are prescribed.  Under these  constraints, the variation $\eta$  must verify $\eta _{ t_0}=\eta _{ t_1}=0,$ so that
\begin{align*}
\mathrm{d} \mathcal{A} _{ \omega}(\eta)
	=\int _{ t_0} ^{ t_1} \big[\partial_q L(t, \omega_t,\dot \omega_t)- \frac{d}{dt}\{\partial_v L(t, \omega_t, \dot \omega_t)\}\big]\cdot  \eta_t \, dt.
\end{align*}
Of course, any minimizer is critical, so that any solution of the least action problem solves the Euler-Lagrange equation.
\\
An important example in classical mechanics is given by the Lagrangian
\begin{align*}
L(t,q,v)= m|v|^2/2- U(t,q),
\end{align*}
because the corresponding Euler-Lagrange equation is the standard equation of motion
\begin{align*}
m \ddot \omega_t=-\nabla U_t( \omega_t).
\end{align*}
A least action \emph{principle} is a law of nature which stipulates that the trajectory of the physical system solves some least action \emph{problem}. As is customary, although it is slightly incorrect,  we keep saying that this system solves a least action principle.

\section{Approximating Poisson by Monge-Ampère} \label{app-MA-P}

The Monge-Ampère equation
\begin{align}\label{eq1}
 \rho=\rb\det\Big( \mathbb{I}+\Hess\big( [4\pi G\,\rb] ^{ -1} \varphi\big)\Big)
\end{align}
writes as
\begin{align}\label{eq3}
\rho=\rb (1+ \alpha)(1+ \beta)(1+ \gamma)
\end{align}
where 
\begin{align*}
\alpha= \frac{a}{4\pi G\,\rb}\ ,\qquad
\beta= \frac{b}{4\pi G\,\rb}\ ,\qquad
\gamma= \frac{c}{4\pi G\,\rb},
\end{align*}
and $a, b$ and $c$ are the eigenvalues of $\Hess \varphi.$
Here, $G$ is the universal constant of gravitation and $\rb$ is the background density of matter. 

\begin{proposition}
\begin{enumerate}[(a)]
\item
The Monge-Ampère equation \eqref{eq1} approximates the Poisson equation
\begin{align}\label{eq5}
\frac{ \Delta \varphi}{4\pi G}= \rho
\end{align}
if,   up to a permutation of the terms $ \alpha, \beta$ and $ \gamma,$ 
\begin{align}\label{eq2}
 \alpha\longrightarrow \infty,\quad \beta\longrightarrow 0,\quad  \gamma\longrightarrow 0.
\end{align}

\item
Similarly, \eqref{eq1} approximates
the cosmological Poisson equation
\begin{align}\label{eq8}
\frac{ \Delta \varphi}{4\pi G}= \rho- \rb
\end{align}
if,   up to a permutation of the terms $ \alpha, \beta$ and $ \gamma,$  
\begin{align}\label{eq9}
\beta\longrightarrow 0,\quad  \gamma\longrightarrow 0.
\end{align}
\end{enumerate}
\end{proposition}

\begin{proof}
Let us start showing (a).
Expanding the product in \eqref{eq3} and noting that $ \Delta \varphi=a+b+c= (4\pi G\ \rb)( \alpha+ \beta+ \gamma),$ we obtain
\begin{align*}
\rho=  \frac{ \Delta \varphi}{4\pi G} + \rb\ (1+ \alpha \beta+ \alpha \gamma+ \beta \gamma+ \alpha \beta \gamma)
	= \frac{ \Delta \varphi}{4\pi G} \  \Big(1+ 
	\frac{1+ \alpha \beta+ \alpha \gamma+ \beta \gamma+ \alpha \beta \gamma}{ \alpha+ \beta+ \gamma}\Big).
\end{align*}
Hence, \eqref{eq1} approximates \eqref{eq5} if 
\begin{align}\label{eq4}
\frac{1+ \alpha \beta+ \alpha \gamma+ \beta \gamma+ \alpha \beta \gamma}{ \alpha+ \beta+ \gamma} \longrightarrow 0.
\end{align}
In the generic case where the numerator does not vanish, this implies that there is only one leading factor, i.e.\ up to a permutation of the terms $ \alpha, \beta$ and $ \gamma,$ we have
$ \alpha\rightarrow \infty,$  $\beta/ \alpha\rightarrow 0,$ and  $ \gamma/ \alpha\rightarrow 0.$
Together with \eqref{eq4}, this implies \eqref{eq2}.
\\
Let us prove (b). In this case, 
 \eqref{eq1} approximates \eqref{eq8} if 
\begin{align*}
\frac{\alpha \beta+ \alpha \gamma+ \beta \gamma+ \alpha \beta \gamma}{ \alpha+ \beta+ \gamma} \longrightarrow 0,
\end{align*}
which happens when \eqref{eq9} holds.
\end{proof}

Going back to \eqref{eq3}, we observe that in  regime \eqref{eq2}, it is necessary that 
\begin{align*}
\rho/\rb \approx \alpha\longrightarrow \infty.
\end{align*}
\subsubsection*{Conclusion}

The Monge-Ampère equation \eqref{eq1}  approximates the cosmological Poisson equations \eqref{eq8} in the regions where  the force field has parallel force lines. This happens in the vicinity of a very large flat wall. The problem is well described by the 1D profile of the distribution of matter in the direction which is orthogonal to the wall. If in addition $ \rho/\rb$ is very large, then \eqref{eq1}  approximates the  Poisson equation  \eqref{eq5}.

\subsection*{\MAG\ in 1D}

In this case, Monge-Ampère and Poisson equations \emph{almost} coincide. Recall that \MAG\ necessitates a target profile   $ \lambda$ with the same  finite mass as the fluid profile $ \rho.$ In particular, \eqref{eq1} is not valid when considered on the whole space $\RR^3,$ but only on the torus $(\RR/L \mathbb{Z})^3$ with $ \lambda(d y)=ML ^{ -3}\1 _{ [0,L]^3}( y)\,d y$ which corresponds to $\rb=ML ^{ -3}$ when $M$ is the total mass of the fluid.

 Let $ \rho$ be a 1D  profile of matter with mass $m,$ to be interpreted as an \emph{area density} since we project a 3D model onto a 1D one.  It is described by a  probability measure $ \rho(dx)$ on the circle $\RR/L \mathbb{Z}$ with perimeter $L$. The target measure $ \lambda(d y) $ is also such that $ \lambda(\RR/L \mathbb{Z})=m $.    While $L$ tends to infinity,  we assume that $\rho$ is supported on a fixed bounded subset $[a,b]$. Therefore, at the limit $L\to \infty,$ one recovers  \MAG\ when looking at $[a,b]$ as a subset of the segment $[-L/2,L/2],$ instead of the circle, i.e.\ without identifying the endpoints. Indeed, on the very large  circle seen as a clock, the small area $[a,b]$ situated near noon looks like a Dirac  point mass at noon, so that up to a tiny fraction of order $L ^{ -1},$ the geodesics arriving at  $[a,b]$ are clockwise when starting from a point between $6$ pm and midnight = noon, and anti-clockwise when starting from a point between noon and $6$ pm. Hence, as $L$ tends to infinity, cutting the circle at $6$ pm does not change anything.
  \\
In this setting, the analog of \eqref{eq1} is
$	
\rho(x)=\1 _{\{ T  (x)\in [-L/2,L/2]\}}( m  L ^{ -1}+ \varphi''(x)/(4\pi G))
$	
where we used $\rb=mL ^{ -1}$ and  $T $ is the optimal transport map from the matter profile   $ \rho(dx)$ to the target measure $ \lambda(d y)=m  L ^{ -1}\1 _{ [-L/2,L/2]}( y)\,d y$.  It  is well known that when $ \rho$ is absolutely continuous 
\begin{align*}
T (x)= [F _{  \lambda} ^{ -1} \circ F_ {  \rho}](x)
	=L( m  ^{ -1}F_ {   \rho}(x)-1/2)\in [-L/2,L/2],\quad x\in \supp \rho,
\end{align*}
 where $F_ {   \lambda} ( y)=   \lambda((- \infty,  y])$ and $F_ {   \rho}(x)=   \rho((- \infty,x])$ are the cumulative functions of $ \lambda$ and $ \rho.$ Therefore, \eqref{eq1} is
\begin{align*}
\rho(x)= m  L ^{ -1}+ \varphi''(x)/(4\pi G),\quad x\in \supp \rho,
\end{align*}
and  \MAG's force field 
 $
4\pi G m L ^{ -1}(x-T  (x))=4\pi G(m L ^{ -1} x+m /2-F _ \rho(x)),
$
tends to
 \begin{align}\label{eq6}
 - \varphi'(x)=4\pi G(m /2-F _ \rho(x)),\quad x\in \supp \rho,
 \end{align}
as $L$ tends to infinity. 
\\
The corresponding Poisson equation \eqref{eq5} is
\begin{align*}
\rho(x)=  \varphi''(x)/(4\pi G),\quad x\in \RR.
\end{align*}
It implies that $ \varphi'(x)= 4\pi G F _ \rho(x)+c$ for some real $c.$ Remark that this last equation also holds when $ \rho$ is not assumed to be absolutely continuous. The  constant $c$  is specified by the physical requirement that $ \varphi'( - \infty)=- \varphi'( \infty).$ That is $c= \varphi'(- \infty)=- \varphi'( \infty)=-4\pi G m -c.$ Hence $c= -2\pi G m$ and 
\begin{align} \label{eq7}
 - \varphi'(x)=4\pi G(m /2-F _ \rho(x)) \in [-2\pi Gm,2\pi Gm],\quad x\in \RR.
\end{align}
Comparing with \eqref{eq6}, we remark that the truly physical solution  \eqref{eq7} extends in a natural way the \MAG\ solution \eqref{eq6}.
\\
Recall that an infinite flat wall with area density $ \sigma$ generates a constant force field with intensity $2\pi G  \sigma$ pointing orthogonally towards the wall. This is  in accordance with formula \eqref{eq7}, applied with $ \sigma(dx)= \rho(dx),$  exhibiting a  force field which is constant on each connected component of the complementary subset of  the support of $ \rho.$

\section{\OW-manifold}\label{app-OW}

\begin{definition}[$\gcb( \mu)$] \label{def-H}
The space $\gcb( \mu)$ is the closure in $L^2 _{ \Rd}( \mu)$ of the space\\ $ \left\{\nabla u; u\in C^1_c(\Rd)\right\} $ of all regular gradient vector fields.
\end{definition}
As will be seen at \eqref{eq-78}, this space is the natural candidate to be the tangent vector space at $ \mu$  in the Otto-Wasserstein manifold.

\begin{proposition}\label{res-01} 
Suppose that the  (normalized) density  $ \mu_t=(Y_t)\pf P$    is the time marginal at time $t$ of some path measure $P\in\PO$ which is supported by absolutely continuous sample paths, i.e.
\begin{align*}
dY_t=\dot Y_t\,dt,\quad P\as,
\end{align*}
where $(Y_t) _{ t_0\le t\le t_1}$ stands for  the canonical process on the path space \[\OO:= C([t_0,t_1],\Rd)\] and $\dot Y_t$ is some random vector, possibly depending (a priori) on the whole history of the path and  satisfying
\begin{align*}
E_P \int _{ t_0} ^{ t_1} |\dot Y_t|^2\, dt < \infty.
\end{align*}
 Then, there exists some measurable vector field $v_t(y)$ such that for almost all $t$, $v_t$  belongs to   $\gcb( \mu_t)$, 
 \[
 \int _{ \Rd\times[t_0,t_1]} |v_t(y)|^2\,\mu_t(dy)dt< \infty,
 \] 
 and the continuity equation
  \begin{align}\label{eq-26b}
 \partial_t \mu_t+\nabla \scal (\mu_t v_t)=0
 \end{align}
holds in the following weak sense: For any function $f$ in $C ^{ 1}_c(\Rd)$  and any $t_0\le t_*\le t_1,$ 
\begin{align}\label{eq-27}
\IZ f &\, d \mu _{t_*}- \IZ f \, d \mu _{ t_0}
	=\int _{ \Rd\times[t_0,t_*]} \nabla f(y)\scal  v_t(y)\, \mu_t(dy)dt. 
\end{align}
\end{proposition}

\begin{proof}
Take  any function $g$ in $C ^{ 1,1}_c(\Rd\times [t_0,t_*])$  with $t_0\le t_*\le t_1.$ Then, 
\begin{align*}
\IZ g _{ t_*}&\, d \mu _{t_*}- \IZ g _{ t_0}\, d \mu _{ t_0}
	= E_P \int _{ t_0} ^{ t_*}[ \partial_tg_t(Y_t)+ \nabla g_t(Y_t)\scal\dot Y_t]\,dt\\
	&= E _{ \bar P}( \partial_t g(\tau,Y_ \tau)+\nabla g _ \tau(Y_ \tau)\scal \dot Y_ \tau)
	=E _{ \bar P}E _{ \bar P} \big[  \partial_t g(\tau,Y_ \tau)+ \nabla g_ \tau(Y_ \tau)\scal \dot Y_ \tau\mid (Y_ \tau, \tau) \big]\\
	&=E _{ \bar P}\big[  \partial_t g(\tau,Y_ \tau)+ \nabla g_ \tau(Y_ \tau)\scal \tilde v_ \tau(Y_ \tau)\big]
	=\int _{ \Rd\times [t_0,t_*]} \big[  \partial_t g_t(y)+ \nabla g_t(y)\scal \tilde v_t(y)\big]\, \mu_t(dy)dt\\
	&=\int _{ \Rd\times [t_0,t_*]} \big[ \partial_t g_t(y)+ \nabla g_t(y)\scal  v_t(y)\big]\, \mu_t(dy)dt, \
\end{align*}
where  $\bar P( d\omega dt):= P(d \omega)dt$, the canonical time is $ \tau$, 
$	
\tilde v_t(y):= E _{ \bar P} \left( \dot Y_ \tau\mid  Y_ \tau=y, \tau=t\right),
$	
and $v$ is the orthogonal projection in $L^2 _{ \Rd}(\mu_t(dy)dt)$ of $\tilde v$ on the closure of the space $ \left\{ \nabla g; g\in C^{1,1}_c(\Rd\times [t_0,t_*])\right\} $. 
This implies   \eqref{eq-27} and the joint measurability of $(t,y)\mapsto v_t(y)$.
\end{proof}

The basic insight of Otto calculus \cite{Otto01,AGS05,Vill09} is to interpret the velocity field  $v_t\in \gcb( \mu_t)$ appearing at Proposition \ref{res-01} as a tangent vector at $\mu_t$   in some Riemannian-like manifold $ \mathrm{P}_2\subset \mathrm{P}(\Rd)$,  called the Otto-Wasserstein manifold, \OW-manifold for short. We denote this velocity by
\begin{align}\label{eq-78}
\dot\mu_t:=v_t \in \mathrm{T} _{ \mu_t} \mathrm{P}_2\subset   \gcb( \mu_t)
\end{align}
where $ \mathrm{T} _{ \mu} \mathrm{P}_2$ stands for the tangent space of $ \mathrm{P}_2$ at $ \mu.$ In particular, the continuity equation \eqref{eq-26b} writes as
\begin{align}\label{eq-26c}
\partial_t \mu+\nabla\scal( \mu\dot\mu)=0.
\end{align}
We emphasize that the vertical variation $ \partial_t \mu_t$ differs from the horizontal variation $\dot \mu_t.$ In particular, $ \partial_t \mu_t$ is a scalar while $\dot \mu_t$ is a vector.

\subsection*{Benamou-Brenier formula}

The relation between quadratic optimal transport and Otto calculus is best illustrated by the Benamou-Brenier formula
\begin{align}\label{eq-BB}
\begin{split}
W_2^2( \alpha,\beta)
	&\overset{(\textrm{i})}=\inf _{ ( \mu,u)} \int_0^1 \Big(\IZ |u_t(y)|^2\, \mu_t(dy)\Big)\,dt\\
	&\qquad \overset{(\textrm{ii})}=\inf _{ ( \mu)} \int_0^1  \Big(\IZ |\dot \mu_t(y)|^2\, \mu_t(dy)\Big)\,dt
	\overset{(\textrm{iii})}=\inf _{ ( \mu)}\int_0^1 \|\dot\mu_t\|^2 _{ \mu_t}\,dt,
\end{split}
\end{align} 
where $ \alpha	, \beta\in \mathrm{P}(\Rd)$ and  
\[W_2^2( \alpha, \beta):= \inf _\pi \int _{ \Rd\times\Rd} |y-x|^2\,\pi(dxdy)\] is the optimal quadratic transport cost from $ \alpha$ to $ \beta.$ 
\\
Let us comment on these infimums:
\begin{enumerate}[-]
\item
in the expression $\inf_\pi,$ $\pi\in \mathrm{P}(\Rd\times\Rd)$ runs through all the couplings of $ \alpha$ and $ \beta,$ that is: $\pi(dx\times\Rd)=  \alpha(dx)$ and $\pi(\Rd\times dy)= \beta(dy);$
\item
in the expression $\inf _{ ( \mu)},$ the infimum runs through all the paths $( \mu)=( \mu_t) _{ 0\le t\le 1}$ in $\PZ$ such that $ \mu_0 = \alpha$ and $ \mu_1= \beta;$ 
\item
in the expression $\inf _{ ( \mu,u)},$ the infimum runs through all the paths  $( \mu,u)=( \mu_t,u_t) _{ 0\le t\le 1}$ where $ \mu_0= \alpha$, $ \mu_1= \beta$   and the  $u_t$'s are vector fields such that the continuity equation $ \partial_t \mu+\nabla\cdot( \mu u)=0$ holds in the weak sense.
\end{enumerate}
Identity (i) is the Benamou-Brenier formula. Identity (ii) relies on Proposition \ref{res-01} which states that one can replace $u_t$ in the continuity equation $ \partial_t \mu+\nabla\cdot( \mu u)=0$ by $ \dot \mu_t$ which, as an element of  $\gcb( \mu_t),$ minimizes $\IZ |u_t|^2\, d \mu_t$ (Hilbertian projection onto  $\gcb( \mu_t)$). The last equality (iii) directly follows from the notation
\begin{align}\label{eq-49}
\| \dot \mu_t\|^2 _{ \mu_t}:=\IZ |\dot \mu_t(y)|^2\, \mu_t(dy)
\end{align}
where $\| \dot \mu_t\|^2 _{ \mu_t}$ should be interpreted as the analogue of the squared Riemannian norm of a tangent vector at $ \mu_t.$ 

\subsection*{Otto's heuristics} \label{sec-OW2}

In the \OW-manifold, the velocity at time $t$ of a moving profile of matter $(q_t)$ is the vector field $\dot\qq_t\in  \gcb(q_t)$ of the continuity equation \eqref{eq-26c}: $ \partial_tq+\nabla \scal(q\dot\qq)=0.$ Staying at a heuristic level mainly consists of replacing $\dot\qq_t\in  \gcb(q_t)$ by   \[\dot\qq_t=\nabla \theta_t\] for some $C ^{ 1,2}$-regular function $(t,x)\mapsto \theta(t,x).$

\subsection*{Acceleration}
 If this wishful thinking is realized, the acceleration is given by
\begin{align}\label{eq-OWacc}
\nabla ^{ \mathrm{OW}} _{ \dot\qq_t}\dot\qq_t
	= ( \partial_t+ \dot\qq_t\scal\nabla) \dot\qq_t
	=\nabla \big( \partial_t \theta_t+\uud |\nabla \theta_t|^2 \big).
\end{align}
We see that, at least informally,  $\nabla ^{ \mathrm{OW}} _{ \dot\qq_t}$ is identified with the convective derivative $ \partial_t+\dot\qq_t\scal \nabla.$ In fact, this works fine when calculating an acceleration, but in the general case $\nabla ^{ \mathrm{OW}} _{\nabla \alpha} \dot\qq_t$  should be  identified with
$\proj _{  \gcb(q_t)}\big[ ( \partial_t+\nabla \alpha\cdot\nabla)\dot\qq_t\big]$ where $\proj _{  \gcb(q_t)}$ is 
 the orthogonal projection onto the space of \emph{gradient} vector fields $  \gcb(q_t).$
 
For a presentation of Otto's heuristics, see the chapter entitled \emph{Otto calculus} in \cite{Vill09}. Rigorous material  is presented in  \cite{LV09,AGS11,Gig12}.

\section{Current velocity as $ \epsilon\to 0$}
\label{app-d}

The pointwise limit as $ \epsilon$ tends to zero  of the current velocity $ \epsilon  \nabla\log \sqrt{r^\epsilon_s}$ of $\Xe$, recall \eqref{eq-45},  does not lead to the right expression for a $ \Gamma$-limit, see  \eqref{eq-24}. Nevertheless, it is interesting to compute it as a first step.

Let $ \gamma\in\PRdk$ be the initial law of $\Xe.$ 
 The Gaussian nature of \[r^\epsilon_s= \gamma \ast \mathcal{G}\big(0, \epsilon (s-s_0)\Id\big)\] permits us to obtain the explicit expression
\begin{align*}
 \epsilon  \nabla\log \sqrt{r^\epsilon_s}(\yy)
 	= -Z_ \epsilon ^{ -1}(\yy) \int _{ \Rdk}  \frac{ (\yy-\xx)}{2(s-s_0)}  \exp\Big(- \frac{|\yy-\xx|^2}{2 \epsilon (s-s_0)}\Big)\,  \gamma(d\xx)
\end{align*}
with the normalizing constant $Z_ \epsilon(\yy):=\int _{ \Rdk}    \exp\big(- {|\yy-\xx|^2}/({2 \epsilon (s-s_0)})\big)\,  \gamma(d\xx).$ It follows from the  Laplace principle that 
\begin{align}\label{eq-y53}
\lime  \epsilon  \nabla\log \sqrt{r^\epsilon_s}(\yy)
	=- \frac{\yy-\txa(\yy) }{2(s-s_0)}, \qquad y\in\ZZ,\ t>0,
\end{align}
where $\txa(\yy)$ is the orthogonal projection of $\yy$ onto the support  of the initial measure $  \gamma$, whenever this projection  is unique. In the general case where the projection might not be unique, we must take
\begin{align*}
\txa(\yy)=\int _{ \Rdk}  \xx\, \gamma(d\xx\mid \proj_{\supp  \gamma}(\yy))
\end{align*}
where  $ \gamma(\sbt\mid \proj_{\supp  \gamma}(\yy))$ is the probability $  \gamma$ conditioned on the set $\proj_{\supp  \gamma}(\yy)$ of all the closest  points to $\yy$ in the support of $  \gamma.$ This result is a special case of next  Lemma \ref{res-y06} applied with $c(x,y)=|y-x|^2/(2(s-s_0))$ and $ \alpha= \gamma.$

\begin{lemma}\label{res-y06}
For each $ \epsilon>0$, let \[\nu ^{ \epsilon}(y)=\int_{\RR^n} k^ \epsilon(x,y)\,  \alpha(dx),\quad y\in\RR^n,\] be a probability density with kernel 
\[
k^ \epsilon(x,y)=\Lambda^ \epsilon(x) ^{ -1}\exp( -c(x,y)/ \epsilon)
\]
 where $  \alpha\in\mathrm{P}(\RR^n)$ and   $c:\RR^n\times\RR^n\to\RR$ is a regular  function such that for any $x\in\RR^n,$ the normalizing constant  $\Lambda^ \epsilon(x)=\int_{\RR^n} e ^{ -c(x,y)/ \epsilon}\, dy< \infty$ is finite.  \\ 
 We also assume  that for each $x$, 
 \begin{enumerate}[(i)]
 \item
  the function $c(x,\sbt)$ admits a unique minimizer $y_x$,  
   \item
   $ \det\nabla_y^2c(x,y_x)>0,$
  \item
   $\inf_y c(x,y)=c(x,y_x)=\inf c$   doesn't depend on $x$.
 \end{enumerate}
  Then, for any $y$,
\begin{align*}
\lime \epsilon\nabla\log \nu^ \epsilon(y)
	=- \frac{\int_{\RR^n} \nabla_y c(x,y)\   \sqrt{\det\nabla_y^2c(x,y_x)}\, \alpha(d\xx\mid \proj^c_{\supp  \alpha}(y))}
		{\int_{\RR^n}  \sqrt{\det\nabla_y^2c(x,y_x)}\, \alpha(d\xx\mid \proj^c_{\supp  \alpha}(y))}
\end{align*}
where  $ \alpha(\sbt\mid \proj^c_{\supp  \alpha}(y))$ is the probability $  \alpha$ conditioned on the set 
\[\proj^c_{\supp  \alpha}(y):= \left\{x\in\supp  \alpha: c(x,y)=\inf _{ x'\in\supp  \alpha} c(x',y)\right\} \] 
of all the points in the support of $  \alpha$ realizing the minimal value of $c(\sbt,y)$.
\end{lemma}

\begin{proof}
Be aware that $y$ is fixed once for all.
Let us disintegrate $  \alpha$ along the values of the function $c(\sbt,y):$ 
\begin{align*}
  \alpha(dx)= \int_{u\in\RR}   \alpha(dx\mid c(\sbt,y)=u)\ (c(\sbt,y)\pf  \alpha)(du).
\end{align*}
With this in our hands,
\begin{align*}
\epsilon\nabla&\log \nu^ \epsilon(y)
	= \epsilon\frac{\int_{\RR^n}\nabla_y k^ \epsilon(x,y)\,  \alpha(dx)}{\int_{\RR^n}  k^ \epsilon(x,y)\,  \alpha(dx)}
	=- \frac{\int_{\RR^n}\nabla_y c(x,y) k^ \epsilon(x,y)\,  \alpha(dx)}{\int_{\RR^n}  k^ \epsilon(x,y)\,  \alpha(dx)}\\
	&=-\frac{\int _{ u\in\RR} e^{-u/ \epsilon}\Big[\int _{ x\in\RR^n}\nabla_y c(x,y) \Lambda^ \epsilon(x) ^{ -1}\,  \alpha(dx\mid c(\sbt,y)=u)\Big](c(\sbt,y)\pf  \alpha)(du)}
	{\int _{ u\in\RR} e^{-u/ \epsilon}\Big[\int _{ x\in\RR^n} \Lambda^ \epsilon(x) ^{ -1}\,  \alpha(dx\mid c(\sbt,y)=u)\Big](c(\sbt,y)\pf  \alpha)(du)}.
\end{align*}
 We obtain with the standard Laplace method
\begin{align*}
\Lambda^ \epsilon(x)= \left( \frac{ (2\pi \epsilon) ^{ d/2}}{\sqrt{\det\nabla_y^2c(x,y_x)}}+o( { \epsilon ^{ d/2}})\right) \ \exp( -\inf c/ \epsilon),
\end{align*}
where  we use the hypothesis $\inf_y c(x,y)=\inf c.$
This implies that $ \epsilon\nabla\log\nu^ \epsilon(y)$ is equal to
\begin{align*}
-\frac{\int _{ u\in\RR} e^{-u/ \epsilon}\Big[\int _{ x\in\RR^n}\nabla_y c(x,y) \big( \sqrt{\det\nabla_y^2c(x,y_x)}+ o _{ \epsilon\to 0}(1)\big)\,  \alpha(dx\mid c(\sbt,y)=u)\Big](c(\sbt,y)\pf  \alpha)(du)}
	{\int _{ u\in\RR} e^{-u/ \epsilon}\Big[\int _{ x\in\RR^n} \big( \sqrt{\det\nabla_y^2c(x,y_x)}+ o _{ \epsilon\to 0}(1)\big)\,  \alpha(dx\mid c(\sbt,y)=u)\Big](c(\sbt,y)\pf  \alpha)(du)}.
\end{align*}
We conclude noting that the sequence of probability measures $ \displaystyle{ \frac{e ^{ -u/ \epsilon} (c(\sbt,y)\pf  \alpha)(du)}{\int_\RR e ^{ -v/ \epsilon} (c(\sbt,y)\pf  \alpha)(dv)}}$ converges as $ \epsilon$ tends to  zero to the Dirac measure   at $\inf _{ x\in\supp  \alpha}c(x,y). $
\end{proof}

The limit \eqref{eq-y53} is purely static. It does not take the dynamics into account. This is the reason why the $ \Gamma$-limit \eqref{eq-24} of the action is expressed in terms of $\widehat{\proj}_S(\yy ):=\proj _{S(\yy )}(\yy )$, see \eqref{eq-19}, instead of $\txa(\yy)$. More details about this are presented below Figure \ref{fig-y1} at Appendix \ref{app-c}.

\section{Concentration of matter}
\label{app-c}

Let us explain informally, using a 2D analogy,  why the expression \eqref{eq-15} at  Definition \ref{def-03} of \MAG's action functional is reasonable. 

\begin{definition}[$\proj^o_S$]
Let  $\Proj_S(\yy)$ be the set of all the closest points to $\yy$ in $S$ and $\cl\cv(\Proj_S(\yy))$ be its closed convex hull. We define
 $\proj^o_S(\yy)$ as the (unique) element in $\cl\cv(\Proj_S(\yy))$ with minimal norm. 
\end{definition}

\begin{figure}[h]
	\includegraphics[width=6cm]{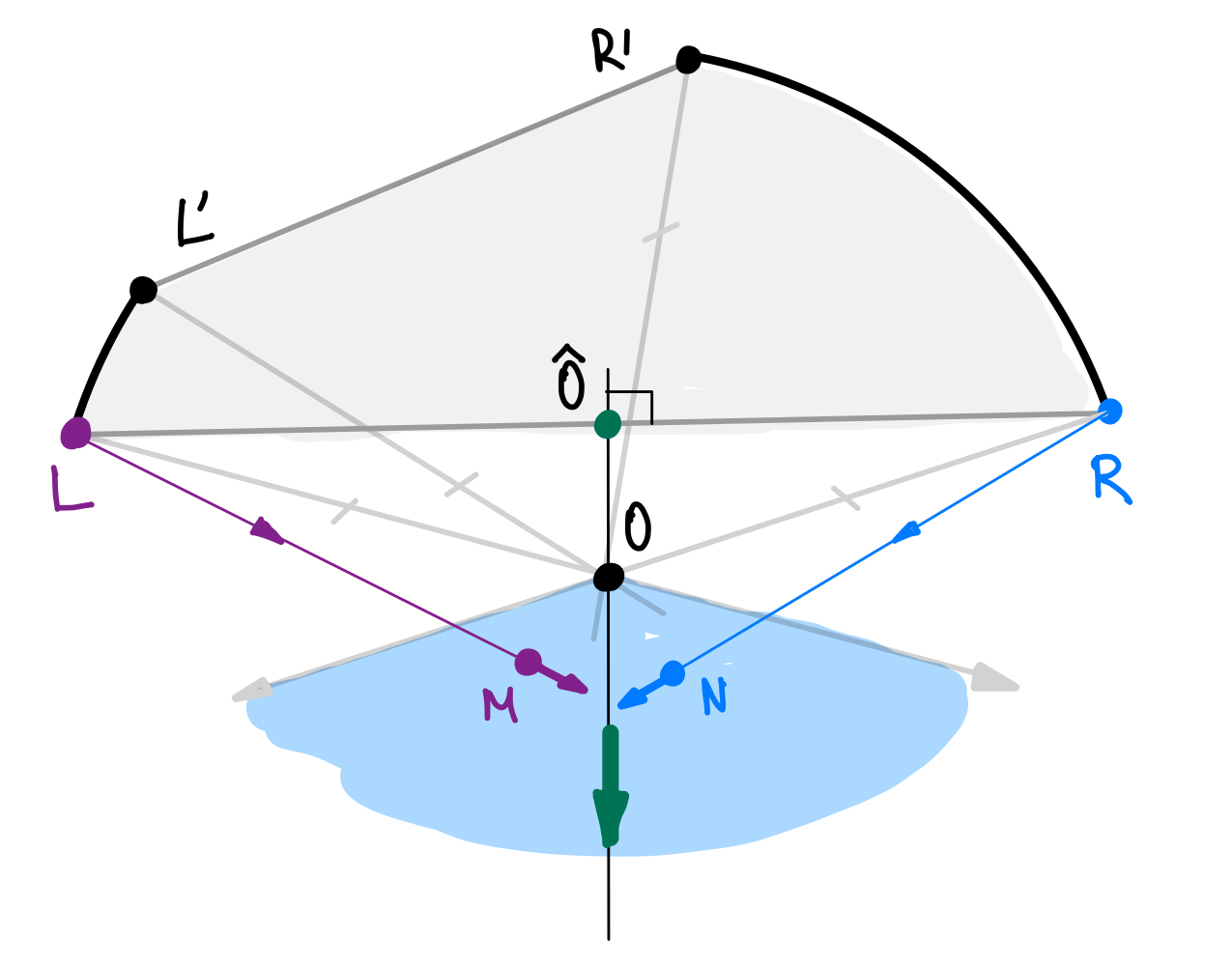}		
	\captionof{figure}{}  \label{fig-y1}
\end{figure}

In Figure \ref{fig-y1},  the point $O$ is the center of the circle passing through the points $L, L', R,  R'$, and the set $S$ is the union of the  arcs: $ \mathrm{arc}_{LL'}$ and   $\mathrm{arc}_{R'R}.$  

We see that $\proj_S^o(M)=\proj_S(M)=L,$ $\proj^o(N)=\proj_S(N)=R$. But $\proj^o(O)=\widehat O$ because $\Proj_S(O)=S$ and $\cl\cv(S)$ is the gray area.

 In order to minimize the kinetic action \eqref{eq-16} the points $M$ and  $N$ must evolve respectively in the direction of the purple and blue arrows. 
\\
The situation for $O$ is a little bit different because its projection on $S$  is the set $\Proj_S(O)=S $.  This implies that $O$ must evolve in the direction of the blue cone. Suppose it moves on the left of the segment bisector of the points $L$ and $R$ which is the line extending $\widehat OO$. Then as happens to $M$, it is instantaneously pushed back onto the right towards the bisector. It can neither move on the right, for the same symmetric  reason. No choice, then: the only option is to move in the direction of the green arrow. It stays on the bisector. 
\begin{center}
\includegraphics[width=6cm]{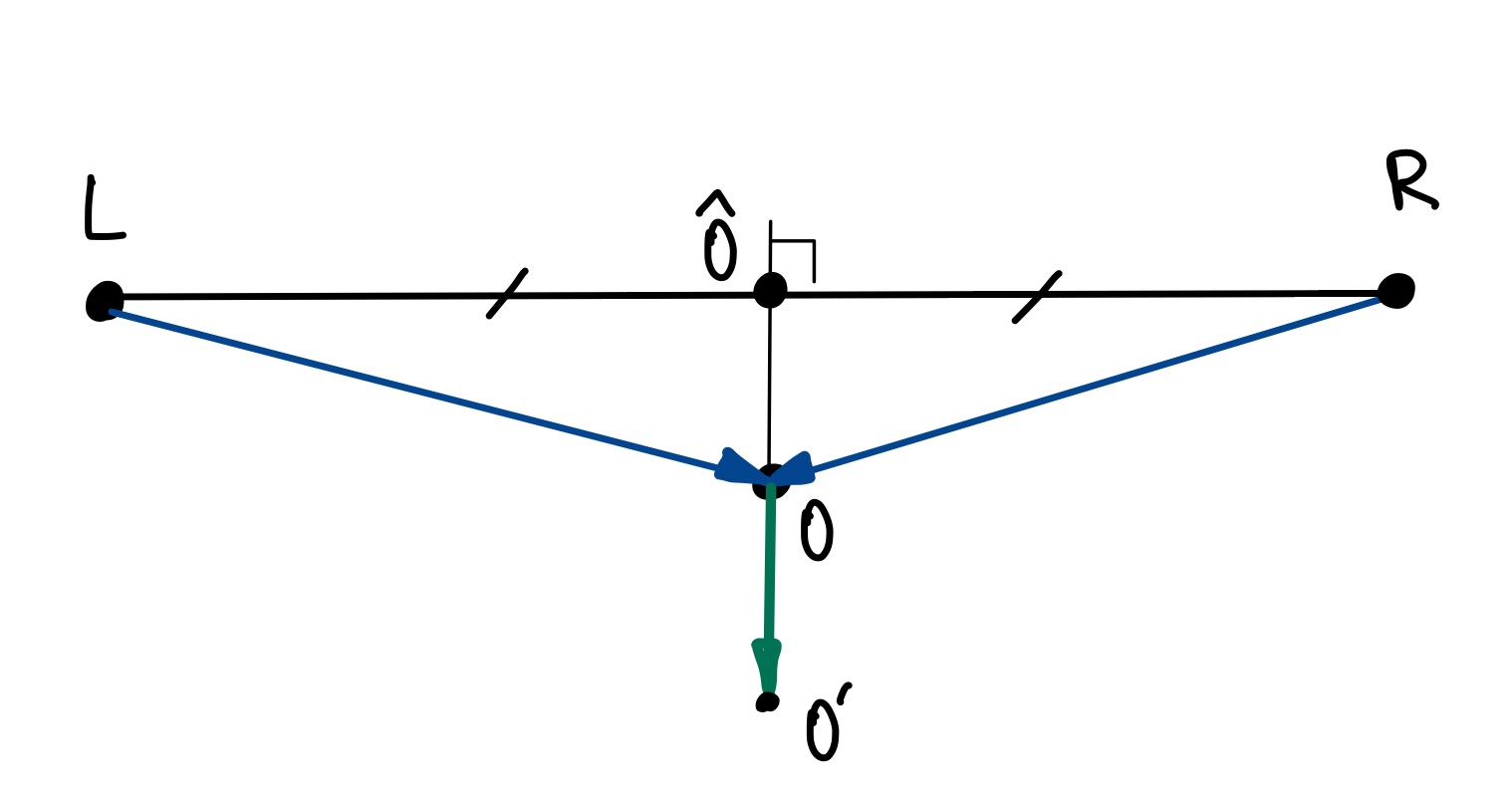}
\end{center}	
By symmetry, the effective velocity $\overrightarrow{O O'}$ is equal to the average of the velocities generated by $L$ on the left and $R$ on the right. That is something proportional to
\begin{align*}
\overrightarrow{OO'}=(\overrightarrow{LO}+\overrightarrow{RO})/2=\overrightarrow{\widehat OO}.
\end{align*}
On the other hand, we have  $\widehat O=\proj _S^ o(O),$ and it is proved at Proposition  \ref{res-05} below that, going back to the notation of Definition \ref{def-03}, we have $\widehat\Phi(0)=0-\widehat\proj_S(0).$
\hfill$\square$

 \subsection*{A convex analytic result}
 
Denoting $ \Psi:=-\Phi$, the action \eqref{eq-16} writes as
\begin{align}\label{eq-A30b}
\IT\uud \|\dot\yy_t+\grad_H\Psi(\yy_t)\|_H^2\ dt,
\end{align}
and any solution of the gradient flow equation
\begin{align}\label{eq-A27}
\dot\yy_t=-\grad_H\Psi(\yy_t), 
\end{align}
is a global minimizer of \eqref{eq-A30b}. 
For any $\yy\in H,$ 
\begin{align*}
\Psi(\yy)= -\inf _{ \xx\in S}\|\yy-\xx\|_H^2/2
	=\iota^* _S(\yy) -\|\yy\|_H^2/2 -r^2/2,
\end{align*}
where $\iota^* _S(\yy):=\sup _{ \xx\in S}\, \xx\cdot\yy,$  is the convex conjugate of $\iota_S(\xx)= \left\{ \begin{array}{ll}
0,&\quad \textrm{if }\xx\in S\\
+ \infty,&\quad \textrm{otherwise}
\end{array}\right.,$ and we used \eqref{eq-112}. 
Since $\iota_S^*$ is convex,  $\Psi$ is a $ \alpha$-convex function (with $ \alpha=-1$). Following \cite[Cor.\,1.4.2, Thm.\,2.4.15]{AGS05}, equation \eqref{eq-A27} admits as a natural extension 
\begin{align}\label{eq-A29}
\dot\yy_t=- \partial^o\Psi(\yy_t), \quad \textrm{for a.e. }t,
\end{align}
\begin{definition}[$ \partial^o\Psi(\yy)$]
We denote by $ \partial^o\Psi(\yy)$ the (unique) element with minimal norm of the  local subdifferential  $ \partial\Psi(\yy)$  defined by
\begin{align*}
\zeta\in \partial\Psi(\yy) \iff
\liminf _{ \yy'\to\yy} \frac{\Psi(\yy')-[\Psi(\yy)+ \langle \zeta,\yy'-\yy\rangle_H ]}{\|\yy'-\yy\|_H}\ge 0,
\qquad \zeta\in H.
\end{align*}
\end{definition}

\begin{proposition}\label{res-05}
For any $ \yy\in H$, the element of $ \partial\Psi(\yy)$   with minimal norm is
\begin{align*}
\partial^o\Psi(\yy)=\proj^o_S(\yy)-\yy.
\end{align*}
\end{proposition}

\begin{proof}
The definition of  $\partial\Psi(\yy)$ clearly implies that it  is a convex set, hence the uniqueness of its element with minimal norm if it is nonempty. But, because $\Psi$ is the sum of the differentiable function $f=-\|\sbt\|^2_H/2$ and the convex function $\iota^*_S,$    its local subdifferentials are nonempty. More precisely,
\[\partial\Psi(\yy)= f'(\yy)+ \partial\iota_S^*(\yy)=-\yy+ \partial\iota_S^*(\yy)\] where 
\begin{align*}
\partial\iota_S^*(\yy)
	= \left\{ \zeta\in H: \iota^*_S(\yy')-\iota^*_S(\yy)- \langle \zeta,\yy'-\yy\rangle \ge0,\ \forall \yy'\in H\right\} \neq \emptyset
\end{align*}
is the (global) subdifferential of $\iota^*_S$ at $\yy.$ The convex conjugate of $ \iota_S^*$  is the convex indicator $ \iota _{ \cl\cv S}$ of the closed convex hull $\cl\cv(S)$ of $S.$ 
Consequently, 
$\zeta \in \partial\iota_S^*(\yy)$ is equivalent to $\yy\in \partial \iota _{ \cl\cv(S)}( \zeta)$, that is: $\zeta\in \cl\cv(S)$ and $\yy$ is in the cone of outer normals of $\cl\cv(S)$ at $ \zeta.$  
\\
In the general case, this property does not imply that $\zeta$ is the  orthogonal projection $\proj _{S}(\yy)$ of $\yy$ on $S.$ But in the present setting where $S$ is a subset of a sphere centered at zero, we obtain 
\[\partial\iota_S^*(\yy)=\cl\cv(\Proj_S(\yy)).\]
In the regular case where $\Proj_S(\yy)$ is reduced to the single point $\proj_S(\yy),$ we clearly obtain: $\partial^o\Psi(\yy)=\proj_S(\yy)-\yy.$ In the general case, we have to replace  $\proj_S(\yy)$ by $\proj^o_S(\yy)$.
 Again, this holds because $S$ is included in a sphere centered at zero. 
\end{proof}

\subsection*{An analogical illustration for the concentration of matter}

In the  figure below illustrating a 2D analogy where the black arcs are $S$ and the gray area is $\cl\cv(S),$ we see that $ \partial\iota_S^*({\red P})=\{\proj_S(P)\},$ $ \partial\iota_S^*({\color{violet} N})=\{{\color{violet} L}\}$, $ \partial\iota_S^*({\color{pink} M})=[L',R']$. These identities are unchanged when replacing $P$ by any red cross point on the picture, $N$ by any purple cross point, and $M$ by any pink cross point. The orange crosses are mapped to ${\color{orange} L'}$ and $\partial\iota_S^*(O)=\cl\cv(S).$
\\
We also see that $\proj^o_S({\color{red}P})=\proj_S(P),$ $\proj^o_S({\color{blue}Q})=\proj_S(Q),$ but $\proj^o_S({\color{pink}M})=M'$ and $\proj^o_S({\color{green}O})=\widehat O.$

\begin{center}
\includegraphics[width=8cm]{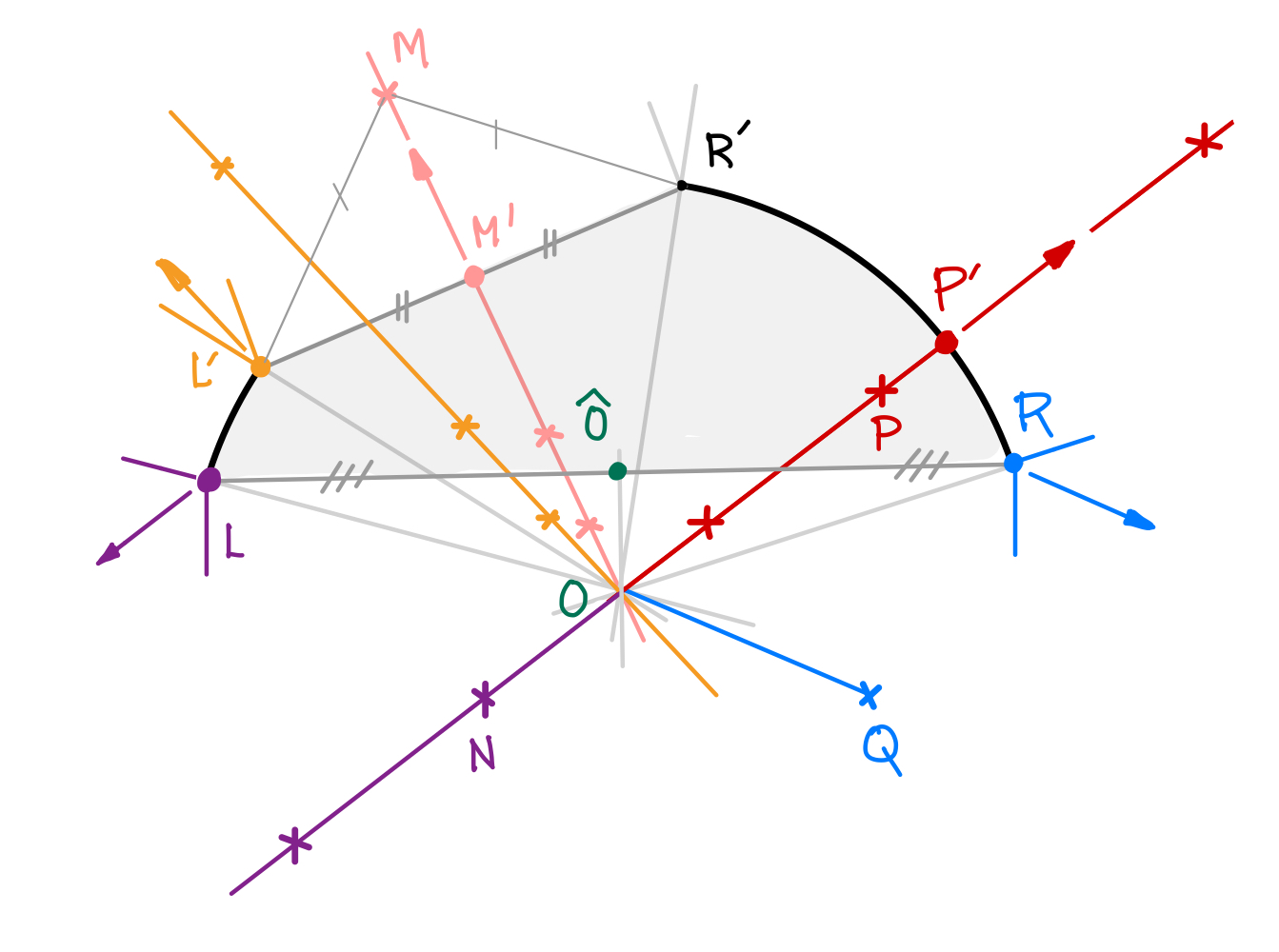}
\end{center}	
Keeping our 2D analogy, the figure below illustrates the motion of a particle moving according to the modified gradient flow  \eqref{eq-A29}. Shocks occur. Successively: $\proj^o_S(P_0)=A,$ $\proj^o_S(P_1)=M',$ $\proj^o_S(P_2=O)=\proj^o_S(P_3)=\widehat O,$  and we see that:
\begin{itemize}
\item
$s\mapsto \textrm{dist}(\yy_t,S)$ increases continuously
\item
 with a discontinuity of its derivative  when passing through $P_1$ and $P_2$.
\item
Furthermore, we observe that at each shock the modulus of the force decreases. Indeed, \eqref{eq-111} tells us that  it is proportional to $|\yy_t-\proj^o_S(\yy_t)|,$ and this quantity jumps from $L'P_1$ to $M'P_1<L'P_1$ at $P_1$, and from $L'O=r$ to $\widehat OO<r$ at $P_2=O.$ Energy  dissipates during each shock.
\end{itemize}
\begin{figure}[h] 
\includegraphics[width=6cm]{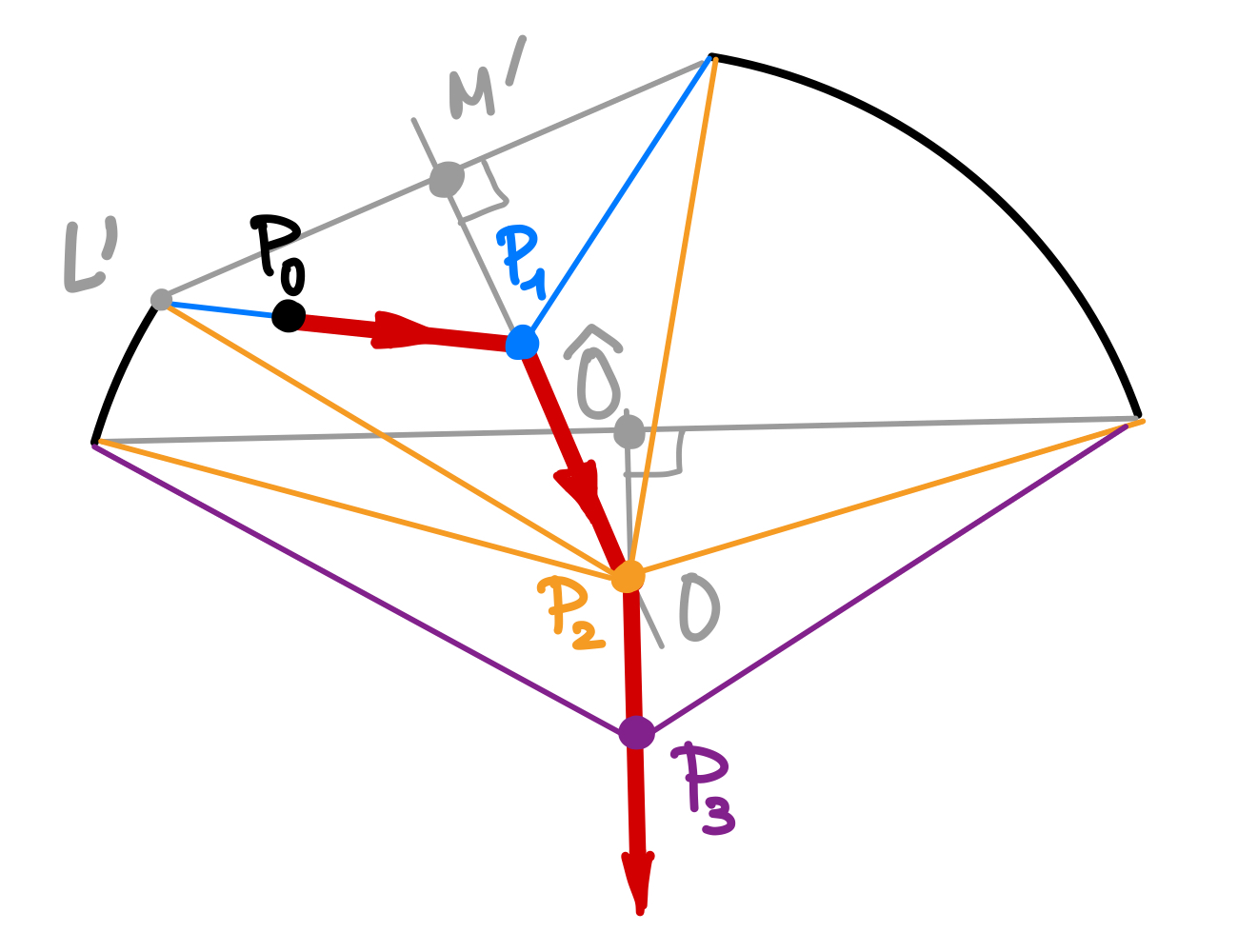}
		\captionof{figure}{Matter concentrates} \label{fig-matcon}
\end{figure}

\subsection*{An illustration of \MAG's non-locality}\label{sec-non-local-2}

It was told at page \pageref{sec-nonlocal} that \MAG\ is a non-local theory. To illustrate this statement, let us go 
 back to Figure \ref{fig-y1} and imagine that point $M$ takes a small step west to become $M'$ (say). This means that only the first coordinate of this $2$-mapping is modified. As $M$ does not stay on the first diagonal, the two particles, whose respective  locations on the line are the coordinates of $M$,  are far apart. We see that, when $M$ moves west to $M'$, the force $\overrightarrow{LM}$ acting on $M$   shrinks a little bit and rotates a little bit in accordance to the transformation of  $\overrightarrow{LM}$  to $\overrightarrow{LM'}.$ As a consequence, the second coordinate of this force is modified (and so is the first one). We see that a small perturbation of the first particle modifies instantaneously  the force acting on the second particle, although these particles are far apart. This is non-locality. 
 
 In addition, if the point $M$ of Figure \ref{fig-y1} were located  much further  south,  the same westward movement would create a weaker force on the second particle. As is the case with Newton's theory,  the   force  is a decreasing function of the distance (at least for large distances).

 \section{Quantum potential}\label{app-e}

Let us consider an abstract setting where the  configuration  space is $\Rn$, the measure
 \[
 m(dx):=m(x)\, dx
 \]
on $\Rn$  is the equilibrium measure of the Markov process with generator
\begin{align*}
A^m:= \nabla \log \sqrt{m} \scal \nabla+   \Delta/2,
\end{align*}
where the  function  $m:\Rn\to(0, \infty)$ is positive everywhere and differentiable. 
 \subsubsection*{A reversible path measure: $R^m$}
 It is assumed that there is a unique path measure $R^m$ that solves the martingale problem with generator $A^m$ and initial measure $m.$  Denoting $(X_t) _{ t_0\le t\le t_1}$ the canonical process, this is equivalent to
\begin{align}\label{eq-89b}
\begin{split}
&dX_t=  \nabla \log \sqrt{m} (X_t)\, dt +  dB ^{ R^m}_t,\quad t_0\le t\le t_1,\quad R^m\ae,\\
&X _{ t_0}\sim m,
\end{split}
\end{align}
where $B ^{ R^m}$ is an $R^m$-Brownian motion. Not only $R^m$ is $m$-stationary, but also  it is reversible.

 \begin{definition}[Quantum potential with respect to $m$]
 For any regular enough probability measure $p$, the quantum potential of $p$ with respect to $m$ is defined by
\begin{multline}\label{eq-90}
\begin{split}
\mathcal{Q}(p|m):= - \frac{A^m \sqrt{\rho}}{ \sqrt{\rho}} 
	&=- \nabla\log \sqrt{m}\cdot \nabla \log \sqrt{\rho}- \frac{ \Delta \sqrt{\rho}}{2 \sqrt{\rho}}\\
	&=- (\nabla\log \sqrt{m}\cdot \nabla \log \sqrt{\rho}
	+ \uud \Delta\log \sqrt{\rho}+\uud|\nabla \log \sqrt{\rho}|^2),
	\end{split}
\end{multline}
where we set $ \rho:=dp/dm.$
 \end{definition}
 
 This notion of \emph{relative} quantum potential seems to be new.  The usual one corresponds to $m=\Leb.$
 
\begin{lemma}\label{res-x4}
The Fisher information of $p$ with respect to $m$   satisfies
\begin{align*}
I(p|m)=\int _{\Rn} \mathcal{Q}(p|m)\,dp,
\end{align*}
meaning that it is the average of the quantum potential,  and
\begin{align*}
I(p|m)=\ud\int _{ \Rn}  | \nabla\sqrt{ dp/dm}|^2\, dm
= \ud \int _{ \Rn}  |\nabla \log \sqrt{ dp/dm}|^2\, dp.
\end{align*}
\end{lemma}

\begin{proof}
Taking $I(p|m):=\int _{\Rn} \mathcal{Q}(p|m)\,dp$ as a definition,  $I(p|m)=\ud\int _{ \Rn}  | \nabla\sqrt{ dp/dm}|^2\, dm$ follows from the  integration by parts formula
\begin{align}\label{eq-87}
\int _{ \Rn} \nabla u\cdot \nabla v\, dm=\int _{ \Rn} \Gamma^m(u,v)\, dm=-2\int _{ \Rn} uA^mv\, dm
\end{align}
which holds because $R^m $ is $m$-reversible. The carré du champ of its generator $A^m$ is denoted by $ \Gamma^m(uv):=A^m(uv)-uA^mv-vA^mu$.
\end{proof}

\subsubsection*{From $m$ to $\Leb$}

An analogy between the thermal evolution of the entropic interpolation and some quantum evolution is put forward at  Section \ref{sec-B-k-fluid}. As Schrödinger  equation refers to densities with respect to Lebesgue measure, it is worth switching from $R^m$ to the new reversible path measure $R ^{ \Leb}$ which is defined as   the unique solution of the  martingale problem \eqref{eq-89b}
with $m=\Leb,$ that is
\begin{align*}
&dX_t=    dB ^{ R^\Leb }_t,\quad t_0\le t\le t_1,\quad R^\Leb \ae,\\
&X _{ t_0}\sim \Leb,
\end{align*}
where $B ^{ R^\Leb }$ is an $R^\Leb $-Brownian motion. Its Markov generator is \[A^\Leb =  \Delta/2.\]
The corresponding quantum potential   is
\begin{align}\label{eq-99}
\mathcal{Q}(m|\Leb):= -   \frac{A^\Leb  \sqrt{m}}{ \sqrt{m}} =- \frac{  \Delta \sqrt{m}}{ 2\sqrt{m}}
	=- \ud( \Delta\log \sqrt{m}+|\nabla \log \sqrt{m}|^2).
\end{align}
  We already encountered it  at  \eqref{eq-103}.

\begin{lemma}\label{res-x95}
For any  measure $m$  with a positive density and any probability measure $p$ also with a positive density 
\begin{align}
\label{eq-95}
\mathcal{Q}(p|m) &= \mathcal{Q}(p|\Leb)
 - \mathcal{Q}(m|\Leb),\\
\label{eq-95b}
I(p|m)&= I(p|\Leb)- \int  \mathcal{Q}(m|\Leb)\,dp.
\end{align}
\end{lemma}

\begin{proof}["Proof"] We skip the regularity problems.\\
With $ \Delta(uv)=u \Delta v+v \Delta u+2 \nabla u\scal\nabla v$ and $p= \rho m$ we obtain
\begin{align*}
-2 \mathcal{Q}(p|\Leb)
	&= \frac{ \Delta \sqrt{p}}{ \sqrt{p}}
	=\frac{ \Delta \sqrt{ \rho m}}{ \sqrt{ \rho m}}
	= \frac{ \sqrt{ \rho} \Delta \sqrt{m}+ \sqrt{m} \Delta \sqrt{ \rho}+2 \nabla \sqrt{ \rho}\cdot\nabla \sqrt{m}}{ \sqrt{ \rho} \sqrt{m}}\\
	&= \frac{ \Delta \sqrt{m}}{ \sqrt{ m}}+ \frac{ \Delta \sqrt{ \rho}}{ \sqrt{ \rho}}+2  \nabla\log \sqrt{m}\cdot \nabla\log \sqrt{ \rho} 
	\overset{ \eqref{eq-90}}=-2 \mathcal{Q}(m|\Leb)-2 \mathcal{Q}(p|m),
\end{align*}
as announced. Finally, integrate \eqref{eq-95} with respect to $p$ and apply Lemma \ref{res-x4} to get \eqref{eq-95b}.
\end{proof}


\end{document}